\documentclass[11pt]{article}

\usepackage{amsfonts}
\usepackage{graphicx}

\newcommand{\ba}{\begin{array}}
\newcommand{\ea}{\end{array}}

\newtheorem{corol}{Corollary}
\newcommand{\bcr}{\begin{corol}}
\newcommand{\ecr}{\end{corol}}
\newtheorem{lemma}{Lemma}
\newtheorem{hypo}{Assumption}
\newcommand{\bhypo}{\begin{hypo}}
\newcommand{\ehypo}{\end{hypo}}
\newtheorem{defi}{Definition}
\newcommand{\ble}{\begin{lemma}}
\newcommand{\ele}{\end{lemma}}
\newcommand{\bde}{\begin{defi}}
\newcommand{\ede}{\end{defi}}

\newtheorem{prop}{Proposition}
\newcommand{\epr}{\end{prop}}
\newcommand{\bpr}{\begin{prop}}
\newtheorem{teo}{Theorem}
\newcommand{\bth}{\begin{teo}}
\newcommand{\eth}{\end{teo}}
\newtheorem{rema}{Remark}
\newcommand{\bre}{\begin{rema} \rm}
\newcommand{\ere}{ \end{rema}}

\newcommand{\ee}{\end{equation}}
\newcommand{\be}{\begin{equation}}

\def\vY{\vec{Y}}
\def\vPsi{\vec{\Psi}}
\def\vpsi{\vec{\psi}}

\def\lp{\lambda^\prime}

\def\sqs#1{\overline{\underline{\left|\matrix{\cr~~ #1~~\cr\cr}\right|}}}

\begin{document}

\centerline{ \Large \bf On Stokes Matrices  in terms of Connection Coefficients}

\vskip 0.5 cm 
\centerline{\Large Davide Guzzetti}

\vskip 0.3 cm 
\centerline{SISSA, International School for Advanced Studies, Trieste, Italy.}
\vskip 1 cm
\noindent
{\bf Abstract:} 
The classical problem of computing a complete system of Stokes multipliers of a linear system of ODEs of rank one   in terms of some connection coefficients of an  associated hypergeometric  system of ODEs, is solved with no genericness assumptions on the residue matrix at zero, by an extension of the method of \cite{BJL}.  
\vskip 1 cm

\section{Introduction}

\label{back}

In the well known paper  \cite{BJL}, among other results, the authors compute a complete system of Stokes multipliers of a linear system of ODEs of  rank one  (at infinity) in terms of some connection coefficients of an associated Fuchsian (or hypergeometric) system.  In \cite{BJL}, this is done  under some assumptions on the system of rank one. One of them is that the leading term at infinity (matrix $A_0$ below) is diagonalizable with distinct eigenvalues. The second assumption, called     {\it assumption (i)}, is  that the diagonal entries of the residue matrix  at zero (matrix  $A_1$  below)  {\it are not integers}. 

In this paper, we are interested in extending the result above when no assumptions on $A_1$ are made. Moreover, we would like to do this using (an extension of) the  method of \cite{BJL}, since it also allows to  obtain results on solutions and monodromy of the associated Fuchsian system.   To this end, we consider the systems of rank one (\ref{01}) below,  and the associated Fuchsian   system (\ref{02}). We are motivated by the fact that these systems appear in some applications, such as the analytic theory of semisimple Frobenius Manifolds \cite{Dub1}, \cite{Dub2}, \cite{Dub3} and  the isomonodromic approach to Painlev\'e equations \cite{MM}. In  these applications, while $A_0$ is diagonalizable with distinct eigenvalues\footnote{After this paper was completed, the work \cite{GGI} appeared on arXiv (April 2014), showing that for the Frobenius manifold given by  the Quantum Cohomology of Grassmannians,  there may be cases (depending on the dimension) when $A_0$ is still diagonalizable, but with some coinciding eigenvalues.}, {\it assumption (i)} fails in  important non generic cases. 

In this paper, we compute  a complete system of Stokes multipliers  in terms of connection coefficients (and we define the connection coefficients) in the  case when no assumptions on $A_1$ are made. Conversely, we express the first monodromy invariants (traces of products of monodromy matrices)  of system (\ref{02}) in terms of Stokes multipliers.  As a side result, the  monodromy of -- and general relations among  -- higher order primitives of vector solutions of (\ref{02})  are obtained.  
 We achieve our results by an extension of the technique of \cite{BJL} when  $A_1$ is any matrix, while $A_0$ is still diagonalizable with distinct eigenvalues\footnote{This allows to consider the normal form (\ref{01}).}. To our knowledge, such extension of the above technique    was not  in the literature. 
 
 As mentioned above, our result applies in particular to  semisimple Frobenius manifold, where $A_1$ has a special form (see {\it Example } below), but still may violate {\it assumption (i)} of \cite{BJL}. For this special  form of $A_1$, the relation between Stokes matrices and connection coefficients was first computed in \cite{Dub3} (and in \cite{Dub1} and \cite{Dub2} when $A_1$ does satisfy {\it assumption (i)}).
 
  From the point of view of the general theory, the  assumption of distinct eigenvalues of  $A_0$ is still  restrictive, however it is enough for  the applications mentioned above. To our knowledge, the case when no assumptions {\it at all} are made on  $A_0$, possibly including a ramified singularity at infinity, and the system of rank one is not is normal form\footnote{If the eigenvalues are not distinct, or $A_0$ is not irreducible, a Birkhoff normal form may not be achieved. See \cite{Balser} for a review}, has not yet been studied. The most general result available is   in \cite{Loday}, where the explicit relation between Stokes-Ramis matrices and connection constants is obtained    for a   general system of rank one   with the only  assumptions of a single level equal to one. The authors of \cite{Loday}  achieve  the result by means of the theory of summation and resurgence. In particular, our Theorem I (Theorem \ref{RelationSC2}) below,  which we obtain by extending the technique of \cite{BJL}, is contained in  the results of Section 4 of \cite{Loday}, which are obtained by the theory of summation and resurgence.

\subsection{Setting} 

We consider a linear system of  rank one in the form  
\be
\label{01}
{dY\over dz}=\left(A_0+{A_1\over z}\right)Y,
\ee
where $A_0$ and $A_1$ are  $n\times n$ matrices. We assume that $A_0$ is diagonalizable, with distinct eigenvalues. Therefore, without loss of generality, we may assume that $A_0$ is already diagonal:  
$$A_0=\hbox{diag}(\lambda_1,...,\lambda_n),~~~\lambda_i\neq \lambda_j \hbox{ for } i\neq j.
$$
Let us denote the diagonal entries of $A_1$ as follows 
$$\hbox{diag}(A_1)=(\lp_1,...,\lp_n) .
$$  
{\it Assumption (i)} of \cite{BJL} is that $\lp_1,...,\lp_n$ are not integers. In this paper we drop the assumption, namely  {\it we allow    any values of $(\lp_1,...,\lp_n)\in\mathbb{C}^n$}. 

 Solutions of (\ref{01}) can be expressed in terms of convergent Laplace-type integrals \cite{Bi}, \cite{Ince}, where the integrands are solutions of the Fuchsian system
\be
\label{02}
(A_0-\lambda){d \Psi\over d\lambda}=  (A_1+I)\Psi,~~~~~I:=\hbox{ identity matrix}
\ee
 Indeed, let $\vPsi(\lambda)$ be a vector valued function and define 
$$
\vY(z)=\int_\gamma e^{\lambda z}\vPsi(\lambda)d\lambda, 
$$
where $\gamma$ is a suitable path. Substituting in (\ref{01}), we obtain
$$
z\int_\gamma \lambda e^{\lambda z}\vPsi(\lambda)d\lambda
=(z A_0+A_1)\int_\gamma e^{\lambda z}\vPsi(\lambda)d\lambda.
$$
This implies that
$$
A_1 \int_\gamma e^{\lambda z}\vPsi(\lambda)d\lambda= \int_\gamma {d(e^{\lambda z}) \over d\lambda} ~(\lambda-A_0)\vPsi(\lambda)d\lambda=
$$
$$
=e^{\lambda z}(\lambda-A_0)\vPsi(\lambda)\Bigl|_\gamma-\int_\gamma
e^{\lambda z} \left[
(\lambda-A_0){d\vPsi(\lambda)\over d\lambda}+\vPsi(\lambda)
\right]d\lambda.
$$
If $\gamma$ is such that $e^{\lambda z}(\lambda-A_0)\vPsi(\lambda)\Bigl|_\gamma=0$, and if the function $\vPsi(\lambda)$ solves (\ref{02}), then $\vY(z)$ solves (\ref{01}).

\begin{figure}
\centerline{\includegraphics[width=0.4\textwidth]{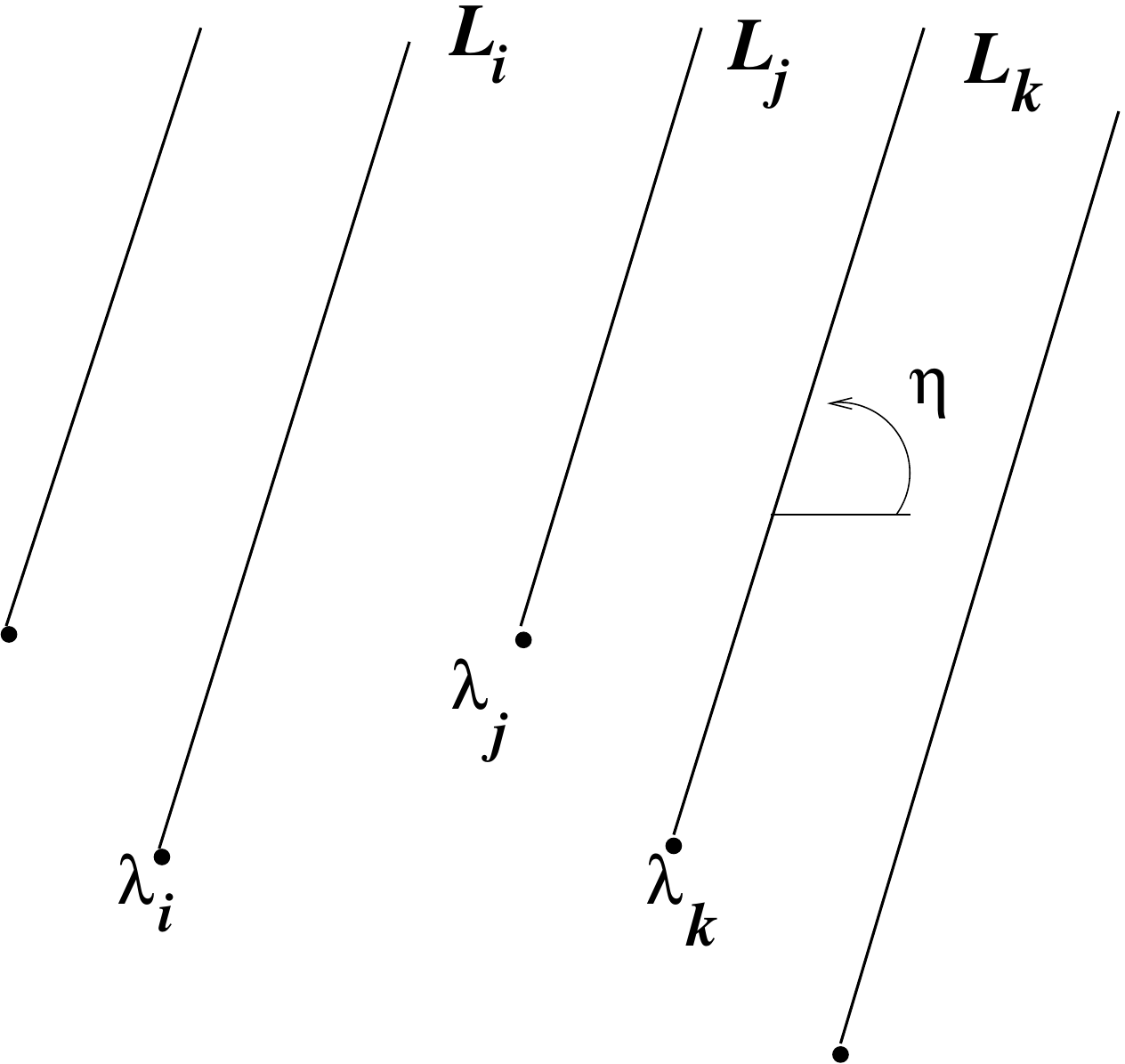}}
\caption{The poles $\lambda_j$, $1\leq j \leq n$ of system (\ref{02}), and branch cuts $L_j$.}
\label{cutss}
\end{figure}

\vskip 0.3 cm 
 In order to generalize the result of  \cite{BJL}, following an analogous method, we need to characterize the solutions of (\ref{02}) without  assumptions on $A_1$.  System (\ref{02}) can be  rewritten as 
\be
\label{03}
 {d\Psi\over d\lambda}=\sum_{k=1}^n {B_k \over \lambda-\lambda_k}\Psi,~~~~~B_k:=-E_k(A_1+I),~~~~1\leq k \leq n,
\ee
where $E_k$ is a $n\times n $ matrix with entries $(E_k)_{kk}=1$ and $(E_k)_{ij}=0$ otherwise.  
A fundamental matrix solution of (\ref{03}) is multivalued in $ \mathbb{C}\backslash\{\lambda_1,...,\lambda_n\}$. 
Let ${\cal U}$ be the universal covering of $ \mathbb{C}\backslash\{\lambda_1,...,\lambda_n\}$. Following \cite{BJL}, we fix parallel cuts $L_k$,  {\it oriented} from $\lambda_k$ to $\infty$ 
$$
L_k:=\{\lambda\in {\cal U}~|~\arg(\lambda-\lambda_k)=\eta\},~~~~~~~~1\leq k \leq n,
$$
 where 
$$
\eta\in \mathbb{R},~~~~~\eta\neq \arg(\lambda_j-\lambda_k) \hbox{ mod } 2\pi,~~~\hbox{ for all }1\leq j,k\leq n.
$$
The above condition means that a cut $L_k$ does not contain another pole $\lambda_j$, $j\neq k$. See figure \ref{cutss}. Such  values of $\eta$ are called {\it admissible}.  Wee fix the branch  $\ln(\lambda-\lambda_k)=\ln|\lambda-\lambda_k|+i\eta-0$ when $\arg(\lambda-\lambda_k)=\eta-0$. The complex plane (as a sheet of ${\cal U}$) with these cuts and choices of the branches of logarithm is denoted 
$$
{\cal P}_\eta:=\left\{ \lambda\in {\cal U}~|~\eta-2\pi< \arg(\lambda-\lambda_k)<\eta ,~~1\leq k\leq n \right\}
.
$$

We prove in Section \ref{locall} that system (\ref{02}), depending on the values of $(\lp_1,...,\lp_n)\in\mathbb{C}^n$, admits a matrix solution (not necessarily fundamental) of the form:
$$
\Psi(\lambda)=\Bigr[\vPsi_1(\lambda)~|~\cdots~|~\vPsi_n(\lambda)\Bigr],~~~~~\lambda\in{\cal P}_\eta
$$
whose columns $\vPsi_k(\lambda)$, $k=1,...,n$, have the following  behaviours  in a neighbourhood of $\lambda_k$:
$$
    \vPsi_k(\lambda)=
\left\{
\matrix{\left(\Gamma(\lp_k+1)\vec{e}_k+\sum_{l\geq 1} \vec{b}_l^{~(k)}(\lambda-\lambda_k)^l\right)~(\lambda-\lambda_k)^{-\lp_k-1} & ~~~\lp_k\not\in\mathbb{ Z}
\cr
\cr
\cr
\left({(-1)^{N_k}\over (-N_k-1)!}\vec{e}_k+\sum_{l\geq 1} \vec{b}_l^{~(k)}(\lambda-\lambda_k)^l\right)(\lambda-\lambda_k)^{-N_k-1}& ~~~\lp_k=N_k\in \mathbb{Z}_{-}
\cr\cr
\cr
\vec{d}_0^{~(k)}+\sum_{l\geq 1} \vec{d}_l^{~(k)}(\lambda-\lambda_k)^l& ~~~\lp_k\in
\mathbb{ N}
}
\right.
$$

$$\mathbb{N}=\{0,1,2,...\} ~\hbox{ integers},~~~~~~~~\mathbb{Z}_{-}=
\{-1,-2,-3,...\}~\hbox{ negative integers}, 
$$
$$
 \vec{e}_k = \hbox{ $k$-th unit column vector in $\mathbb{C}^n$}.
$$
The  Taylor series in $(\lambda-\lambda_k)$ converge  in a neighbourhood of $\lambda_k$. The coefficients $ \vec{b}_l^{~(k)}\in\mathbb{C}^n$ are uniquely determined by the choice of the normalizations $\Gamma(\lp_k+1)\vec{e}_k$ and ${(-1)^{N_k}\over (-N_k-1)!}\vec{e}_k$. The coefficients  $ \vec{d}_l^{~(k)}\in\mathbb{C}^n$ are uniquely determined by the existence of a singular vector  solution at $\lambda_k$ with behaviour 
$$
\vPsi_k(\lambda)\ln(\lambda-\lambda_k)+{N_k!\vec{e}_k +O(\lambda-\lambda_k) \over (\lambda-\lambda_k)^{N_k+1}},~~~~~~N_k=\lp_k\in \mathbb{N}.
$$
 
We will show (Definition \ref{ofconnection}, Section \ref{locall}) that there exist unique {\it connection coefficients } $c_{jk}\in\mathbb{C}$  such that, in a neighbourhood of any $\lambda_j\neq \lambda_k$:  

\be 
\label{psik}
 \vPsi_k(\lambda) =
\left\{
\matrix{
\vPsi_j(\lambda)c_{jk}+\hbox{reg}(\lambda-\lambda_j),& \lp_j\not\in \mathbb{ Z}
\cr
\cr
 \vPsi_j(\lambda)\ln(\lambda-\lambda_j)c_{jk}+\hbox{reg}(\lambda-\lambda_j),& \lp_j\in \mathbb{Z}_{-}
\cr
\cr
 \left(\vPsi_j(\lambda)\ln(\lambda-\lambda_j)+{P_{N_j}^{(j)}(\lambda)\over (\lambda-\lambda_j)^{N_j+1}}\right)c_{jk}+\hbox{reg}(\lambda-\lambda_j),& \lp_j=N_j\in \mathbb{ N}
}
\right.
\ee
Here $P_{N_j}^{(j)}$ is a polynomial in $(\lambda-\lambda_j)$ of degree $N_j$, and $\hbox{reg}(\lambda-\lambda_j)$ is a vector function analytic (regular) in a neighbourhood of $\lambda_j$.   
 We will characterize the connection coefficients and the solutions $\vPsi_k$ in  Section \ref{locall}. In particular, 
 $$
 c_{kk}=1 ~~\hbox{ if } \lp_k\not \in \mathbb{Z},~~~~~c_{kk}=0~~ \hbox{ if } \lp_k\in \mathbb{Z}.
 $$

Let  $\tau = 3\pi/2-\eta$. There are three unique fundamental matrices of $(\ref{01})$, say $Y_I(z)$, $Y_{II}(z)$ and $Y_{III}(z)$, with canonical asymptotic behaviour $(I+O(1/z))\exp\{A_0z+A_1\ln z\}$ in the three sectors $\{z~|~\tau-\pi\leq \arg z\leq \tau\}$,  $\{z~|~\tau \leq \arg z\leq \tau+\pi\}$ and  $\{z~|~\tau+\pi \leq \arg z\leq \tau+2\pi\}$ respectively. They are related by two Stokes matrices $S_+$ and $S_-$   such that   
 $$Y_{II}(z)=Y_I(z)S_{+},~~\arg z=\tau;~~~~~~~~Y_{III}(z)=Y_{II}(z)S_{-},~~\arg z=\tau+\pi.
 $$
 Introduce in $\{1,2,...,n\}$ the partial ordering $\prec$ given by 
 $$
 j\prec k ~~\Longleftrightarrow~~ \Re(z(\lambda_j-\lambda_k))<0 ~\hbox{ for } \arg z=\tau,~~~~~i\neq j,~~i,j\in\{1,...,n\}.
 $$

 \subsection{Main Results}

 \vskip 0.3 cm 
 \noindent
 {\bf Proposition I (Proposition \ref{promod}):} {\it Let the branch cuts $L_1,...,L_n$ be fixed.  The monodromy matrix  $M_k=(m_{ij}^{(k)})_{i,j=1...n}$  of $\Psi(\lambda)$ representing a small loop in  anticlockwise direction around $\lambda_k$, not encircling all the other points $\lambda_j\neq \lambda_k$, $j=1,...,n$ is:
  $$
  m_{jj}^{(k)}=1~~~ 1\leq j\leq n, ~~j\neq k;~~~~~~~ m_{kk}^{(k)}=e^{-2\pi i\lp_k};
  $$ 
 $$
m_{kj}^{(k)}=\alpha_k ~c_{kj},~~~1\leq j \leq n, ~~~j\neq k;~~~~~~m_{ij}^{(k)}=0~\hbox{ otherwise}.
$$
where 
$$
\left\{
\matrix{
\alpha_k:=(e^{-2\pi i \lp_k}-1), & \hbox{ if }\lp_k\not\in \mathbb{ Z},
\cr
\cr
\alpha_k:=2\pi i ,& \hbox{ if }\lp_k\in\mathbb{Z}.
}
\right.        
$$
Equivalently, the effect of the loop on $\Psi(\lambda)$ is
$$
\vPsi_k(\lambda)\longmapsto e^{-2\pi i \lp_k} \vPsi_k(\lambda);~~~~~\vPsi_j(\lambda)\longmapsto \vPsi_j(\lambda)+\alpha_k c_{kj} \vPsi_k(\lambda),~~~j\neq k.
$$
}

\vskip 0.3 cm 
\noindent
{\bf Theorem I (Theorem \ref{RelationSC2}):} {\it 
The Stokes  matrices of system (\ref{01}) are given in terms of the connection coefficients $c_{jk}$ of system (\ref{02}) according to the formulae
$$
[ S_+]_{jk}
=
\left\{
\matrix{
e^{2\pi i \lp_k}\alpha_k ~c_{jk}& ~~~\hbox{ for } j\prec k,
\cr
\cr
             1   & ~~~\hbox{ for } j =k,
\cr
\cr
             0   &  ~~~\hbox{ for } j\succ k,
}
\right.
~~~~~~~~
[ S_{-}^{-1}]_{jk}
=
\left\{
\matrix{
             0   &  ~~~\hbox{ for } j\prec k,
\cr
\cr
             1   &~~~ \hbox{ for } j =k,
\cr
\cr
        -e^{2\pi i (\lp_k-\lp_j)}\alpha_k~c_{jk} & 
                                         ~~~\hbox{ for } j\succ k.
}
\right.
$$

}

\vskip 0.3 cm 
\noindent
{\bf Corollary I (Corollary \ref{CHICH}):} {\it  
The following equalities hold for the monodromy matrices of $\Psi(\lambda)$:
$$
 \hbox{\rm Tr}(M_k)=n-1+e^{-2\pi i \lp_k}
$$

$$
 \hbox{\rm Tr}(M_jM_k)=
\left\{
\matrix{
n-2 +e^{-2\pi i \lp_j}+e^{-2\pi i \lp_k}-e^{-2\pi i \lp_j}~[S_+]_{jk}[S_{-}^{-1}]_{kj} & \hbox{ if } j\prec k,
\cr
\cr
n-2 +e^{-2\pi i \lp_j}+e^{-2\pi i \lp_k}-e^{-2\pi i \lp_k}~[S_{-}^{-1}]_{jk}[S_+]_{kj} & \hbox{ if } j\succ k.
}
\right.
$$

}

\vskip 0.3 cm 
 \noindent
 {\bf Proposition II (Propositions \ref{propoC} and \ref{PropoE}):} 
{\it If $A_1$ has no integer eigenvalues, then $\Psi(\lambda)$ is a fundamental matrix and $M_1,...,M_n$  generate the monodromy group of system (\ref{02}). Moreover,  the matrix $C:=(c_{jk})$ is invertible if and only if $A_1$ has no integer eigenvalues. 
}

\vskip 0.3 cm 
\noindent
{\bf Remark:} There are cases when $A_1$ has integer eigenvalues and $\Psi$ is fundamental. We prove that in these cases, necessarily, some $\lp_k\in \mathbb{Z}$.

\vskip 0.5 cm 
\noindent
{\it Example:} 
When system (\ref{01}) is associated to Frobenius Manifolds \cite{Dub1}, \cite{Dub2}, \cite{Dub3}, the matrix $A_1$ has a special form, namely it is expressed in terms of a skew symmetric matrix $V$ as follows:  
$$
A_1= V-\left({1\over 2}+\nu\right)I,~~~~~\nu\in \mathbb{C},~~~~~V^T=-V
$$  
 We show how our general results above apply to this case. Since $\lp_k= -\nu-{1\over 2}$, $1\leq k \leq n$, 
 it follows that 
 $$
  \alpha_1=\alpha_2=\cdots= \alpha_n =\alpha, ~~~\hbox{ where }~~\alpha:=\left\{
  \matrix{
  -(1+e^{2\pi i \nu}) & \hbox{ if } \nu\not\in \mathbb{Z}+{1\over 2},
  \cr
  2\pi i &  \hbox{ if } \nu \in \mathbb{Z}+{1\over 2}.
  }
  \right.
  $$
  From Theorem I above (and the fact that the $c_{kk}=0$ when $\lp_k\in \mathbb{Z}$), we deduce that  
  $$
e^{2\pi i \nu} S_{+}+S_{-}^{-1}= -\alpha C,~~~~~\hbox{ where }~C:=(c_{ij}).
$$  
  Since $V$ is  a $n\times n$ skew symmetric matrix, it can be easily verified that 
 $$
 S_+^T=S_-^{-1}.
 $$
 Thus  
  \be 
  \label{dubroveq}
e^{2\pi i \nu} S_{+}+S_{+}^T= -\alpha C. 
\ee
  The above, and Proposition II, 
  allow us to conclude that if $e^{2\pi i \nu} S_{+}+S_{+}^T$ is invertible, then $A_1$ has no integer eigenvalues and so $\Psi(\lambda)$ is invertible. 
  This is part of the first assertion of Theorem 4.3  of \cite{Dub3}, namely  if 
  $$
  \det(e^{2\pi i \nu} S_{+}+S_{+}^T)\neq 0,
  $$
   then system (\ref{02}) has $n$ linearly independent solutions $\vPsi_1,...,\vPsi_n$. From (\ref{dubroveq}) and Proposition I, it follows that for an anticlockwise loop around $\lambda_i$, the monodromy of the above solutions is 
$$
 \vPsi_i\longmapsto - e^{2\pi i \nu} \vPsi_i,~~~~~~~~\vPsi_j\longmapsto \vPsi_j -e^{i\pi \nu}
 \left( e^{i\pi \nu} S_++e^{-i\pi \nu} S_+^T   \right)_{ij} \vPsi_i,~~~j\neq i.
 $$  
  The above is formula (4.11) in Theorem 4.3  of \cite{Dub3}. $\Box$
  
  \vskip 0.2 cm

The paper is organized as follows:

-- Section \ref{locall}: We construct vector solutions $\vPsi_k(\lambda)$, $1\leq k \leq n$, to system (\ref{02})-(\ref{03}), and define the connection coefficient, with no assumptions on $A_1$. 

-- Section \ref{MV}:  We construct two matrix solutions  $\Psi$ and $\Psi^*$ to system (\ref{02})-(\ref{03}), discuss when they are fundamental, and compute their monodromy in terms of connection coefficients (with no assumptions on $A_1$).  

-- Section \ref{RC}: We discuss the dependence of $\Psi$ and $\Psi^*$ on the choice of the branch cuts $L_1$, ..., $L_n$ (with no assumptions on $A_1$). 

-- Section \ref{FL}: We define a complete set of Stokes multipliers for (\ref{01}). We write the columns of the fundamental matrix of system (\ref{01}), having canonical asymptotics in a wide sector,  as Laplace integrals of  the $\vPsi_k$, $1\leq k \leq n$, and express the latter in terms of the the coefficients of the former asymptotics.  

-- Section \ref{MT}: We state the main theorem (Theorem \ref{RelationSC2}), which gives  Stokes matrices and Stokes factors of (\ref{01}) in terms of connection coefficients of (\ref{02})-(\ref{03}), and  express the first monodromy invariants of system (\ref{02})-(\ref{03}) in terms of Stokes matrices (Corollary \ref{CHICH}).

-- Section \ref{THEPROOF}: we prove Theorem \ref{RelationSC2}, and find relations and monodromy for  $q$-primitives of vector solutions of (\ref{02})-(\ref{03}).  

-- In the Appendix, we prove some propositions  which generalize similar results of \cite{BJL}  when no assumptions on $A_1$ are made. 

\section{Local Solutions of System  (\ref{03}) (equivalently, of (\ref{02}))}
\label{locall}

The matrix $B_k$ in system (\ref{03}) has zero entries, except for the $k$-th row. Indeed,  letting $A_1=(A_{ij})_{i,j=1,...,n}$, a straightforward computation yields 
$$
B_k=\pmatrix{0   &&&0&&&0\cr
                        \vdots &&&\vdots&&&\vdots\cr
 -A_{k1}&\cdots& -A_{k,k-1}&-\lp_k-1&-A_{k,k+1}& \cdots&-A_{kn}\cr
\vdots&&&\vdots&&&\vdots\cr
0&&&0&&&0
}
$$
A fundamental matrix solution of (\ref{03}) is multivalued in $ \mathbb{C}\backslash\{\lambda_1,...,\lambda_n\}$ and single-valued in ${\cal P}_\eta$, for and admissible direction $\eta$.
 If $\lambda$ is in a neighbourhood of a $\lambda_k$ not containing the other poles, there exists a fundamental matrix solution
$$
\Psi^{(k)}(\lambda)=[\vPsi_1^{(k)}(\lambda)~|~\cdots~|~\vPsi_n^{(k)}(\lambda)],
$$
which can be computed in a standard way, depending on the value of $\lp_k$  (see  \cite{Wasow}). In \cite{BJL}, only the case $\lp_k\not \in \mathbb{Z}$ is considered (point 1) below). Here we  need to analyse also the case $\lp_k\in \mathbb{Z}$ (points 2), 3) and 4) below).  

\vskip 0.3 cm 
{\large \bf 1) [Generic Case, as in \cite{BJL}].} If $\lp_k\not\in \mathbb{Z}$, then $B_k$ is diagonalizable, with diagonal form
$$
      T^{(k)}=[G^{(k)}]^{-1}B_kG^{(k)}=\hbox{diag}(0,...,0,-\lp_k-1,0,...,0),
$$
 where the non zero entry is at the $k$-th position. 
The $k$-th column of the diagonalizing matrix $G^{(k)}$ can be chosen to be a multiple of the  $k$-th vector $\vec{e}_k$ of the canonical basis of $\mathbb{C}^n$. As in \cite{BJL} we choose normalization
$ \Gamma(\lp_k+1)\vec{e}_k$. Any other column of $G^{(k)}$ has two non zero entries. A fundamental matrix is then
$$
 \Psi^{(k)}(\lambda)= G^{(k)}(I+O(\lambda-\lambda_k))(\lambda-\lambda_k)^{T^{(k)}},
$$
 Here $O(\lambda-\lambda_k)$ is a matrix valued Taylor series, converging in the neighbourhood of $\lambda_k$ and vanishing as $\lambda \to \lambda_k$. 
Write 
$$G^{(k)}(I+O(\lambda-\lambda_k))=[\vpsi_1^{(k)}(\lambda)|~...~|\vpsi_n^{(k)}(\lambda)],
$$
 where the columns $\vpsi_j^{(k)}(\lambda)$ are analytic functions in a neighbourhood of $\lambda_k$, expanded in convergent Taylor series. 
Then:
$$
 \Psi^{(k)}(\lambda)=\Bigl[\vpsi_1^{(k)}(\lambda)~|~...~|~\vpsi_{k-1}^{(k)}(\lambda)~|
~\vpsi_k^{(k)}(\lambda)(\lambda-\lambda_k)^{-\lp_k-1}~|~
\vpsi_{k+1}^{(k)}(\lambda)~|~...~|~\vpsi_{n}^{(k)}(\lambda)\Bigr].
$$
The columns are $n$ independent vector solutions,  $n-1$ being analytic and the $k$-th  singular.  We assign the   symbol $\vPsi_k$ to the singular solution, as follows
\be
\label{asso1}
  \sqs{
\vPsi_k(\lambda):=\vPsi_k^{(k)}(\lambda)=\vpsi_k^{(k)}(\lambda)(\lambda-\lambda_k)^{-\lp_k-1}.
}
\ee
where 
$$
\sqs{
\vpsi_k^{(k)}(\lambda)=\Gamma(\lp_k+1){\vec{e}_k}+\sum_{l\geq 1} \vec{b}_l^{~(k)}(\lambda-\lambda_k)^l.
}
$$
 The vector coefficients $\vec{b}_l^{(k)}$ can be computed rationally from the matrix coefficients $B_l$'s of system (\ref{03}). See \cite{Wasow}. The above $\vPsi_k$ is called {\it associated function} in \cite{BJL}.

\vskip 0.3 cm 
{\large \bf 2) [Jordan Case].} If $\lp_k=-1$, then $B_k$ has Jordan form 
$$
J^{(k)}=[G^{(k)}]^{-1}B_kG^{(k)}=\pmatrix{0 & & & & & \cr
&\ddots& &     &           &\cr 
&          & 0&1&            &\cr
&          & 0&0&             &\cr
&          &   &   & \ddots&\cr
&           &   &   &          &0\cr 
}.
$$
Entry 1 is at row $(k-1)$ and column $k$. The column $k-1$ of $G^{(k)}$ can be normalized to be $-\vec{e}_k$. The $k$-th column only has a non zero entry, and the other columns have two non zero entries. There exist a fundamental matrix solution with local representation 
$$
 \Psi^{(k)}(\lambda)= G^{(k)}(I+O(\lambda-\lambda_k))(\lambda-\lambda_k)^{J^{(k)}}~~~~~~~~~~~~~~
$$
$$
=\Bigl[\vpsi_1^{(k)}(\lambda)|~...~|\vpsi_n^{(k)}(\lambda)\Bigr](\lambda-\lambda_k)^{J^{(k)}}
$$
$$
=\Bigl[\vpsi_1^{(k)}(\lambda)~|~...~|~\vpsi_{k-1}^{(k)}(\lambda)~|
~\vpsi_{k-1}^{(k)}(\lambda)\ln(\lambda-\lambda_k)+\vpsi_k^{(k)}(\lambda)~|~
\vpsi_{k+1}^{(k)}(\lambda)~|~...~|~\vpsi_{n}^{(k)}(\lambda)\Bigr],
$$
where the columns $\vpsi_j^{(k)}$ are analytic in a neighbourhood of $\lambda_k$. The columns are $n$ independent vector solutions,  $n-1$ being analytic and the $k$-th  singular. We assign the   symbol $\vPsi_k$ to the non-singular factor of $\ln(\lambda-\lambda_k) $,  
as follows 
\be
\label{asso12}
\sqs{
  \vPsi_k(\lambda):= \vpsi_{k-1}^{(k)}(\lambda)=-\vec{e}_k+\sum_{l\geq 1} \vec{b}_l^{~(k)}(\lambda-\lambda_k)^l.
}
\ee
Note that this is a solution of (\ref{03}). Then, the $k$-th column of $\Psi^{(k)}$ is 
\be
\label{asso2}
\sqs{
\vPsi^{(k)}_k(\lambda)=\vPsi_k(\lambda)\ln(\lambda-\lambda_k)+\hbox{reg}(\lambda-\lambda_k),
}
\ee
where $\hbox{reg}(\lambda-\lambda_k)$ means an analytic (vector) function in a neighbourhood of $\lambda_k$.

\vskip 0.3 cm
{\large \bf 3) [First Resonant Case]} If $\lp_k=N_k\geq 0$ is integer, then $B_k$ is diagonalizable as in case 1), but now a fundamental solution has the form 
$$
 \Psi^{(k)}(\lambda)= G^{(k)}(I+O(\lambda-\lambda_k))~(\lambda-\lambda_k)^{T^{(k)}}~(\lambda-\lambda_k)^{R^{(k)}},
$$
where $R^{(k)}$ is a matrix with zero entries expect for $R^{(k)}_{jk}$,  $j=1,..,n$, and  $j\neq k$, because $(Eigen(B_k))_j-Eigen(B_k))_k=N_k+1>0$. Thus, only the $k$-th column of $R^{(k)}$ may be non zero. Let $r_j^{(k)}:=R^{(k)}_{jk}$, so that the $k$-th column is
 $$
 \vec{r}^{~(k)}=(r^{(k)}_1,...,r^{(k)}_{k-1},0,r^{(k)}_{k+1},...,r^{(k)}_n)^T,
 $$
  where $T$ means transposition. The entries $r^{(k)}_j$ are computed as rational functions of the entries of the  matrices $B_l$, $l=1,..,n$ (see \cite{Wasow}). From the above, it follows that  
$$
\Psi^{(k)}(\lambda)=\Bigl[\vpsi_1^{(k)}(\lambda)|~...~|\vpsi_n^{(k)}(\lambda)\Bigr]~(\lambda-\lambda_k)^{T^{(k)}}~(I+R^{(k)}\ln(\lambda-\lambda_k))=
$$
$$
=\Bigl[\vpsi_1^{(k)}(\lambda)~|~...~|~\vpsi_{k-1}^{(k)}(\lambda)~|
~\Psi_k^{(k)}(\lambda)~|~\vpsi_{k+1}^{(k)}(\lambda)
~|~...~|~\vpsi_{n}^{(k)}(\lambda)\Bigr],
$$
where
$$\Psi_k^{(k)}(\lambda)=\Bigl\{\sum_{j\neq k}r^{(k)}_j\vpsi_{j}^{(k)}(\lambda)\Bigr\}\ln(\lambda-\lambda_k)+{\vpsi_k^{(k)}(\lambda)\over (\lambda-\lambda_k)^{N_k+1} },
$$
$$
\vpsi^{(k)}_k(\lambda)=N_k!~\vec{e}_k+O(\lambda-\lambda_k),
$$
the factor $N_k!$ coming from a chosen normalization of $G^{(k)}$. 
The columns are $n$ independent vector solutions,  $n-1$ being analytic (i.e. the $\vpsi_j^{~(k)}$, $j\neq k$) and the $k$-th  singular (i.e. $\Psi_k^{(k)}$) .   We assign the  symbol $\vPsi_k$  to the non-singular factor of $\ln(\lambda-\lambda_k) $ 
as follows
$$
\sqs{
  \vPsi_k(\lambda):=\sum_{j\neq k}r^{(k)}_j\vpsi_{j}^{(k)}(\lambda)=\sum_{l\geq 0} {\vec{d}_l}^{~(k)}(\lambda-\lambda_k)^l.
}
$$
Note that this is a solution of (\ref{03}), being linear combination of regular solutions $\vpsi_{j}^{(k)}$. Special cases can occur when $\vec{r}^{~(k)}=0$, so that $\vPsi_k(\lambda)\equiv 0$.  We conclude that  the $k$-th column of $\Psi^{(k)}$ is 
\be
\sqs{
\label{asso3}
\vPsi^{(k)}_k(\lambda)=\vPsi_k(\lambda)\ln(\lambda-\lambda_k)+{ P_{N_k}^{(k)}(\lambda)
\over (\lambda-\lambda_k)^{N_k+1}}
+\hbox{reg}(\lambda-\lambda_k),
}
\ee
where 
$$
\sqs{
P_{N_k}^{(k)}(\lambda)=N_k!~\vec{e}_k+\sum_{l=0}^{N_k}\vec{b}_l^{~(k)}(\lambda-\lambda_k)^l,
}
$$
represents the first $N_k+1$ terms in the expansion of $\vpsi_k^{(k)}$. The vector coefficients $\vec{b}_l^{(k)}$ are computed rationally from the coefficients of (\ref{03}).  The solution (\ref{asso3}) is not uniquely determined, because we can add a linear combination of regular solutions $\vpsi_{j}^{(k)}$, but the singular part is uniquely determined by the normalization of $P^{(k)}(\lambda)$. Consequently, also  $\vPsi_k(\lambda)$ is  uniquely determined.

\vskip 0.3 cm 
{\large \bf 4) [Second Resonant Case]} If $\lp_k=N_k\leq -2$ is integer, then $B_k$ is diagonalizable as in case 1), but now a fundamental solution has the form 
$$
 \Psi^{(k)}(\lambda)= G^{(k)}(I+O(\lambda-\lambda_k))~(\lambda-\lambda_k)^{T^{(k)}}~(\lambda-\lambda_k)^{R^{(k)}},
$$
where $R^{(k)}$ is a matrix with zero entries expect for $R^{(k)}_{kj}$,  $j=1,..,n$, and $j\neq k$, because $(Eigen(B_k))_k-Eigen(B_k))_j=-N_k-1>0$. Thus, only the $k$-th   row  of $R^{(k)}$ may be non zero. Let $r_j^{(k)}:=R^{(k)}_{kj}$, so that the $k$-th   row is
$$
\underline{r}^{(k)}=[r^{(k)}_1,...,r^{(k)}_{k-1},0,r^{(k)}_{k+1},...,r^{(k)}_n],
$$ 
where the entries $r^{(k)}_j$ are computed as rational functions of the entries if the  matrices $B_l$, $l=1,..,n$ (see \cite{Wasow}).  Thus, 
$$
\Psi^{(k)}(\lambda)=\Bigl[\vpsi_1^{(k)}(\lambda)|~...~|\vpsi_n^{(k)}(\lambda)\Bigr]~(\lambda-\lambda_k)^{T^{(k)}}~~(I+R^{(k)}\ln(\lambda-\lambda_k))
$$
where the $\vpsi_j^{~(k)}(\lambda)$ are analytic and Taylor expanded in a neighbourhood of $\lambda_k$. 
The columns of the above matrix are
$$
 \vPsi_j^{(k)}(\lambda)=r_j^{(k)}\vpsi_k^{(k)}(\lambda)~(\lambda-\lambda_k)^{-N_k-1}\ln(\lambda-\lambda_k) +\vpsi_j^{(k)}(\lambda),~~~~~j=1,...,n,~~j\neq k,
$$
$$
 \vPsi_k^{(k)}(\lambda)=\vpsi_k^{(k)}(\lambda)~(\lambda-\lambda_k)^{-N_k-1}.
$$
There are  at most $n-1$ independent singular solutions at $\lambda_k$,  and at least one analytic  solution $\vPsi_k^{(k)}$. In special cases, it may happen that $\underline{r}^{(k)}=0$, so that there are $n$ independent solutions analytic at $\lambda_k$. We show below (Lemma \ref{onanalytic}) that in fact we can always find $n-1$ independent   solutions  analytic at  $\lambda_k$, whatever $\underline{r}^{(k)}$ is. 

\noindent
We assign the symbol $\vPsi_k$ to the $k^{th}$ column:
$$
\sqs{
\vPsi_k(\lambda):= \vPsi_k^{(k)}(\lambda)=\vpsi_k^{(k)}(\lambda)~(\lambda-\lambda_k)^{-N_k-1}
}
$$
with normalization 
\be 
\label{ORN}
\sqs{\vpsi_k^{(k)}(\lambda)={(-1)^{N_k}\over (-N_k-1)!}\vec{e}_k+\sum_{l\geq 1}b_l^{(k)}(\lambda-\lambda_k)^l
}
\ee
where the convergent Taylor series has coefficients determined rationally by the 
matrices  $B_l$'s of (\ref{03}).  The logarithmic solutions are rewritten as 
 $$
 \vPsi_j^{(k)}(\lambda)=r_j^{(k)}\vPsi_k(\lambda)\ln(\lambda-\lambda_k)+\vec{\psi}_j^{(k)}(\lambda), ~~~~~j\neq k,~~1\leq j\leq n.
 $$  
It follows that if at least one $r_j^{(k)}\neq 0$, we can pick up the singular solutions 
\be
\label{STARE}
\sqs{
\vPsi_k(\lambda)\ln(\lambda-\lambda_k)+\hbox{reg}(\lambda-\lambda_k)
},
\ee
The  regular part is an arbitrary  linear combination of  the $\vpsi_{j}^{(k)}$'s, $1\leq j\leq n$, $j\neq k$. The singular part is determined uniquely by the normalization (\ref{ORN}).

\ble
\label{onanalytic}
 Let $\lp_k$ be an integer $N_k \leq -2$.  If $\underline{r}^{(k)}=0$, system (\ref{03}) has $n$ independent solutions analytic at $\lambda_k$. . If $\underline{r}^{(k)} \neq 0$, system (\ref{03}) has $n$ independent solutions, of which  $n-1$ are  analytic and one is log-singular  at $\lambda_k$. 
\ele

\vskip 0.2 cm 
\noindent
{\it Proof:} Let  $0 \leq s\leq n-1$ be the number of non zero values  $r_{i_1},...,r_{i_s}$. 
 If  $\underline{r}^{(k)}=0$, then $s=0$ and by the preceding construction there exist $n$ independent solutions 
 $$
\vpsi_{1}^{(k)},~...~,\vpsi_{k-1}^{(k)},~\vPsi_k,~\vpsi_{k+1}^{(k)},~...~,\vpsi_{n}^{(k)}.
$$
If $s>0$, then consider the following partition of $\{1,2,...,k-1,k+1,...,n\}$:
$$\{i_1,i_2,...,i_s\}\cap\{j_1,j_2,...,j_l\}=\emptyset,~~~~~l+s=n-1,
$$
$$
 \{i_1,i_2,...,i_s\}\cup\{j_1,j_2,...,j_l\}=\{1,2,...,k-1,k+1,...,n\}.
$$  There are  $s$ singular (at $\lambda_k$) solutions 
$$ 
 \vPsi_{i_1}^{(k)},~ \vPsi_{i_2}^{(k)},~...~,~ \vPsi_{i_s}^{(k)},
$$ 
 and the remaining analytic (at $\lambda_k$) solutions
$$
\vpsi_{j_1}^{(k)},~\vpsi_{j_2}^{(k)},~...~,~\vpsi_{j_l}^{(k)}
$$
We construct another set of  $s-1$  independent analytic (at $\lambda_k$) solutions: 
$$ 
 \varphi_{i_1}^{(k)}:={1\over r_{i_1}^{(k)}}\vPsi_{i_1}^{(k)}-{1\over r_{i_s}^{(k)}}\vPsi_{i_s}^{(k)},
$$
$$ 
 \varphi_{i_2}^{(k)}:={1\over r_{i_2}^{(k)}}\vPsi_{i_2}^{(k)}-{1\over r_{i_s}^{(k)}}\vPsi_{i_s}^{(k)},
$$
$$
\vdots
$$
$$ 
 \varphi_{i_{s-1}}^{(k)}:={1\over r_{i_{s-1}}^{(k)}}\vPsi_{i_{s-1}}^{(k)}-{1\over r_{i_s}^{(k)}}\vPsi_{i_s}^{(k)},
$$
It follows  that there always exist $n-1$ linearly independent vector solution which are analytic at $\lambda_k$, namely
$$
 \varphi_{i_1},~...~,\varphi_{i_{s-1}};~~~\vpsi_{j_1}^{(k)},~\vpsi_{j_2}^{(k)},~...~,~\vpsi_{j_l}^{(k)};~~~\vPsi_k,
$$
Moreover, there also exists the singular solution $
\vPsi_k(\lambda)\ln(\lambda-\lambda_k)+\hbox{reg}(\lambda-\lambda_k)
$. 
This proves the lemma. $\Box$.

\vskip 0.2 cm 
\noindent
{\large \bf Conclusion:} The four cases above are summarized below  (letting $0!=1$):
 \be
\label{psik0}
    \vPsi_k(\lambda)=
\left\{
\matrix{\left(\Gamma(\lp_k+1)\vec{e}_k+O(\lambda-\lambda_k)\right)~(\lambda-\lambda_k)^{-\lp_k-1} & ~~~\hbox{ Case 1): }\lp_k\not\in\mathbb{ Z}.
\cr
\cr
\left({(-1)^{N_k}\over (-N_k-1)!}\vec{e}_k+O(\lambda-\lambda_k)\right)(\lambda-\lambda_k)^{-N_k-1}& ~~~\hbox{ Case 2), 4): }\lp_k=N_k\in \mathbb{Z}_{-}.
\cr
\cr
\sum_{j\neq k}r_j^{(k)}\psi_j^{(k)}(\lambda)=\hbox{reg}(\lambda-\lambda_k)& ~~~\hbox{ Case 3): }\lp_k\in
\mathbb{ N}.
}
\right.
\ee

\noindent 
Moreover, there exists a singular solution given by 
\be 
\label{GNIS}
 \vPsi_k^{(sing)}(\lambda)
:=
\left\{
\matrix{\vPsi_k(\lambda),& \lp_k\not\in\mathbb{Z}&  \hbox{ i.e. (\ref{asso1})},
\cr
\cr
\vPsi_k(\lambda)\ln(\lambda-\lambda_k)+\hbox{reg}(\lambda-\lambda_k), & \lp_k=-1 & \hbox{ i.e.  (\ref{asso2})},
\cr
\cr
\vPsi_k(\lambda)\ln(\lambda-\lambda_k)+{ P_{N_k}^{(k)}(\lambda)
\over (\lambda-\lambda_k)^{N_k+1}}
+\hbox{reg}(\lambda-\lambda_k),& \lp_k\in\mathbb{N} & \hbox{ i.e.  (\ref{asso3})},
\cr
\cr
\left\{
\matrix{
\vPsi_k(\lambda)\ln(\lambda-\lambda_k)+\hbox{reg}(\lambda-\lambda_k),
\cr
\cr
\vPsi_k^{(sing)}\equiv 0,~~\hbox{ if } \underline{r}^{(k)}=0,
}
\right.
& \lp_k\in -\mathbb{N}-2 & \hbox{ i.e.   (\ref{STARE})}.
}
\right.
\ee

\vskip 0.3 cm 
\noindent
The singular part of $ \vPsi_k^{(sing)}(\lambda)$   is uniquely determined. In logarithmic case of (\ref{asso2}), (\ref{asso3}) and (\ref{STARE}),  $ \vPsi_k^{(sing)}(\lambda)$  is defined modulo the  addition of a linear combination of regular solutions.

\bde
\label{ofconnection}
The connection coefficients $c_{jk}$, $1\leq j,k\leq n$, are uniquely defined by
$$
\left\{
\matrix{
\vPsi_k(\lambda) =\vPsi_j^{(sing)}(\lambda)c_{jk}+\hbox{\rm reg}(\lambda-\lambda_j),~~~~~~~~~~~~~~~~~~~~~~~~~~~
\cr
\cr
~c_{jk}:=0,~1\leq k \leq n,~\hbox{ when } ~ \vPsi_j^{(sing)}(\lambda)\equiv 0 \hbox{ for }\lp_j\in-\mathbb{N}-2.
}
\right.
$$
\ede

\noindent
Observe that:

a) $
c_{kk}=1 $ for $ \lp_k\not\in \mathbb{ Z}$,  $c_{kk}=0$ for $\lp_k\in \mathbb{Z}.
$

 b) In case $\lp_k\in \mathbb{ N}$, it may happen that $\vPsi_k\equiv 0$. This occurs when $\vec{r}^{~(k)}=0$. In this case $c_{jk}=0$ for any $j=1,..,n$, namely the $k$-th column of the matrix $C=(c_{jk})$ is zero.

 c) In case $\lp_j\in -\mathbb{ N}-2$, it may happen that the there is no logarithmic singularity, namely $\vPsi_j^{(sing)}\equiv 0$. This occurs if $\underline{r}^{(j)}=0$. In such a case, we need to {\it define}
 $
c_{jk}:=0$, for any $k$, 
so that the matrix $C=(c_{jk})$ has zero $j$-th row. 

d) Letting $
c_{jk}:=0$, for any $k$, when $\underline{r}^{(j)}=0$, a more explicit way to write the definition of connection coefficients is (\ref{psik}).


\section{Matrix Solutions $\Psi$ and $\Psi^*$ of System (\ref{02})-(\ref{03}), Monodromy and Invertibility} 
\label{MV}

In the previous section, we have constructed a matrix solution 
\be 
\label{Pl}
 \Psi(\lambda):=[\vPsi_1(\lambda)~|~\cdots~|~\vPsi_n(\lambda)].
\ee
 In Section \ref{in} we will establish under which conditions it is fundamental. 

\bre
\label{remark1}
 System (\ref{02}), (\ref{03}) may have vector solutions that are analytic at all $\lambda_1,...,\lambda_n$. Such solutions must be polynomials in $\lambda$, because $\infty$ is a Fuchsian singularity. 
\ere

The following holds:

\ble
\label{lemma1}
  System (\ref{02}),  (\ref{03}) has no polynomial vector solution if and only if $A_1$ has no  negative integer  eigenvalues. Equivalently (see Remark \ref{remark1}),  System (\ref{02}),  (\ref{03}) has a singular solution at any $\lambda_k$, $1\leq k \leq n$, if and only if $A_1$ has no negative integer eigenvalues. 
\ele

\noindent
{\it Proof:} This Lemma is proved in remark 1.1 of \cite{BJL}. $\Box$

\vskip 0.3 cm 
In \cite{BJL} it is proved, under the {\it assumption (i)} of non integer $\lp_k$'s,  that (\ref{02}) admits a matrix solution $\Psi^*(\lambda)$, whose $k^{th}$ column, $k=1,...,n$, is  analytic at all poles $\lambda_j\neq \lambda_k$.   We prove existence of $\Psi^*$  {\it without  any assumption on $\lp_1,...,\lp_n$}. 

\bpr
\label{prop1}
Let the matrix $A_1$ be any (no assumptions). Then

i)  There exists a  matrix solution $\Psi^*=[\vPsi_1^*(\lambda)~|~\cdots~|~\vPsi_n^*(\lambda)]$ such that 
\be
\label{star1}
 \vPsi_k^*(\lambda)=\hbox{\rm reg}(\lambda-\lambda_j)~~~\forall j\neq k.
\ee
 
ii)  $\Psi^*(\lambda)$ is a {\rm fundamental matrix solution} if  and only if none of the eigenvalues of $A_1$ is a negative integer. In this case,  $\vPsi_k^{(sing)}(\lambda)\neq 0$ for any $k$, and $\vPsi_k^*(\lambda)$  has the following behaviour for $\lambda$ close to $\lambda_k$
$$
\vPsi_k^*(\lambda)=\vPsi_k^{(sing)}(\lambda)+\hbox{\rm reg}(\lambda-\lambda_k)~~~~~~~~~~~~~~~~~~~~~~~~~~~~~~~~~~~~~~~~~~~~~~~~~~~~~~~~~
$$

\be
\label{star2}
=
\left\{
\matrix{\vPsi_k(\lambda)+\hbox{\rm reg}(\lambda-\lambda_k) & \hbox{if $\lp_k\not\in \mathbb{Z}$},
\cr
\cr
 \vPsi_k(\lambda)\ln(\lambda-\lambda_k)+\hbox{\rm reg}(\lambda-\lambda_k)& \hbox{if $\lp_k\in\mathbb{Z}_{-}$},
\cr
\cr
 \vPsi_k(\lambda)\ln(\lambda-\lambda_k)+{P_{N_k}^{(k)}(\lambda)\over (\lambda-\lambda_k)^{N_k+1}}+\hbox{\rm reg}(\lambda-\lambda_k)& \hbox{if $\lp_k\in  \mathbb{N}$}.
}
\right.
\ee
  $\Psi^*(\lambda)$ is uniquely defined by (\ref{star1}) and (\ref{star2}), and 
\be
\label{C}
 \Psi(\lambda)=\Psi^*(\lambda)C,~~~~~C:=(c_{jk}).
\ee

\epr

\vskip 0.2 cm 
\noindent
{\it Proof:} See the Appendix. $\Box$

\bre 

From the above proposition, we see that if none of the eigenvalues of $A_1$ is a negative integer and   $\lp_k\in-\mathbb{N}-2$, then    $\underline{r}^{(k)}\neq 0$, namely $\vPsi_k^{(sing)}\neq 0$. 
For any $k$, the solution $\vPsi_k^*$ is always singular at $\lambda_k$. Indeed, if $\lp_k\in -\mathbb{N}-2$, by statement  above,  $\underline{r}^{(k)}\neq 0$, so there is a log-singular solution; if $\lp_k=-1$ there always is a log-singular solution; if $\lp_k\in \mathbb{N}$, there always  is a solution with at least the pole from   $P^{(k)}(\lambda)/(\lambda-\lambda_k)^{N_k+1}$.  

\ere

\subsection{Monodromy of $\Psi$ and $\Psi^*$ associated to a loop around $\lambda_k$}

Consider a small loop in ${\cal P}_\eta$ around a pole $\lambda_k$ in counter-clockwise direction, not encircling the other poles; for example $\lambda-\lambda_k\mapsto (\lambda-\lambda_k)e^{2\pi i }$, $|\lambda-\lambda_k|$ small. Monodromy of $\Psi=[\vPsi_1|...|\vPsi_n]$ is easily computed from (\ref{psik0}), which immediately implies 
$$
\vPsi_k(\lambda)\longmapsto
\left\{
\matrix{
\vPsi_k(\lambda)~e^{-2\pi i \lp_k}, & \lp_k\not\in\mathbb{Z}  \cr
\cr 
\vPsi_k(\lambda), & \lp_k\in\mathbb{Z}
}
\right.
$$ 
and from (\ref{psik}), which implies (note that $j$ and $k$ are exchanged here):
$$
\vPsi_j(\lambda)\longmapsto
\left\{
\matrix{
\vPsi_j(\lambda)+(e^{-2\pi i \lp_k}-1)c_{kj}~\vPsi_k(\lambda), & \lp_k\not\in\mathbb{Z}\cr
\cr 
\vPsi_j(\lambda)+2\pi i c_{kj}~\vPsi_k(\lambda), & \lp_k\in\mathbb{Z}
}\right.
$$ 
These formulae make sense also when $c_{kj}=0$ for any $k$ in the special case $\vPsi_j= 0$, possibly occurring  when  $\lp_k\in \mathbb{N}$, and when  $c_{kj}=0$ for any $j$ in the special case $\vPsi_k^{(sing)}\equiv 0$, possibly occurring when   $\lp_k\in -\mathbb{N}-2$. 

\vskip 0.3 cm 
\noindent
{\small {\it Proof:} Indeed, we have: 

\noindent
a) in case $\lp_k\not\in \mathbb{Z}$:
$$
\vPsi_j=\vPsi_kc_{kj}+\hbox{reg}(\lambda-\lambda_k)\longmapsto \vPsi_k e^{-2\pi i \lp_k} c_{kj}+\hbox{reg}(\lambda-\lambda_k)
$$
$$\equiv \vPsi_k e^{-2\pi i \lp_k} c_{kj}+\vPsi_j-\vPsi_kc_{kj}. 
$$
b) In case $\lp_k\in\mathbb{N}$  we have 
$$\vPsi_j=\left(\vPsi_k\ln(\lambda-\lambda_k)+{P^{(k)}\over (\lambda-\lambda_k)^{N_k+1}}\right)c_{kj}+\hbox{reg}(\lambda-\lambda_k)\longmapsto 
$$
$$
\longmapsto
\left( 2\pi i \vPsi_k + \vPsi_k\ln(\lambda-\lambda_k)+{P^{(k)}\over (\lambda-\lambda_k)^{N_k+1}}\right)c_{kj}+\hbox{reg}(\lambda-\lambda_k)
=
2\pi i c_{kj}\vPsi_k +\vPsi_j
$$
c) In case $\lp_k\in \mathbb{Z}_{-}$ ,   we have 
$$\vPsi_j=\vPsi_k\ln(\lambda-\lambda_k)c_{kj}+\hbox{reg}(\lambda-\lambda_k)\longmapsto 
$$
$$
\longmapsto
\left( 2\pi i \vPsi_k + \vPsi_k\ln(\lambda-\lambda_k)\right)c_{kj}+\hbox{reg}(\lambda-\lambda_k)
=
2\pi i c_{kj}\vPsi_k +\vPsi_j
$$
In  c) with  $\lp_k\leq -2$, in case it happens that   $\vPsi_k^{(sing)}\equiv 0$, we have 
$$\vPsi_j=0+\hbox{reg}(\lambda-\lambda_k)\longmapsto \vPsi_j
$$
This last fits into the general formula $\vPsi_j
\mapsto \vPsi_j(\lambda)+2\pi i c_{kj}~\vPsi_k(\lambda)$, because by definition  $c_{kj}:=0$ for any $j$ in this case. 
$\Box$}
 
\vskip 0.3 cm 
\noindent
Next, we compute the monodromy of  $\Psi^*=[\vPsi_1^*|\cdots|\vPsi_n^*]$, which exists when $A_1$  has no negative integer eigenvalues. We consider again a small loop around $\lambda_k$ as above. We have
$$
 \vPsi_j^*(\lambda)\longmapsto
\vPsi_j^*(\lambda)~~~~~\forall j=1,...,n,~~~ j\neq k,~~~~~~~~~~~~~~~~~~~~~~~~~~~~~~~~
$$

$$ \vPsi_k^*(\lambda)\longmapsto
\left\{
\matrix{
e^{-2\pi i\lp_k}\vPsi_k^*(\lambda)+(e^{-2\pi i\lp_k}-1)\sum_{j\neq k}c_{jk}\vPsi_j^*(\lambda),& 
~~~\lp_k\not\in\mathbb{Z}
\cr
\cr
\vPsi_k^*(\lambda)+2\pi i \sum_{j\neq k} c_{jk}\vPsi_j^*(\lambda), &~~~ \lp_k\in \mathbb{Z}
}
\right.
$$
{\small {\it Proof:} Invariance of $\vPsi_j^*$  follows from (\ref{star1}). The only singular at $\lambda_k$ solution is $\vPsi_k^*$. We use (\ref{C}) and invariance of $\vPsi_j^*$. For $\lp_k\not\in\mathbb{Z}$:
$$
\vPsi_k^*=\vPsi_k-\sum_{j\neq k}\vPsi_j^* c_{jk}~\longmapsto~ e^{-2\pi i\lp_k}\vPsi_k-\sum_{j\neq k}\vPsi_j^* c_{jk}
$$
$$
\equiv
e^{-2\pi i\lp_k}(\vPsi_k^*+\sum_{j\neq k}\vPsi_j^*c_{jk})-\sum_{j\neq k}\vPsi_j^* c_{jk}.
$$
For $\lp_k\in\mathbb{Z}$,  we use the behaviour of $\vPsi_k^*$ at $\lambda_k$ and  (\ref{C}) with $c_{kk}=0$ (in the formula below
 $P^{(k)}_{N_k}\equiv 0$ when $N_k\leq -1$): 
$$
 \vPsi_k^*=\vPsi_k\ln(\lambda-\lambda_k)+{P_{N_k}^{(k)}\over (\lambda-\lambda_k)^{N_k+1}}+\hbox{reg}(\lambda-\lambda_k)\longmapsto 
$$
$$
\longmapsto~ 2\pi i \vPsi_k+\left\{ \vPsi_k\ln(\lambda-\lambda_k)+{P^{(k)}\over (\lambda-\lambda_k)^{N_k+1}}+\hbox{reg}(\lambda-\lambda_k)\right\}
$$
$$
\equiv  2\pi i \vPsi_k+\vPsi_k^*~=2\pi i  \sum_{j\neq k} \vPsi_j^* c_{jk} + \vPsi_k^*.
$$
$\Box$
}
\vskip 0.2 cm
\noindent 
We summarize in the following
\bpr
\label{promod}
The monodromy matrices representing the monodromy of $\Psi$ and $\Psi^*$ for a small counter-clockwise loop around $\lambda_k$ in ${\cal P}_\eta$ are as follows. 

a) The matrix $\Psi$ is always defined. The monodromy is
$$
\Psi\mapsto \Psi M_k,
~~~~~~~~
 M_k=I+\alpha_k\pmatrix{0 &   0    & \cdots &   0   & \cdots  &0       \cr
        \vdots&\vdots & &   \vdots     &     & \vdots      \cr
       c_{k1} & c_{k2} & \cdots& c_{kk} & \cdots &c_{kn}      \cr
        \vdots &\vdots  &            & \vdots    &        &\vdots      \cr
             0    &      0   &\cdots   &       0         &   \cdots &0              },~~~~~1\leq k \leq n,
$$
where $I$ is the $n\times n$ identity matrix,  only the $k$-th row  in the second matrix is non zero, and 
$$\left\{
\matrix{
\alpha_k:=(e^{-2\pi i \lp_k}-1), & \hbox{ if }\lp_k\not\in \mathbb{Z}
\cr
\cr
\alpha_k:=2\pi i ,& \hbox{ if }\lp_k\in \mathbb{Z}
}
\right.        
$$

 b) If $A_1$ has no negative integer eigenvalues, then $\Psi^*$ exists. The monodromy is
$$
\Psi^*\mapsto \Psi^* M_k^*,
~~~~~
M_k^*=I+\alpha_k
\pmatrix{0 &   0    & \cdots &   c_{1k}   & \cdots  &0       \cr
               0 &   0     &       &    c_{2k}        &           & 0         \cr
        \vdots&\vdots &\ddots &\vdots   &           & \vdots      \cr
             0      & 0 & \cdots& c_{kk} & \cdots & 0      \cr
        \vdots &\vdots  &            & \vdots    &    \ddots    &\vdots      \cr
             0    &      0   &\cdots   &       c_{nk}       &   \cdots &0               }
$$
where only the $k$-th column in the second matrix is non zero.

\epr
\bre The matrix $M_k$ is the matrix $(m_{ij}^{(k)})$ in Proposition I of the Introduction. 
 For a clockwise loop, we analogously find that 
$$
[M_k^{-1}]_{kj}=\beta_k ~c_{kj},~~~j\neq k;~~~~~~~[M_k^{-1}]_{kj}=0~~\hbox{ otherwise};
$$
$$
[M_k^{-1}]_{jj}=1,~~~~~j\neq k;~~~~~~~[M_k^{-1}]_{kk}=e^{2\pi i\lp_k};
$$
where 
$$~~~~~\left\{
\matrix{
\beta_k:=(e^{2\pi i \lp_k}-1), & \hbox{ if }\lp_k\not\in \mathbb{Z},
\cr
\cr
\beta_k:=-2\pi i ,& \hbox{ if }\lp_k\in \mathbb{Z},
}
\right. ~~~~~~\Longrightarrow~~~ \beta_k=-e^{2\pi i \lp_k} \alpha_k.      $$
Moreover 
$$
[(M_k^*)^{-1}]_{jk}=\beta_k~ c_{jk},~~~j\neq k;~~~~~~[(M_k^*)^{-1}]_{ji}=0~~~\hbox{ otherwise};
$$
$$
[(M_k^*)^{-1}]_{jj}=1,~~~~~j\neq k;~~~~~~
[(M_k^*)^{-1}]_{kk}=[M_k^{-1}]_{kk}.
$$

\ere

\bre  One can define the coefficients $c_{kj}$, for $j\neq k$, starting from $\Psi$ and its monodromy matrices, as  $c_{kj}:=m_{kj}/\alpha_k$. 
\ere

\bcr 
The first invariants of the monodromy matrices in Proposition \ref{promod} are
$$
\hbox{\rm Tr}(M_k)=n-1+e^{-2\pi i \lp_k}
$$
$$
\hbox{\rm Tr}(M_jM_k)=
n-2 +e^{-2\pi i \lp_j}+e^{-2\pi i \lp_k}+\alpha_j\alpha_k~c_{jk}c_{kj}
$$
$$
=n-2 +e^{-2\pi i \lp_j}+e^{-2\pi i \lp_k}+e^{-2\pi i (\lp_j+\lp_k)}\beta_j\beta_k~c_{jk}c_{kj}
$$
If $A_1$ has no negative integer eigenvalues, then 
$$\hbox{\rm Tr}(M_k^*)=\hbox{\rm Tr}(M_k)
,~~~~~~
\hbox{\rm Tr}(M_j^*M_k^*)=\hbox{\rm Tr}(M_jM_k).
$$
\ecr
From Proposition \ref{prop1} we know that $\Psi^*$ is fundamental if and only if $A_1$ has no negative integer eigenvalues. Thus:
\bcr 
 Suppose that $A_1$ has no negative integer eigenvalues; then $M_1^*,...,M_n^*$ generate the monodromy group of equation (\ref{02}--\ref{03}).
\ecr

\subsection{On the Invertibility of $C$ and $\Psi(\lambda)$}
\label{in}

We establish necessary and sufficient conditions for the matrices $\Psi(\lambda)$ and $C=(c_{jk})$ to be invertible. Let $\lambda\in {\cal P}_\eta$. 

\bre
\label{teverde}
If $\vec{r}^{~(k)}=0$ (case $\lp_k\in \mathbb{N}$) then $C$ has  zero $k$-th column and also 
$\Psi(\lambda)$ has zero $k$-th column, so it is not a fundamental matrix. If $\underline{r}^{(k)}=0$ (case $\lp_k\in -\mathbb{N}-2$), then $C$ has zero $k$-th row. In both cases, $C$ is not invertible. 
\ere

\ble
\label{troppote}

{\bf i)} If $A_1$ has no negative integer eigenvalues and $\Psi(\lambda)$ is fundamental, then $C $ is invertible. {\bf ii)} Conversely, if $C$ is invertible, then: 

-- $A_1$ has no negative integer eigenvalues, 

-- $\Psi(\lambda)$ is fundamental, 

-- the matrix defined by $\Psi^*(\lambda):=\Psi(\lambda)C^{-1},
$
is the unique  fundamental solution $\Psi^{*}$ of Proposition \ref{prop1}, 

-- in case $\lp_k\in \mathbb{N}$ then $\vec{r}^{~(k)}\neq 0$, in case $\lp_k\in -\mathbb{N}-2$, then   $\underline{r}^{(k)}\neq 0$.  
\ele

\noindent
{\it Proof:} 
{\bf i)}  If $A_1$ has no negative integer eigenvalues, then the fundamental matrix $\Psi^*(\lambda)$ exists form Proposition \ref{prop1}. If $\Psi(\lambda)$ is invertible, then $C=\Psi^*(\lambda)^{-1}\cdot\Psi(\lambda)$ is invertible. 

{\bf ii)} From Remark \ref{teverde}, we see that $C$ invertible implies that $\vec{r}^{~(k)}\neq 0$ and $\underline{r}^{(k)}\neq 0$, when defined. In particular, in any row and any column of $C$ there is a $c_{ij}\neq 0$ for some $i\neq j$.  Write $\Psi(\lambda)$ at $\lambda_k$:
$$
\Psi=\Bigl[~\vPsi_1~|~\cdots~| ~\vPsi_n~\Bigr]~=\Bigl[~\vPsi_k^{(sing)}c_{k1}~|~\cdots~|~\vPsi_k^{(sing)}c_{kn}~\Bigr]+\hbox{reg}(\lambda-\lambda_k)
$$
$$
= \left\{ \Bigl[~0~\left|~...~\left|~0~\left|~\vPsi_k^{(sing)}~\right|~0~\right|~...~\right|~0~\Bigr]+\hbox{reg}(\lambda-\lambda_k)\right\}C.
$$
 The last step is possible because existence of $C^{-1}$ allows to write 
 $$\hbox{reg}(\lambda-\lambda_k)=\hbox{reg}(\lambda-\lambda_k)C^{-1}C \equiv \hbox{reg}(\lambda-\lambda_k) C.$$
Thus
$$
\Psi C^{-1}=  \Bigl[~0~\left|~...~\left|~0~\left|~\vPsi_k^{(sing)}~\right|~0~\right|~...~\right|~0~\Bigr]+\hbox{reg}(\lambda-\lambda_k).
$$
This is equivalent to (\ref{star1}) and (\ref{star2}), which implies that there exist the unique fundamental matrix  $\Psi^* \equiv\Psi C^{-1}$. From Proposition \ref{prop1} we conclude that $A_1$ has no negative integer eigenvalues. Obviously, it follows also that $\Psi=\Psi^*C$ is invertible. $\Box$.

\bpr
\label{propoC}
$C$ is invertible $\Longleftrightarrow$ $A_1$ has no integer eigenvalues.
\epr

\noindent
{\it Proof:} The "$\Longrightarrow$" is proved in the previous lemma, point ii). The proof of  "$\Longleftarrow$" is analogous to that of proposition 2 in \cite{BJL}, which we repeat here without assumptions on $\lp_1,...,\lp_n$.  Since $A_1$ has no negative integer eigenvalues,  there exists the unique fundamental $\Psi^*$. Therefore, the monodromy group is generated by $M_1^*,...,M_n^*$. We consider the monodromy $M^*_\infty$ at infinity, for a counterclockwise loop encircling all the poles. For purpose of this proof we can numerate the poles in such a way that the ray  $L_{k+1}$ is to the left of the ray $L_k$. Thus, 
$$ 
 M_\infty^*=M_n^*\cdots M_1^*.
$$
The behavior of system (\ref{03}) at $\infty$ is
$$ 
 {d\Psi\over d\lambda}= -{A_1+1\over \lambda} \left[I+O\left({1\over \lambda}\right)\right]\Psi.
$$
This implies that $A_1$ has no integer eigenvalues if and only if $M_\infty^*$ has no eigenvalue $=1$. We show that this is equivalent to the fact that $C$ is invertible, namely has no zero eigenvalue. Indeed,  
existence of an eigenvalue equal to 1 means  that there exists a non zero row vector $\hat{w}=[w_1,...,w_n]$, such that $\hat{w}M_\infty^*=\hat{w}$.  Using the explicit expression of the $M_k^*$ in terms of the $\alpha_k c_{jk}$'s,  we compute
$$
 \hat{w}~M_n^*\cdots M_1^* ~=\hat{w}~+\sum_{j=1}^n b_j \hat{e}_j,
$$
where the $\hat{e}_j$'s are the basis rows
$$
\hat{e}_1=[1,0,...,0],~~~\hat{e}_2=[0,1,...,0],~~~\hat{e}_n=[0,...,0,1],
$$
and
$$b_n=\alpha_n~(\hat{w} C)_n,
$$
$$
b_{n-1}=\alpha_{n-1}\Bigl[ (\hat{w}C)_{n-1}+c_{n,n-1}b_n \Bigr],
$$
$$
\vdots
$$
$$
b_i=\alpha_i\left[  (\hat{w}C)_i+\sum_{j=i+1}^n~c_{ji}~b_j  \right]
$$
for $i=1,2,...,n-1.$  Since all the $\alpha_i$, for $1\leq i \leq n$, are not zero, we conclude that $ \hat{w} M_\infty^*=\hat{w}$ if and only if  $\hat{w} C=0$.  Thus,  $A_1$ has integer eigenvalues if and only if  $C$ has  zero eigenvalue, namely is not invertible.
$\Box$.

\bpr
\label{PropoE}
 {\bf i)} If $A_1$ has no integer eigenvalues, then  $\Psi(\lambda)$ is a fundamental matrix solution. 

{\bf ii)} With the additional assumption that  $\lp_k\not \in\mathbb{Z}$, $\forall k =1,..., n$, also the converse holds: if $\Psi(\lambda)$ is a fundamental matrix solution, then $A_1$ has no integer eigenvalues.

\epr

\noindent
{\it Proof:} i)  If $A_1$ has no integer eigenvalues, $C$ is invertible (Proposition \ref{propoC}). Therefore, the statement follows from Lemma \ref{troppote}, point ii).  

ii)  Let $\Psi=[\vPsi_1|\cdots|\vPsi_n]$  be fundamental. Observe that under the hypothesis that  $\lp_k\not\in \mathbb{Z}$ for any $k$, the columns are singular. Namely: 
 $$\vPsi_k(\lambda)\equiv \vPsi_k^{(sing)}(\lambda)=
 (\Gamma(\lp_k+1)\vec{e}_k+O(\lambda-\lambda_k))(\lambda-\lambda_k)^{-\lp_k-1},
 ~~~~~1\leq k \leq n.
 $$ 
  The monodromy of $\Psi(\lambda)$  at infinity is $M_\infty=M_nM_{n-1}\cdots M_1$. Suppose that there is an integer eigenvalues of $A_1$. It follows that there exists a non zero  column vector  $\vec{v}=(v_1,...,v_n)^t$ ($t$ means transpose)  such that $M_\infty\vec{v}=\vec{v}$. As in the proof of Proposition \ref{propoC}, making use of the explicit form of the $M_k$'s in terms of the $c_{jk}$'s, we see that $M_\infty\vec{v}=\vec{v}$ is equivalent to $C\vec{v}=\vec{0}$. Take the vector $\vpsi(\lambda)= \sum_{l=1}^nv_l\vPsi_l(\lambda)$. At every $\lambda_k$ it behaves like
$$
\vpsi(\lambda)= \sum_{l=1}^n v_l\vPsi_k c_{kl}+\hbox{reg}(\lambda-\lambda_k)= \Bigl(\sum_{l=1}^n c_{kl}v_l \Bigr)\vPsi_k+\hbox{reg}(\lambda-\lambda_k)
$$
But $\sum_{l=1}^n c_{kl}v_l=0$, thus 
$$
 \vpsi(\lambda)=\hbox{reg}(\lambda-\lambda_k),~~\hbox{ close to any } \lambda_k,~~k=1,\cdots,n.
$$
 This implies that $\vpsi(\lambda)$ is a polynomial solution. This contradicts the fact that $\vPsi_1(\lambda)$,...,   $\vPsi_n(\lambda)$ is a basis, each  $\vPsi_k(\lambda)$ being singular at $\lambda_k$, $1\leq k \leq n$. 
$\Box$
\vskip 0.2 cm

\bcr
 If $A_1$ has no integer eigenvalues, then $M_1,...,M_n$ of Proposition \ref{promod} generate the monodromy group of system (\ref{02}). 
\ecr

\bcr
Suppose $A_1$ has some integer eigenvalues and $\Psi(\lambda)$ is a fundamental matrix solution (consequently,  $M_1$, ...,$M_n$ generate the monodromy group). In such cases, at least some $\lp_k$ is necessarily integer. 

\ecr


\section{Relation between Matrices $\Psi^*(\lambda;\eta)$ as $\eta$ changes}
\label{RC} 

 Following \cite{BJL}, we call {\it critical values}  the inadmissible values for $\eta$,  namely  
 $$
\arg(\lambda_j-\lambda_k)~\hbox{ mod } 2\pi.
$$
 We numerate them as in \cite{BJL}, as follows. In the angular interval $\left(-{\pi\over2},{3\pi\over 2}\right]$ there is an even number $m$  of critical values, ordered as
$$
 {3\pi \over 2} \geq \eta_0>\eta_1<\cdots>\eta_{m-1}>-{\pi\over 2}.
$$
All the possible critical values are then 
$$
 \eta_{\nu+hm}:=\eta_n-2h\pi,~~~\nu=0,...,m-1;~~~h\in \mathbb{Z}.
$$
In each interval $(\theta-2\pi,\theta]$ there are $m$ such values. In other words, $\{\eta_\nu~|~\nu\in \mathbb{Z}\}$ is the set of all critical values. 

  There is an  ordering of the poles with respect to an admissible $\eta$, given by  the {\it dominance relation}  $\prec$ below: 

\bde 
{\rm [as in \cite{BJL}]:} 
Let $\eta$ be admissible. We say that 
   $j \prec k $, whenever in the plane ${\cal P}_\eta$  the cut $L_j$ lies to the right of the cut  $L_k$. 
 Equivalently,  choose the determinations 
$$
\eta_{jk}:= \hbox{ determination of } \arg(\lambda_j-\lambda_k)~~ s.t.~~ \eta-2\pi <\eta_{jk}<\eta,~~~ j\neq k, ~~1\leq j, k\leq n.
$$ 
Then
\be
\label{domexam}
 j\prec k ~~~~~~\Longleftrightarrow ~~~ -\pi+\eta<\eta_{jk}<\eta.
\ee 
\ede
 The reason for the nomenclature "dominance" will be explained in section \ref{StoKM}. 

\vskip 0.3 cm
\noindent 
{\bf Remark:}  $\lambda_1,...,\lambda_n$ are in {\it lexicographical order with respect to the admissible $\eta$} when the labelling order $j<k$ coincides with the dominance order $j\prec k $.

\vskip 0.3 cm 

The matrices $\Psi^{(k)}(\lambda)$,  $\Psi(\lambda)$ and $\Psi^*(\lambda)$ defined in the plane ${\cal P}_\eta$, with $\eta$ admissible, and the connection matrix $C$, depend on $\eta$. Therefore we write 
 $$
\Psi^{(k)}(\lambda)=\Psi^{(k)}(\lambda,\eta),~~~~~  \Psi(\lambda)=\Psi(\lambda,\eta),$$
$$
\Psi^*(\lambda)=\Psi^*(\lambda,\eta),~~~~~~~C=C(\eta).
$$

\begin{figure}
\centerline{\includegraphics[width=0.7\textwidth]{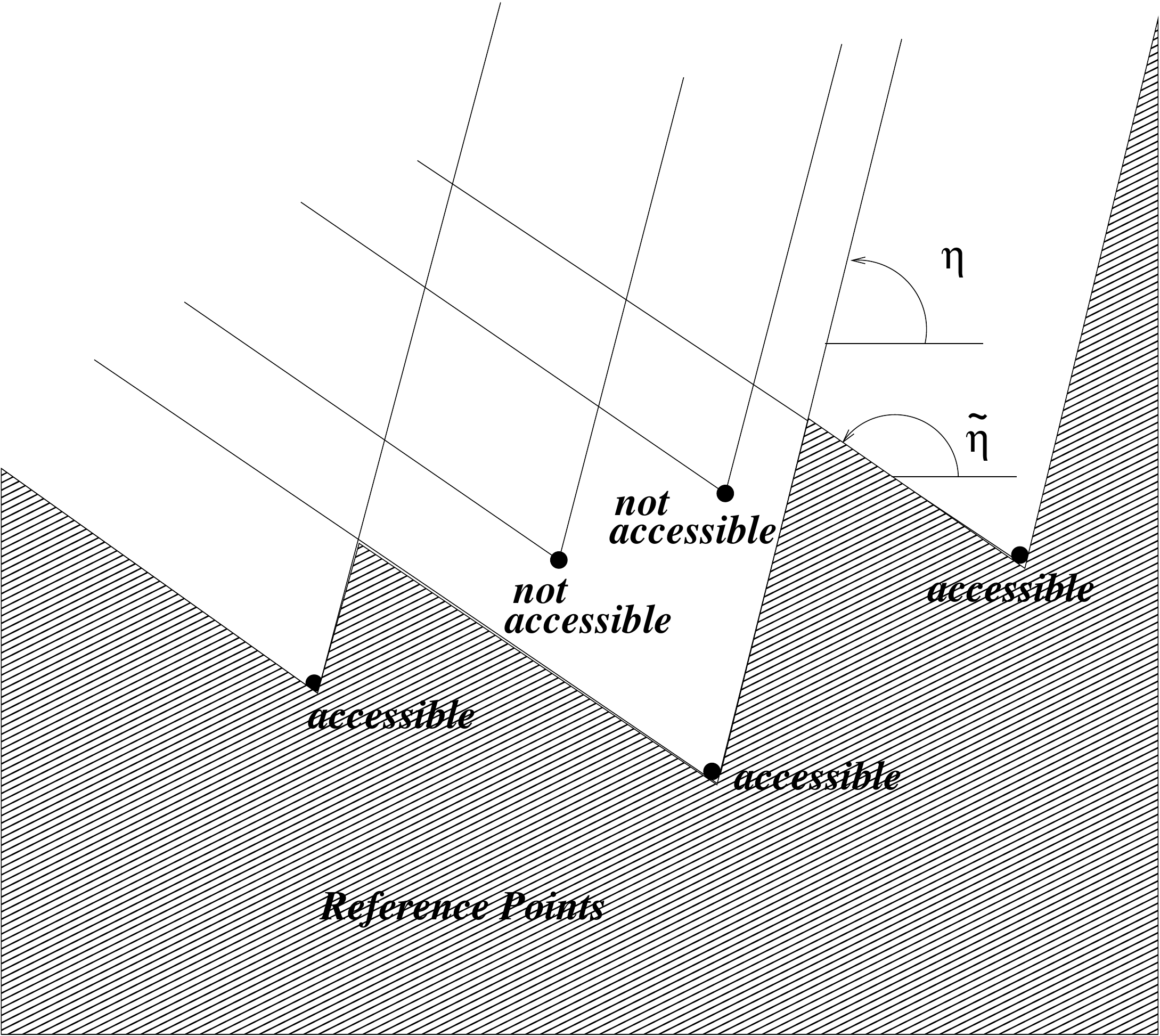}}
\caption{Picture of ${\cal P}_\eta\cap {\cal P}_{\tilde{\eta}}$.}
\label{fig1}
\end{figure}

For two values  $\eta<\tilde{\eta}$, we 
 consider the plane with both the cuts of
  ${\cal P}_\eta$ and ${\cal P}_{\tilde{\eta}}$.
 We denote by ${\cal P}_\eta\cap{\cal P}_{\tilde{\eta}}$ the simply connected set of {\it reference points w.r.t.  ${\cal P}_\eta$ and ${\cal P}_{\tilde{\eta}}$}, namely the points 
in the doubly cut plane such that $\arg(\lambda-\lambda_k)\not\in [\eta,\tilde{\eta}]$, $\forall k=1,...,n$. A pole $\lambda_j$ is called {\it accessible} if it is on the boundary of ${\cal P}_\eta\cap{\cal P}_{\tilde{\eta}}$. See figure \ref{fig1}.

We  generalize propositions 3  of \cite{BJL} without any assumptions on  the values of $\lp_1,...,\lp_n$.

\bpr 
\label{etto}
\noindent
i) Let  $\lambda_k$ be accessible w.r.t.  some admissible $\eta$ and $\tilde{\eta}$. Then  
$$\vPsi_k^{(k)}(\lambda,\eta)=\vPsi_k^{(k)}(\lambda,\tilde{\eta}) ~~\hbox{ and  }~~\vPsi_k(\lambda,\eta)=\vPsi_k(\lambda,\tilde{\eta}), ~~~\forall \lambda\in{\cal P}_\eta\cap{\cal P}_{\tilde{\eta}}$$

\noindent
ii) Let  $\eta$ and $\tilde{\eta}$ lie between two consecutive critical values: namely $\eta_{\nu+1}<\eta<\tilde{\eta}<\eta_{\nu}$. Then 
$$ 
C(\eta)=C(\tilde{\eta})$$

\noindent
iii) Let again   $\eta_{\nu+1}<\eta<\tilde{\eta}<\eta_{\nu}$. Then 
$$
 \Psi^*(\lambda,\eta)=\Psi^*(\lambda,\tilde{\eta}),~~~\forall \lambda\in{\cal P}_\eta\cap{\cal P}_{\tilde{\eta}}.
$$
whenever $\Psi^*$ is uniquely  defined (namely, when $A_1$ has no negative integer eigenvalues).
\epr

\vskip 0.2 cm 
\noindent
{\it Proof:} See the Appendix. $\Box$

\vskip 0.3 cm

The above implies  that the dependence on  $\eta$ is discrete, namely it changes when a critical value is crossed. Thus, if $\eta_{\nu+1}<\eta<\eta_\nu$, we follow \cite{BJL} and write 
$$
 \Psi_\nu(\lambda):=\Psi(\lambda,\eta),~~~\Psi_\nu^*(\lambda):=\Psi^*(\lambda,\eta),~~~C_\nu=(c_{jk}^{(\nu)}):=C(\eta).
$$
We now compute how $\Psi^*(\lambda,\eta)$ changes when $\eta $ changes, so  generalizing  proposition 4 of \cite{BJL}, without assumptions on $\lp_1,..,\lp_n$. 
\bpr
\label{propoW}
 Suppose that  $A_1$ has no negative integer eigenvalues, so that the $\Psi^*_\nu(\lambda)$'s exist, for $\nu\in \mathbb{Z}$. 
Let  $\eta_{\nu+1}<\eta<\eta_\nu<\tilde{\eta}<\eta_{\nu-1}$. Then 
\be  
\label{BEC0}
 \Psi^*_{\nu-1}(\lambda)=\Psi^*_\nu(\lambda) ~W_\nu,~~~~\forall \lambda\in{\cal P}_\eta\cap{\cal P}_{\tilde{\eta}},
\ee
where the invertible matrix  $W_\nu=(W_{jk}^{(\nu)})$ is
\be 
\label{BEC1}
W_{jk}^{(\nu)}= -\alpha_k c_{jk}^{(\nu)} ,~~\hbox{ for } j\succ k ~\hbox{ such that }~ \arg(\lambda_j-\lambda_k)=\eta_\nu 
 \ee
\be 
\label{BEC2}
W_{jj}^{(\nu)}=1,~~~j=1,\cdots,n;~~~~~~~~W_{jk}^{(\nu)}=0 \hbox{ otherwise}.
\ee
where $\prec$ is the dominance relation w.r.t. $\eta$.  In the same way, 
$$
 \Psi^*_{\nu}(\lambda)=\Psi^*_{\nu-1}(\lambda) ~W_\nu^{-1},~~~~\forall \lambda\in{\cal P}_\eta\cap{\cal P}_{\tilde{\eta}},
$$
where $W_\nu^{-1}$ 
has zero entries except for 
$$
[W_\nu^{-1}]_{jj}=1,~~~j=1,\cdots,n,
$$
$$
[W_\nu^{-1}]_{jk}= -\beta_k c_{jk}^{(\nu-1)} ,~~\hbox{ for } j \succ k ~\hbox{ s.t. }~ \arg(\lambda_j-\lambda_k)=\eta_\nu .
 $$
\epr
 Note that $W_{jk}^{(\nu)}= -\alpha_k c_{jk}(\eta) $ for $\arg(\lambda_j-\lambda_k)=\eta_\nu$ implies that  $~j\succ k$, 
\vskip 0.2 cm 
\noindent
{\it Proof:} See the Appendix. $\Box$

\vskip 0.3 cm 
In an angular interval $(\theta-2\pi,\theta]$, there are $m$ critical values,  $m$  even. Let $\mu=m/2$. Let $\eta_{\nu+1}<\eta<\eta_\nu$ and  introduce, as in \cite{BJL}, the matrices $C_\nu^+$ and $C_\nu^{-}$ such that
\be
\label{cpiu}
\Psi_{\nu+\mu}^*(\lambda)=\Psi_\nu^*(\lambda)C_\nu^+, ~~~~~~~\lambda\in{\cal P}_\eta\cap{\cal P}_{\eta-\pi}.
\ee
\be
\label{cmeno}
\Psi_{\nu+\mu}^*(\lambda)=\Psi_{\nu+m}^*(\lambda)C_\nu^-, ~~~~~~~\lambda\in{\cal P}_{\eta-\pi}\cap{\cal P}_{\eta-2\pi}
\ee
Immediately it follows that 
\be 
\label{ciclo}
 C_\nu^+=(W_{\nu+\mu}\cdots W_{\nu+1})^{-1},~~~C_\nu^-=W_{\nu+m}\cdots W_{\nu+\mu+1}.
\ee
Note\footnote{
Equivalent way to write (\ref{cpiu}) and (\ref{cmeno}):
$$\Psi^*(\lambda,\eta-\pi)=\Psi^*(\lambda,\eta)C_\nu^+,~~~\lambda\in{\cal P}_\eta\cap{\cal P}_{\eta-\pi},
$$
$$
\Psi^*(\lambda,\eta-\pi)=\Psi^*(\lambda,\eta-2\pi)C_\nu^-,~~~\lambda\in{\cal P}_{\eta-\pi}\cap{\cal P}_{\eta-2\pi}.
$$
}
\bre
 ${\cal P}_\eta\cap{\cal P}_{\eta-\pi}$ is the half plane to the left hand side of all lines whose positive parts are the cuts of direction $\eta$, while ${\cal P}_{\eta-\pi}\cap{\cal P}_{\eta-2\pi}$ is the half plane to the right hand side of all lines whose positive parts are the cuts of direction $\eta-2\pi$. 
\ere
We restate remark 3.3 of \cite{BJL} with no assumptions on $A_1$:  
\ble
\label{pippo}
Let $\Lambda^\prime:=\hbox{diag}(\lp_1,...,\lp_n)$, $\lp_k\in\mathbb{C}$, $1\leq k \leq n$. Then
$$
\Psi_{\nu+m}(\lambda)=\Psi_\nu(\lambda)e^{2\pi i \Lambda^\prime},
$$
for any $\lambda$ in the universal covering of $\mathbb{C}\backslash \{\lambda_1,...,\lambda_n\}
$. Moreover 
$$
 C_\nu=e^{2\pi i\Lambda^\prime}~C_{\nu+m}~e^{-2\pi i\Lambda^\prime}
, ~~~~(\hbox{namely: } C(\eta)=e^{2\pi i\Lambda^\prime}~C(\eta-2\pi)~e^{-2\pi i\Lambda^\prime}).
$$
\ele

\vskip 0.2 cm 
\noindent
{\it Proof:} See the Appendix. $\Box$

\vskip 0.3 cm 
\noindent
We generalize proposition 5 of \cite{BJL}, with no assumptions on  diag$(A_1)=(\lp_1,...,\lp_k)$. 

\bpr 
\label{cpmPro} 
Let $\eta_{\nu+1}<\eta<\eta_{\nu}$, and let $c_{jk}^{(\nu)} ~=c_{jk}(\eta)$. Consider the relation (\ref{cpiu}) and (\ref{cmeno}). The connection matrices $C_\nu^+$, $C_\nu^-$ are 
\be 
\label{BEC3}
[ C_\nu^+]_{jk}
=
\left\{
\matrix{
-\beta_k~c_{jk}^{(\nu)} &=&e^{2\pi i \lp_k}\alpha_k ~c_{jk}^{(\nu)}& ~~~\hbox{ for } j\prec k
\cr
\cr
             1   &&& ~~~\hbox{ for } j =k
\cr
\cr
             0   &&&  ~~~\hbox{ for } j\succ k
}
\right.
\ee
\vskip 0.3 cm 
\be 
\label{BEC4}
[ C_\nu^-]_{jk}
=
\left\{
\matrix{
             0   &&&  ~~~\hbox{ for } j\prec k
\cr
\cr
             1   &&&~~~ \hbox{ for } j =k
\cr
\cr
       e^{-2\pi i \lp_j}\beta_k~c_{jk}^{(\nu)}&= & -e^{2\pi i (\lp_k-\lp_j)}\alpha_k~c_{jk}^{(\nu)} & 
                                         ~~~\hbox{ for } j\succ k
}
\right.
\ee
\vskip 0.2 cm 
\noindent
where $\alpha_k=e^{-2\pi i \lp_k}-1$ if $\lp_k\not\in \mathbb{Z}$,  $\alpha_k=2\pi i$ if $\lp_k\in \mathbb{Z}$, $\beta_k=e^{-2\pi i \lp_k} \alpha_k$. 
\epr

\vskip 0.2 cm
\noindent
{\it Proof:} See the Appendix. $\Box$

\vskip 0.3 cm 
The matrices $C_\nu^+$ and $C_\nu^-$ {\it can be defined} by formulae (\ref{BEC3}) and (\ref{BEC4}), independently of the fact that $A_1$ has no negative integer eigenvalues, namely independently of (\ref{cpiu}) and (\ref{cmeno}). The following corollary is a direct computation.

\bcr
\label{mana} 
Let the matrices $C_\nu^+$ and $C_\nu^-$ {\rm be defined} by formulae (\ref{BEC3}) and (\ref{BEC4}). Then
$$
\hbox{\rm Tr}(M_k)=n-1+e^{-2\pi i \lp_k}
$$
$$
\hbox{\rm Tr}(M_jM_k)=
\left\{
\matrix{
n-2 +e^{-2\pi i \lp_j}+e^{-2\pi i \lp_k}-e^{-2\pi i \lp_j}~[C_\nu^+]_{jk}[C_\nu^-]_{kj} & \hbox{ if } j\prec k,
\cr
\cr
n-2 +e^{-2\pi i \lp_j}+e^{-2\pi i \lp_k}-e^{-2\pi i \lp_k}~[C_\nu^-]_{jk}[C_\nu^+]_{kj} & \hbox{ if } j\succ k.
}
\right.
$$
If  moreover $A_1$ has no integer eigenvalues, then
$$
\hbox{\rm Tr}(M_k^*)=\hbox{\rm Tr}(M_k),~~~~~~~~
\hbox{\rm Tr}(M_j^*M_k^*)=\hbox{\rm Tr}(M_jM_k)
.
$$
\ecr

\section{ Fundamental Solutions of (\ref{01}) as Laplace Integrals}
\label{FL}

\subsection{Fundamental solutions of (\ref{01}) and Stokes Matrices}
\label{StoKM}

\begin{figure}
\centerline{\includegraphics[width=0.4\textwidth]{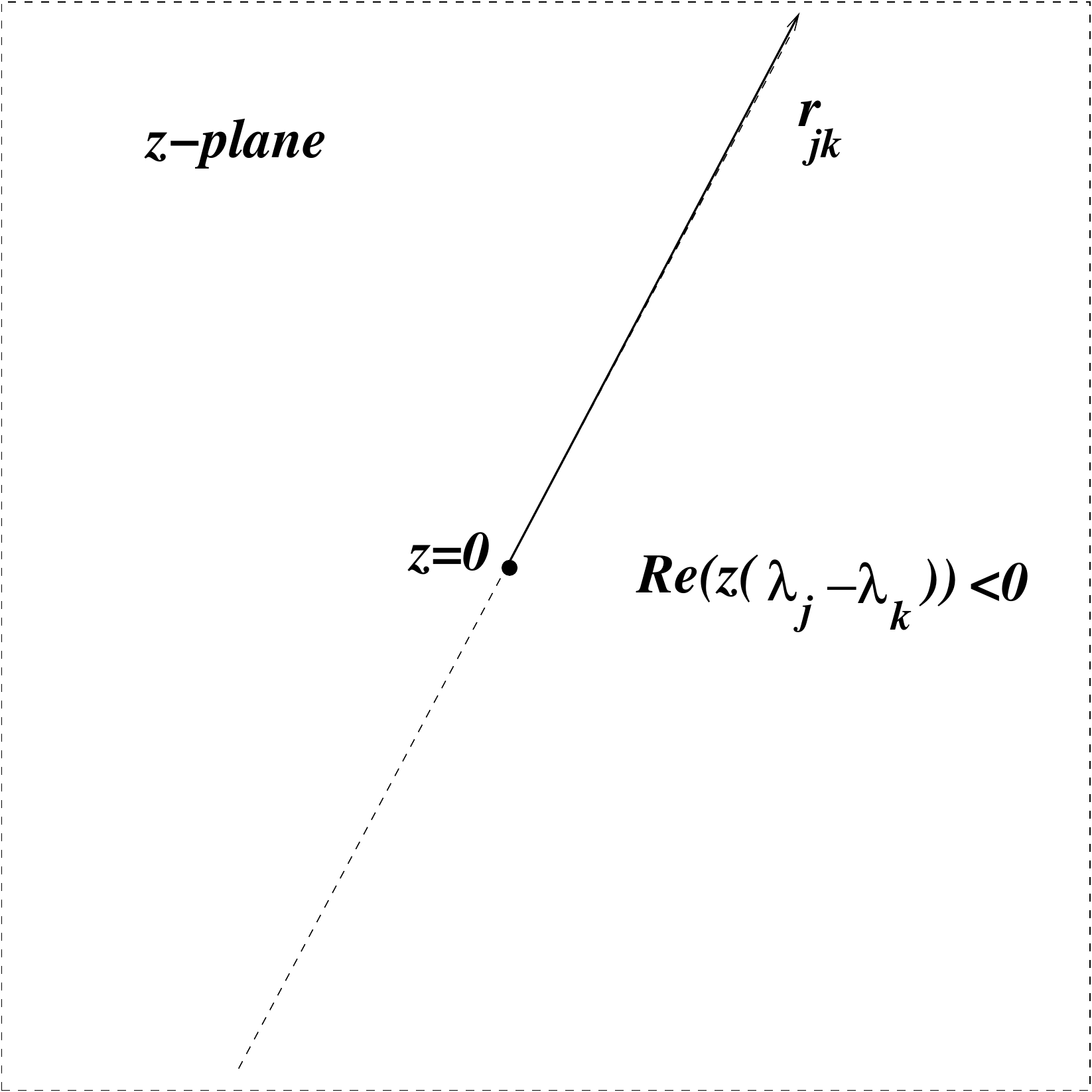}}
\caption{Stokes' ray $r_{jk}$, defined by  $\arg z={3\pi\over 2}-\eta_{jk}$ mod $2\pi$. To the right of $r_{jk}$,   $z={3\pi\over 2}-\eta$, with $\eta_{jk}<\eta<\eta_{jk}+\pi$, and  $\Re(z(\lambda_j-\lambda_k))<0$. }
\label{stokesRay}
\end{figure}

\bde
 Stokes rays are the oriented rays from 0 to $\infty$ contained in the universal covering of $\mathbb{C}\backslash \{0\}$, defined by  the condition
$$
\Re(z(\lambda_j-\lambda_k))=0,~~~i\neq j,~~~1\leq i,j\leq n,~~~~~~~z\in \widetilde{\mathbb{C}\backslash \{0\}}.
$$ 
\ede
Let $\eta\in \mathbb{R}$ be admissible, namely  $\eta\neq \arg(\lambda_i-\lambda_j) \hbox{ mod } 2\pi $, for any $i\neq j$. We choose the Stokes rays 
$$
r_{jk}:=\Bigl\{z\in \mathbb{C}~|~z=\rho \exp\left\{i\left({3\pi\over 2} -\eta_{jk}\right)\right\},~\rho>0\Bigr\},~~~j\neq k,~~~1\leq j,k\leq n,
$$
where
$$
\eta_{jk}=\hbox{ determination of } \arg(\lambda_j-\lambda_k)~~\hbox{ s.t. } \eta_{jk}\in (\eta-2\pi,\eta].
$$
If follows that
$$\left\{
\matrix{
\Re(z(\lambda_j-\lambda_k))=0, & \hbox{for $z$ belonging to the ray,}\cr
 \Re(z(\lambda_j-\lambda_k))<0, & \hbox{for $z$ in the half plane to the right of $r_{jk}$.}
}
\right.
$$
  According to the definition,  all the Stokes rays are characterized by 
$$
\arg z= {3\pi\over 2}-\eta_{jk}~\hbox{ mod }2\pi.
$$
Moreover, for any $(j,k)$ such that  $\eta_{jk}<\eta<\eta_{jk}+\pi$ we have 
$$
 \Re(z(\lambda_j-\lambda_k))<0 ~~~\hbox{\rm  if }~~ \arg z = {3\pi \over 2} -\eta \hbox{ mod }2\pi.
$$
This means that when we fix an {\it admissible} $\eta$ for system (\ref{02})-(\ref{03}), we have 
\be
\label{dominance}
\Re(z(\lambda_j-\lambda_k))<0  \hbox{ for }\arg z = {3\pi \over 2} -\eta \hbox{ mod }
2\pi~~~\Longleftrightarrow~~~j \prec k
\ee
Indeed, by (\ref{domexam}),  $j \prec k $ means that $-\pi+\eta<\eta_{jk}<\eta$. Relation (\ref{dominance}) explains why we have called $\prec$ a {\it dominance relation}: in the half plane to the right of $r_{jk}$ the eigenvalue $\lambda_k$ is dominant w.r.t. $\lambda_j$, in the usual sense of asymptotic theory of ODE with singularities of the II kind. A Stokes ray is represented in figure \ref{stokesRay}.

In the same way as all the critical values $\eta_\nu$, $\nu\in\mathbb{Z}$, are obtained from $\eta_{jk}$ by adding multiples of $2\pi$,  so  all the Stokes rays are given  by  $\arg z=\tau_\nu$, where
$$
\tau_\nu:={3\pi \over 2}-\eta_\nu,
$$
$$
0\leq \tau_0<\tau_1<\cdots<\tau_n<2\pi;~~~~~\tau_{\nu+mh}=\tau_n+2\pi h,~~~h\in \mathbb{Z}.
$$

\vskip 0.3 cm 
\noindent
Let us denote a sector in the universal covering  $\widetilde{\mathbb{C}\backslash\{0\}}$ of $\mathbb{C}\backslash\{0\}$ in the following way
$$S(\alpha,\beta):=\{z\in\widetilde{\mathbb{C}\backslash\{0\}} ~|~~~\alpha<\arg(z)<\beta\}.
$$ 
The following result is known \cite{Wasow}, \cite{BJL1}: in any sector 
$$
 {\cal S}_\nu:=S(\tau_\nu-\pi,\tau_{\nu+1}), ~~~~~\nu\in \mathbb{Z},
$$
equation (\ref{01}) has a fundamental matrix solution 
\be
\label{funda}
Y_\nu(z)=\hat{Y}_\nu(z)~e^{A_0z+\Lambda^\prime\ln z},~~~~~z\in  {\cal S}_\nu
\ee
where $\Lambda^\prime=\hbox{ diag}(A_1)$, and 
$\hat{Y}_\nu(z)$ is an invertible matrix, analytic in a neighbourhood of $\infty$, with {\it asymptotic} expansion 
\be
\label{asint}
   \hat{Y}_\nu(z)~\sim~I+{F_1\over z}+{F_2\over z^2}+\cdots~=I +\sum_{k=1}^\infty {F_k\over z^k}, ~~~~~~\hbox{ for } z\to\infty 
\hbox{ in }  {\cal S}_\nu.
\ee
 The sector ${\cal S}_\nu$ is the maximal sector where the asymptotic behavior holds, and $Y_\nu(z)$ is unique, namely it is uniquely  determined by its asymptotic behavior. The $n\times n$ matrices $F_k$ are determined as rational functions of $A_0$ and $A_1$, by formal substitution into (\ref{01}) (see \cite{Wasow}, \cite{BJL1}).

\bde [Stokes Matrices]
\label{definsto}
Given two fundamental matrices $Y_\nu$ and $Y_{\nu^\prime}$ as above, whose maximal sectors ${\cal S}_\nu$ and ${\cal S}_{\nu^\prime}$ intersect in such a way that no Stokes rays are contained in ${\cal S}_\nu\cap{\cal S}_{\nu^\prime}$, then the connection matrix $S$ such that $Y_{\nu^\prime}(z)=Y_\nu(z)S$, $z\in {\cal S}_\nu\cap{\cal S}_{\nu^\prime}$, is called a {\rm Stokes matrix}. 
\ede

\begin{figure}
\centerline{\includegraphics[width=0.6\textwidth]{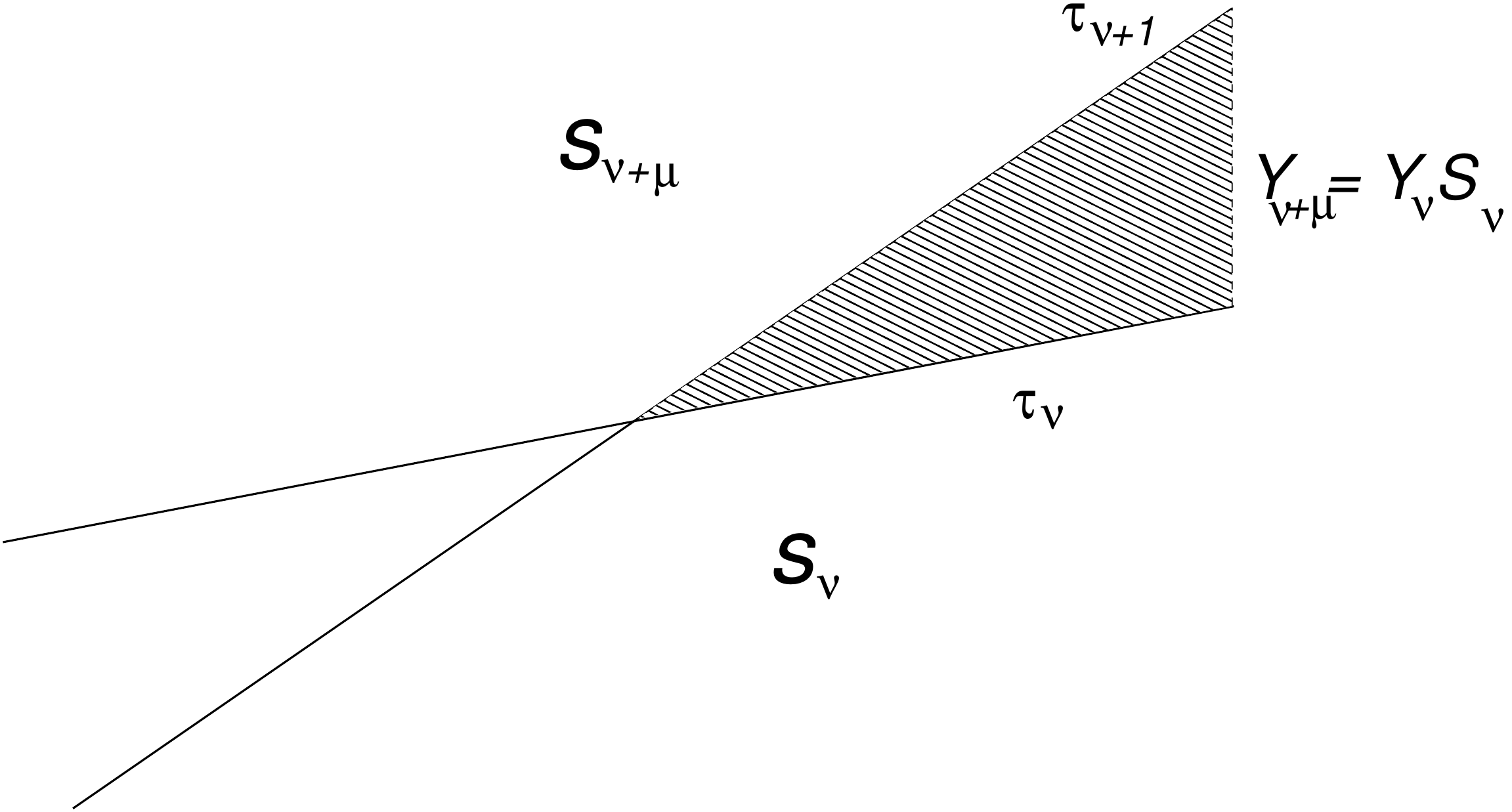}}
\caption{$
 Y_{\nu+\mu}(z)=Y_\nu(z)S_\nu$, $z\in {\cal S}_{\nu}\cap{\cal S}_{\nu+\mu}=S(\tau_\nu,\tau_{\nu+1})$, $\nu\in \mathbb{Z}.
$}
\label{fig4}
\end{figure}

\noindent
Recall that in a sector $(\theta-2\pi, \theta]$ there are $m$ (even) critical values $\eta_\nu$, $\nu=\nu_0,\nu_0+1,~...~,\nu_0+m-1$, therefore in $[{3\pi\over 2}-\theta,{3\pi\over 2}-\theta+2\pi)$ there are $m$ Stokes rays with directions $\tau_\nu$, $\nu=\nu_0,\nu_0+1,~...~,\nu_0+m-1$.  Again, let $\mu={m\over 2}$. Observe that $\eta_{\nu+\mu}=\eta_\nu-\pi$, therefore $\tau_{\nu+\mu}=\tau_\nu+\pi$. It follows that two matrices $Y_\nu$ and $Y_{\nu^\prime}$ satisfying the conditions above are precisely  $Y_\nu(z)$ and $Y_{\nu+\mu}(z)$. Therefore (see figure \ref{fig4}):

\vskip 0.3 cm 
\noindent
 {\bf Definition \ref{definsto}${}^\prime$ (Stokes Matrices)} {\it  For any  $\nu\in \mathbb{Z}$  is defined the  {\rm Stokes matrix} $S_\nu$, which is  the connection matrix such that 
$$
 Y_{\nu+\mu}(z)=Y_\nu(z)S_\nu,~~~~~z\in {\cal S}_{\nu}\cap{\cal S}_{\nu+\mu}=S(\tau_\nu,\tau_{\nu+1}).
$$
}

\vskip 0.2 cm 

\noindent
Observe that $\Re(z(\lambda_j-\lambda_k))<0$ for $j\prec k$ when $z\in S(\tau_\nu,\tau_{\nu+1})$, where the dominance relation is referred to a any $\eta\in(\eta_{\nu+1},\eta_\nu)$. From the asymptotic behaviours of $ Y_{\nu+\mu}(z)$ and $Y_\nu(z)$, it follows that $ 
\delta_{jk}\sim e^{(\lambda_j-\lambda_k)z}(S_\nu)_{jk}$ and thus 
$$
 (S_\nu)_{jj}=1,~~~~~(S_\nu)_{jk}=0 ~\hbox{ for } j\succ k.
$$

\bde [Stokes Factors] 
The {\rm Stokes factors} are the connection matrices $V_\nu$ such that 
$$
 Y_{\nu-1}(z)=Y_\nu(z)V_\nu,~~~~~z\in{\cal S}_{\nu-1}\cap{\cal S}_{\nu}=S(\tau_\nu-\pi,\tau_\nu).
$$
\ede

\noindent
It follows that $Y_{\nu+1}(z)=Y_{\nu+\mu}(z)V_{\nu+\mu}\cdots V_{\nu+1}$, (the r.h.s. is seen as the analytic continuation of the l.h.s.),  and thus
\be 
\label{ciclo1}
 S_\nu=\Bigl(V_{\nu+\mu}
\cdots V_{\nu+1}\Bigr)^{-1}.
\ee
Observe that $\tau_\nu$ and $\tau_{\nu-\mu} =\tau_\nu-\pi=-\tau_\nu$ are not contained in   $S(\tau_\nu-\pi,\tau_\nu)$, while $\tau_{\nu-1},...,\tau_{\nu-\mu+1}$ are. 
It follows that  $\Re(z(\lambda_j-\lambda_k))$ change sign in $S(\tau_\nu-\pi,\tau_\nu)$ except for $(j,k)$ such that $\arg(\lambda_j-\lambda_k)=\eta_\nu$, or $\eta_\nu-\pi$. 
Precisely,   $\Re(z(\lambda_j-\lambda_k))<0$ if $\arg(\lambda_j-\lambda_k)=\eta_\nu$, and $j\succ k$ with respect to $\eta\in(\eta_{\nu+1},\eta_\nu)$.  As above, we conclude that 
$$
 (V_\nu)_{jj}=1,
$$
$$
(V_\nu)_{jk}=0 ~~~\forall~j\neq k, \hbox{  except possibly for }  (j,k) \hbox{ s.t. } 
\left\{
\matrix{\arg(\lambda_j-\lambda_k)={3\pi\over 2} -\tau_\nu.
\cr
j\succ k \hbox{ w.r.t. } \eta\in(\eta_{\nu+1},\eta_\nu).
}
\right.
$$

\vskip 0.2 cm 
\bre 
i) The monodromy of  $Y_\nu(z)$ is completely described by the {\it monodromy data} $S_\nu$, $S_{\nu+\mu}$ and $\Lambda^\prime$, because the following holds
\be
\label{temporale}
   Y_\nu(z e^{2\pi i })=Y_\nu(z) e^{2\pi i \Lambda^\prime} (S_\nu S_{\nu+\mu})^{-1}, ~~~~~z\in {\cal S}_\nu.
  \ee
For the above reason, $S_\nu$ and  $S_{\nu+\mu}$ are a {\it complete set of Stokes multiplies}.  
Any other Stokes matrix can be expressed in terms of entries of  $S_\nu$, $S_{\nu+\mu}$ and $\Lambda^\prime$. 

\noindent
\noindent
ii) $S_\nu$ and $S_{\nu+\mu}$ are completely determined by $Y_{\nu}(z)$, $Y_{\nu+\mu}(z)$  and $\Lambda^\prime$, because
\be
\label{sole}
    Y_{\nu+\mu}(z)  = Y_\nu(z) S_\nu,~~~~~z\in {\cal S}_\nu\cap {\cal S}_{\nu+\mu},
   \ee
\be
\label{le12}
   Y_\nu(z e^{-2\pi i})  = Y_{\nu+\mu}(z) S_{\nu+\mu}e^{-2\pi i \Lambda^\prime},~~~~~z\in {\cal S}_{\nu+\mu}\cap {\cal S}_{\nu +m}.
  \ee

\noindent
\noindent
iii) $Y_{\nu+m}(z)$ is determined by $Y_\nu(z)$ and $\Lambda^\prime$, because
\be
\label{le123}
 Y_{\nu+m}(ze^{2\pi i})= Y_\nu (z)e^{2\pi i \Lambda^\prime},~~~~~z\in {\cal S}_\nu.
\ee

\ere

   \bre 
   In \cite{Dub1}, \cite{Dub2}, \cite{Dub3}, and oriented {\it admissible line $l$}, namely a line not containing Stokes rays, is introduced to define the analytic continuation of  $Y_\nu$ and $Y_{\nu +\mu}$. The branch cut is taken to be the negative part of $l$. We can take $l=l(\eta)$ such that its positive part is the ray from 0 to $\infty$ with angle $\tau={3\pi\over 2}-\eta$. 
   
   \ere

\bre If $A_1$ has some symmetries, then there  may be a relation between $S_\nu$ and $S_{\nu+\mu}$. For example, in \cite{Dub1}, \cite{Dub2} the  case $
   A_1^T=-A_1
   $ is considered (where $T$ means transposition). 
    This  implies that
    $$
   S_{\nu+\mu}^T=S_{\nu}^{-1}.
   $$  
   \ere


\subsection{Solutions of (\ref{01}) as Laplace Integrals}

\begin{figure}
\centerline{\includegraphics[width=0.3\textwidth]{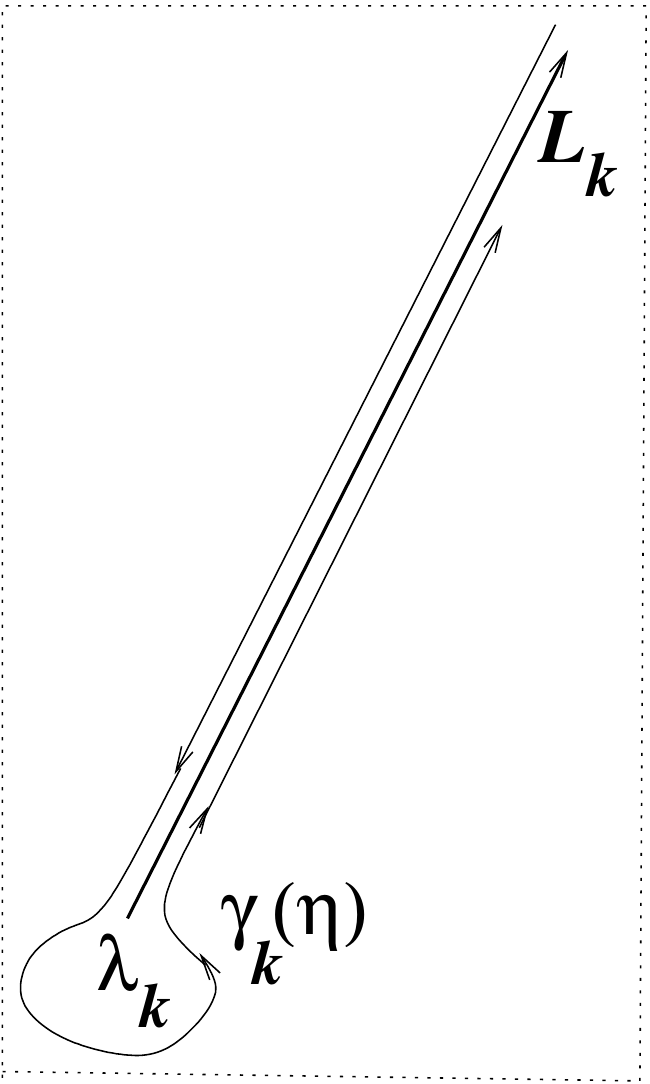}}
\caption{The path $\gamma_k(\eta)$
}
\label{fig5}
\end{figure}

We consider a path $\gamma_k(\eta)$ which comes from infinity along the left side of the cut  $L_k$ of direction $\eta$, encircles $\lambda_k$ with a small loop excluding all the other poles, and goes back to infinity along the right side of $L_k$ (where $L_k$ is oriented from  $\lambda_k$ to $\infty$). See figure \ref{fig5}.

\vskip 0.2 cm 
 {\bf 1) Case of $\lp_k\not \in \mathbb{Z}$.} We define 
\be
\label{chocho1}
\sqs{
 \vec{Y}_k(z,\eta):={1\over 2\pi i } \int_{\gamma_k(\eta)}e^{z\lambda}~\vPsi_k(\lambda,\eta)~ d\lambda~\equiv {1\over 2\pi i } \int_{\gamma_k(\eta)}e^{z\lambda}~\vPsi_k^{(k)}(\lambda,\eta)~ d\lambda
}
\ee

Due to the fact that $\lambda_k$ is a regular singularity, the exponential ensures that  the integral converges in the sector 
\be
\label{seccho1}
{\cal S}(\eta):=\left\{
z\in \widetilde{\mathbb{C}\backslash\{0\}}~|~\Re(ze^{i\eta})<0
\right\}
~~~\Longrightarrow~~~
{\pi\over 2}-\eta<\arg z <{3\pi \over 2}-\eta.
\ee

 The asymptotic behaviour of (\ref{chocho1}) can be computed by expanding the integrand $\vPsi_k(\lambda)$ in series at $\lambda_k$ and then formally exchanging integration and series (see \cite{Doe}). 
 Namely, for any ${\cal N}>0$ integer, 
$$
 \vec{Y}_k(z,\eta)= {1\over 2\pi i }  \int_{\gamma_k(\eta)}e^{z\lambda}~
\left[
\Gamma(\lp_k+1)\vec{e}_k+\sum_{l\geq 1} \vec{b}_l^{~(k)}(\lambda-\lambda_k)^l
\right](\lambda-\lambda_k)^{-\lp_k-1}~d\lambda=(*)
$$
Now write $\sum_{l\geq 1}=\sum_{l=1}^{\cal N}+\sum_{l> {\cal N}}$, and use the formula (see \cite{Doe}):
$$
 \int_{\gamma_k(\eta)}(\lambda-\lambda_k)^{-a} e^{z\lambda}d\lambda=
{z^{a-1} e^{\lambda_k z}\over \Gamma(a)}. 
$$  We obtain
\be 
(*)=\left(\vec{e}_k+\sum_{l=1}^{\cal N}{\vec{b}_l^{~(k)}\over \Gamma(\lp_k+1-l)}~{1\over z^l}+{\cal R}(z)\right)~e^{\lambda_k z}z^{\lp_k},
\label{Rz}
\ee
where ${\cal R}(z)$ is the integral of $\sum_{l> {\cal N}}$. 
It is standard computation to show that ${\cal R}(z)=O(z^{\cal N})$.  Thus, formula (\ref{Rz}) allows us to write the asymptotic expansion
$$
\vec{Y}_k(z,\eta)~e^{-\lambda_k z}z^{-\lp_k}\sim \vec{e}_k+\sum_{l=1}^{\infty}{\vec{b}_l^{~(k)}\over \Gamma(\lp_k+1-l)}~{1\over z^l},~~~~~~~~z\to \infty,~~z\in {\cal S}(\eta).
$$

\ble 
\label{lem3}
Assume $\lp_k\not\in \mathbb{Z}$. Let $\eta\in(\eta_{\nu+1},\eta_\nu)$, and  $\tau_\nu:=3\pi/2-\eta_\nu$. Then, $\vec{Y}_k(z,\eta)$ defined by  (\ref{chocho1}) is the $k$-th column of the unique fundamental solution of (\ref{01}) identified by the asymptotic behavior (\ref{funda}), (\ref{asint}) in the sector 
$$
 {\cal S}_\nu=S(\tau_\nu-\pi,\tau_{\nu+1}).
$$
\ele

\vskip 0.2 cm 
\noindent
{\it Proof:} If $\eta_{\nu+1}<\eta<\tilde{\eta}<\eta_\nu$, then $Y_k(z,\eta)=Y_k(z,\tilde{\eta})$. This defines the analytic continuation of (\ref{chocho1}) to 
$$
S(\tau_\nu-\pi,\tau_{\nu+1})=\bigcup_{\eta_{\nu+1}<\eta<\eta_\nu}{\cal S}(\eta),
$$
with the required asymptotic behaviour. It remains to prove that $Y_k(z,\eta)$ is a vector solution of (\ref{01}). This follows from  integration by parts, as shown in the Introduction, since  $\gamma_k(\eta)$ is such that $e^{\lambda z}(\lambda-A_0)\vPsi_k(\lambda)\Bigl|_{\gamma_k}=0$.
$\Box$

\vskip 0.3 cm 
We write coefficients in (\ref{asint}) as 
$$
F_k=[~\vec{f}_1^{~(k)}~|~\cdots~|~\vec{f}_n^{~(k)}~]
$$
 Thus, we can write (\ref{asso1}) for $\eta_{\nu+1}<\eta<\eta_\nu $ as 
$$
\sqs{
 \vPsi_k(\lambda,\eta)\equiv \vPsi_k^{(k)}(\lambda)= 
\sum_{l\geq 0} \Gamma(\lp_k+1-l)~\vec{f}_l^{~(k)}~(\lambda-\lambda_k)^{l-\lp_k-1},
~~~~~~~
\vec{f}_0^{~(k)}=\vec{e}_k
}
$$

\vskip 0.3 cm 

{\bf 2) Case  of $\lp_k=-1$.} We define 
\be
\label{chocho2}
 \vec{Y}_k(z,\eta):= \int_{L_k}e^{z\lambda}~\vPsi_k(\lambda,\eta)~ d\lambda~=-\int_{-L_k}e^{z\lambda}~\vPsi_k(\lambda,\eta)~ d\lambda
\ee
along the cut $L_k$  from $\lambda_k$ to infinity. This is convergent in ${\cal S}(\eta)$ as before. Its asymptotic behaviour is obtained as before by expanding $\vPsi_k$ in the convergent series (\ref{asso12}), and then exchanging integration and series, the result having meaning of asymptotic series. We obtain, by elementary integration:
$$
Y_k(z,\eta)= \vec{e}_k\int_{-L_k}e^{z\lambda}d\lambda- \sum_{l\geq 1} \vec{b}_l^{~(k)}  \int_{-L_k}(\lambda-\lambda_k)^l 
e^{z\lambda}~d\lambda
=
$$
$$
=
{e^{\lambda_k z}\over z} \left[\vec{e}_k+\sum_{l=1}^\infty (-1)^{l+1}~ l!~ ~\vec{b}_l^{~(k)}{1\over z^l}\right].
$$
where we have used the fact that 
$$
\int_{-L_k}(\lambda-\lambda_k)^l e^{\lambda z} d\lambda= {e^{\lambda_k z}\over z^{l+1}}\int_{+\infty e^{i\phi}}^0 \xi^l e^\xi d\xi={e^{\lambda_k z}\over z^{l+1}}~l!~(-1)^l ,~~~~~{\pi \over 2} <\phi<{3\pi \over 2}.
$$
The same proof of Lemma \ref{lem3} yields the following
\ble 
Assume $\lp_k=-1$. Let $\eta\in(\eta_{\nu+1},\eta_\nu)$, and  $\tau_\nu:=3\pi/2-\eta_\nu$. Then, $Y_k(z,\eta)$ defined by  (\ref{chocho2}) is the $k$-th column of the unique fundamental solution of (\ref{01}) identified by the asymptotic behavior (\ref{funda}), (\ref{asint}) in the sector 
$$
 {\cal S}_\nu=S(\tau_\nu-\pi,\tau_{\nu+1}).
$$
\ele

\noindent
By virtue of the lemma, we rewrite (\ref{asso12}),  for $\eta_{\nu+1}<\eta<\eta_\nu$, as follows
$$
\sqs{
 \vPsi_k(\lambda,\eta)=
-\vec{e}_k+\sum_{l\geq 1} {(-1)^{l+1}\over l!} \vec{f}_l^{~(k)}(\lambda-\lambda_k)^l
}
$$

\ble

\vskip 0.3 cm 
In case $\lp_k=-1$, the solution (\ref{chocho2}) has also the representation
\be
\label{chocho23}
\sqs{
 \vec{Y}_k(z,\eta)= \int_{L_k}e^{z\lambda}~\vPsi_k(\lambda,\eta)~ d\lambda~
\equiv
{1\over 2\pi i} \int_{\gamma_k(\eta)} \vPsi_k^{(k)}(\lambda) e^{z\lambda}~d\lambda
}
\ee
where $\gamma_k(\eta)$ is the same of (\ref{chocho1}). 
\ele

\vskip 0.2 cm 
\noindent
{\it Proof:} Recall that $\vPsi_k^{(k)}= \vPsi_k(\lambda)\ln(\lambda-\lambda_k)+ \hbox{reg}(\lambda-\lambda_k)$. Since 
$$
\int_{\gamma_k} \hbox{reg}(\lambda-\lambda_k) ~e^{z\lambda}~d\lambda=0
$$
we have 
$$
\int_{\gamma_k(\eta)} \vPsi_k^{(k)}(\lambda) ~e^{z\lambda}~d\lambda=\int_{\gamma_k(\eta)} \vPsi_k(\lambda)\ln(\lambda-\lambda_k)~e^{z\lambda}~d\lambda
$$
Indicate with $L_k^L$ and $L_k^R$ the left and right sides of $L_k$ (oriented from $\lambda_k$ to $\infty$), and with $(\lambda-\lambda_k)_{R/L}$ the branch of $(\lambda-\lambda_k)$ to the right/left of $L_k$. Then 
$$
 \int_{\gamma_k(\eta)} =\int_{-L_k^L}~+\int_{L_k^R}~= \int_{L_k^R}~-\int_{L_k^L}~=(*)
$$
Moreover 
$$ 
(\lambda-\lambda_k)_L=e^{-2\pi i} (\lambda-\lambda_k)_R ,~~~\hbox{ where } \arg( (\lambda-\lambda_k)_R)=\eta.
$$
Therefore
$$
(*)=
\int_{L_k^R}\vPsi_k(\lambda)\ln(\lambda-\lambda_k)_R~e^{z\lambda}~d\lambda+\left\{2\pi i   \int_{L_k}\vPsi_k(\lambda) ~e^{z\lambda}~d\lambda  -\int_{L_k^R}\vPsi_k(\lambda)\ln(\lambda-\lambda_k)_R~e^{z\lambda}~d\lambda\right\}
$$
$$
\equiv 2\pi i   \int_{L_k}\vPsi_k(\lambda) ~e^{z\lambda}~d\lambda
.
$$
$\Box$

\vskip 0.3 cm 
{\bf 3) Case  of $\lp_k\in \mathbb{N}$. }  
Define the convergent in ${\cal S}(\eta)$ integral 
\be
\label{chocho3}
\sqs{
 \vec{Y}_k(z,\eta):={1\over 2\pi i } \int_{\gamma_k(\eta)}e^{z\lambda}~\vPsi_k^{(k)}(\lambda,\eta)~ d\lambda
}
\ee
where 
$$
\vPsi^{(k)}_k(\lambda)=\vPsi_k(\lambda)\ln(\lambda-\lambda_k)+{ P_{N_k}^{(k)}(\lambda)
\over (\lambda-\lambda_k)^{N_k+1}}
+\hbox{reg}(\lambda-\lambda_k),
$$
$$
  \vPsi_k=\sum_{l\geq 0} {\vec{d}_l}^{~(k)}(\lambda-\lambda_k)^l,~~~
P_{N_k}^{(k)}=N_k!~\vec{e}_k+\sum_{l=0}^{N_k}\vec{b}_l^{~(k)}(\lambda-\lambda_k)^l,
$$
As before, we prove that 
$$
{1\over 2\pi i } \int_{\gamma_k(\eta)}e^{z\lambda}  \vPsi_k(\lambda)\ln(\lambda-\lambda_k)d\lambda
=\int_{L_k} e^{z\lambda} \vPsi_k(\lambda) d\lambda
$$
$$=\sum_{l=0}^\infty (-1)^{l+1}~l!~d_l^{~(k)}~{1\over z^{l+1}}~e^{z\lambda_k}.
$$
where "$=$"  means asymptotic for $z\to \infty$.  
On the other hand, by Cauchy theorem,  
$$
{1\over 2\pi i } \int_{\gamma_k(\eta)}e^{z\lambda} { P^{(k)}(\lambda)
\over (\lambda-\lambda_k)^{N_k+1}}
=
{1\over N_k!} {d^{N_k}\over d\lambda^{N_k}}\left.
\left( P^{(k)}(\lambda)e^{z\lambda}\right)\right|_{\lambda=\lambda_k}
$$
$$
= e^{\lambda_k z}~ \sum_{q=0}^{N_k}{\vec{b}_{N_k-q}^{~(k)}\over q!} ~z^q= \left[\vec{e}_k +\cdots+{\vec{b}_{N_k}^{(k)}\over z^{N_k}}\right]z^{N_k}e^{\lambda_k z},~~~~~\vec{b}_{0}^{~(k)}=N_k!~\vec{e}_k.
$$
We conclude that 
$
\vec{Y}_k(z,\eta)$ has the correct asymptotics. The same proof of Lemma \ref{lem3} yields the following

\ble 
Assume $\lp_k\in \mathbb{N}$. Let $\eta\in(\eta_{\nu+1},\eta_\nu)$, and  $\tau_\nu:=3\pi/2-\eta_\nu$. Then, $\vec{Y}_k(z,\eta)$ defined by  (\ref{chocho3}) is the $k$-th column of the unique fundamental solution of (\ref{01}) identified by the asymptotic behavior (\ref{funda}), (\ref{asint}) in the sector 
$$
 {\cal S}_\nu=S(\tau_\nu-\pi,\tau_{\nu+1}).
$$
\ele
Accordingly, 
 we rewrite for $\eta_{\nu+1}<\eta<\eta_\nu$:
$$
\sqs{\vPsi_k^{(k)}(\lambda,\eta)= \sum_{l=0}^{N_k}(N_k-l)!~\vec{f}_l^{~(k)}~u^{l-N_k-1}
+\left[
\sum_{l=0}^\infty {(-1)^{l+1}~\vec{f}_{N_k+l+1}^{~(k)}\over l!} ~u^l
\right]
\ln(u)+\hbox{reg}(u)}
$$
where $u:=\lambda-\lambda_k$.

\vskip 0.3 cm

{\bf 4) case of $\lp_k=N_k\in -\mathbb{N}-2$.} We define 
\be
\label{chocho9}
\sqs{\vec{Y}_k(z,\eta):=\int_{L_k}e^{z\lambda} \vPsi_k(\lambda) d\lambda}
\ee
The asymptotic behaviour of the above is readily computed:
$$
\int_{L_k}e^{z\lambda} \vPsi_k(\lambda) d\lambda=\int_{L_k}e^{z\lambda} \sum_{l\geq 0} \vec{b}_l^{(k)}(\lambda-\lambda_k)^{l-N_k-1} d\lambda
$$
$$
= e^{\lambda_k z}(-1)^{N_k}\sum_{l\geq 0}{(-1)^l (l-N_k-1)! ~\vec{b}_l^{(k)} \over z^{l-N_k}}
$$
$$
= \left[
\vec{e}_k +(-1)^{N_k}\sum_{l\geq 1} (-1)^l (l-N_k-1)! ~\vec{b}_l^{(k)} ~{1\over z^{l}}\right]z^{N_k}e^{\lambda_k z}
.
$$ 
where we have used the normalization $\vec{b}_0^{(k)}=(-1)^{N_k}\vec{e}_k /(-N_k-1)!$. 
We conclude that $\vec{Y}_k(z,\eta)$ has the correct asymptotics. The same proof of Lemma \ref{lem3} yields the following

\ble 
Assume $\lp_k=N_k\in -\mathbb{N}-2$. Let $\eta\in(\eta_{\nu+1},\eta_\nu)$, and  $\tau_\nu:=3\pi/2-\eta_\nu$. Then, $\vec{Y}_k(z,\eta)$ defined by  (\ref{chocho9}) is the $k$-th column of the unique fundamental solution of (\ref{01}) identified by the asymptotic behaviour (\ref{funda}), (\ref{asint}) in the sector 
$$
 {\cal S}_\nu=S(\tau_\nu-\pi,\tau_{\nu+1}).
$$
\ele
Accordingly:
$$
\sqs{
\vPsi_k(\lambda)= 
\sum_{l\geq 0} {(-1)^{l-N_k}\over (l-N_k-1)!}~ \vec{f}_l^{(k)}~(\lambda-\lambda_k)^{l-N_k-1}
}
$$
Also in this case we have
$$
{1\over 2\pi i } \int_{\gamma_k(\eta)}e^{z\lambda}  \vPsi_k(\lambda)\ln(\lambda-\lambda_k)d\lambda
=\int_{L_k} e^{z\lambda} \vPsi_k(\lambda) d\lambda
$$
Therefore, {\it when the singular solution $\vPsi_k(\lambda)\ln(\lambda-\lambda_k)+\hbox{reg}(\lambda-\lambda_k)$ exists}, we have
$$
\sqs{
 \vec{Y}_k(z,\eta)=\int_{\gamma_k(\eta)}e^{z\lambda} \Bigl(\vPsi_k(\lambda)\ln(\lambda-\lambda_k)+\hbox{reg}(\lambda-\lambda_k)\Bigr)d\lambda.
 }
 $$
\vskip 0.3 cm 
\bpr
\label{FFUND}
The following are the fundamental matrix solutions of (\ref{01}) uniquely identified by the asymptotic behaviour (\ref{funda}),  (\ref{asint}) in ${\cal S}_\nu$, $\nu\in \mathbb{Z}$:
$$
 Y_\nu(z)=\left[ ~\vY_1(z,\eta)~|~\cdots~|~ \vY_n(z,\eta)~\right]
 $$
 \be
 \label{chocho4}
 \vY_k(z,\eta)={1\over 2\pi i } \int_{\gamma_k(\eta)}e^{z\lambda}~\vPsi_k^{(sing)}(\lambda,\eta)~ d\lambda,~~~1\leq k \leq n,~~~~\eta_\nu<\eta<\eta_{\nu+1},
 \ee
 where $\vPsi_k^{(sing)}$ is defined in (\ref{GNIS}). In case it happens  that  $\vPsi^{(sing)}=0$,  $\lp_k\in -\mathbb{N}-2$, then (\ref{chocho4}) is replaced by (\ref{chocho9}).
 \epr 
 {\it Proof:} The above is a consequence of the preceding discussion. Linear independence of the columns of $Y_\nu(z)$ follows from the independence of the first term of the asymptotic behaviour of each column. Uniqueness follows from the maximality of the sector. $\Box$
 
 \vskip 0.3 cm 
 
\ble If $A_1$ has no negative integer eigenvalues, then, $\vPsi_k^{(sing)}$ in the integral (\ref{chocho4}) can be replaced  by $\vPsi_k^*$ . Namely:
$$
 \vY_k(z,\eta)={1\over 2\pi i } \int_{\gamma_k(\eta)}e^{z\lambda}~\vPsi_k^{*}(\lambda,\eta)~ d\lambda,~~~1\leq k \leq n,~~~~\eta_\nu<\eta<\eta_{\nu+1}.
$$
\ele
\vskip 0.2 cm 
\noindent
{\it Proof:} Recall that if $A_1$ has no negative integer eigenvalues, then $\vPsi^{(sing)}\neq 0$ also for $\lp_k\in -\mathbb{N}-2$. Since $\vPsi_k^{(sing)}-\vPsi_k^*=\hbox{reg}(\lambda-\lambda_k)$, we have
$$
\int_{\gamma_k(\eta)} \left( \vPsi_k^{(sing)}(\lambda)-\vPsi_k^*(\lambda)\right) e^{z\lambda}d\lambda=0.
$$
$\Box$


\section{Stokes Factors and Matrices in terms of $C$ - Main Theorem (Th. \ref{RelationSC2})}
\label{MT}

In this section we state the main result of the paper,  which is  Theorem \ref{RelationSC2} and Corollary \ref{CHICH}.

Consider as in \cite{BJL} a new path of integration $\gamma(\eta)$ homotopic to the product  $\gamma_{k_n}(\eta)\cdots\gamma_{k_1}(\eta)$, $k_1\prec  k_2\prec ... \prec k_n$,  namely a path  coming from $\infty$ in direction $\eta$ to the left of all the poles $\lambda_1,...,\lambda_n$, encircling all the poles, and going back to $\infty$ in direction $\eta$ to the right of all the poles. 
The following Proposition is the generalization of Theorem $2^\prime$ of \cite{BJL}  when no assumptions are made on  diag$(A_1)$. 
\bpr
\label{RelationSC1}
 If $A_1$ has no negative integer eigenvalues the fundamental matrix of Proposition \ref{FFUND} is 
$$
 Y_\nu(z)={1\over 2\pi i} \int_{\gamma(\eta)}~e^{z\lambda}~\Psi^*(\lambda,\eta)~d\lambda,~~~~~\eta_{\nu+1}<\eta<\eta_\nu.
 $$
 and 
 $$
  Y_{\nu-1}(z)=Y_\nu(z) W_\nu,~~~~~z\in{\cal S}_{\nu-1}\cap{\cal S}_\nu = S(\tau_\nu-\pi,\tau_\nu).
  $$
  where the $W_\nu$'s are given in Proposition \ref{propoW}. Therefore, the Stokes factors are 
  $$
  V_\nu = W_\nu
  $$
    Moreover,
  $$
  S_\nu=C_\nu^+,~~~~~S_{\nu+m}^{-1}=C_\nu^-.
  $$
\epr

\begin{figure}
\centerline{\includegraphics[width=0.5\textwidth]{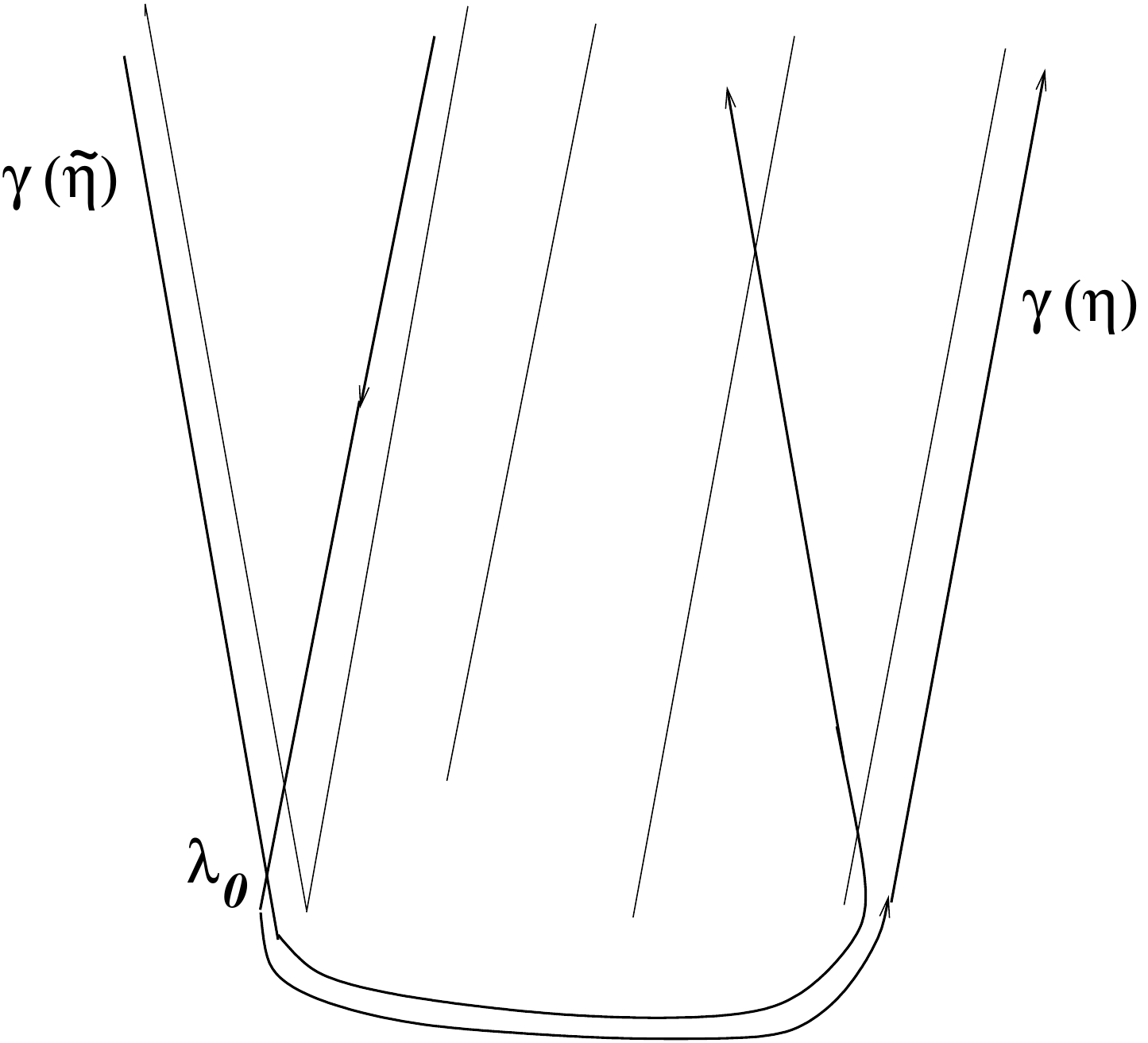}}
\caption{The paths $\gamma(\eta)$ and $\gamma(\tilde{\eta})$.}
\label{fig6}
\end{figure}

\vskip 0.2 cm 
\noindent
{\it Proof:}  The first statement is easy. Indeed 
$\Psi_k^*(\lambda)=\hbox{reg}(\lambda-\lambda_j)$ for any $j\neq k$ implies
$$
\int_{\gamma_k(\eta)}~e^{z\lambda}~\vPsi_k^*(\lambda,\eta)~d\lambda
=
\int_{\gamma(\eta)}~e^{z\lambda}~\vPsi_k^*(\lambda,\eta)~d\lambda
$$
Then,we  consider $\eta_{\nu+1}<\eta<\eta_\nu<\tilde{\eta}<\eta_{\nu-1}$ and 
$$
2\pi i Y_{\nu-1}(z)= \int_{\gamma(\tilde{\eta})}~e^{z\lambda}~\Psi^*(\lambda,\tilde{\eta})~d\lambda=(*)
$$
We can define the analytic continuation of $\Psi^*(\lambda,\tilde{\eta})$ along $\gamma(\eta)$ as follows. We consider on  $\gamma(\tilde{\eta})$ a reference point $\lambda_0$ w.r.t. both $\tilde{\eta}$ and $\eta$, and deform $\gamma(\tilde{\eta})$ by keeping $\lambda_0$ fixed, until we obtain $\gamma(\eta)$. See fifure \ref{fig6}.  Since $\lambda_0$ is a reference point, $\Psi^*(\lambda_0,\tilde{\eta})=\Psi^*(\lambda_0,\eta)$. Thus  the analytic continuation of $\Psi^*(\lambda,\tilde{\eta})$ along $\gamma(\eta)$ is $\Psi^*(\lambda,\eta)W_\nu$. Consequently
$$
(*)=  \int_{\gamma(\eta)}~e^{z\lambda}~\left(\Psi^*(\lambda,\eta)W_\nu\right)~d\lambda= \left(\int_{\gamma(\eta)}~e^{z\lambda}~\Psi^*(\lambda,\eta)~d\lambda\right)~W_\nu
~\equiv 
2\pi i Y_\nu(z) W_\nu.
$$
The last statement follows from (\ref{ciclo}) and (\ref{ciclo1}).
$\Box$

\vskip 0.2 cm

We are going to  prove that the statement of the above Proposition holds also {\it for any $A_1$, without assumptions}. This result is the following 
%

\bth
\label{RelationSC2} Let  $A_1$ be any $n\times n$ matrix, with no assumptions.
The Stokes multipliers and matrices of system (\ref{01})  are given in terms of the connection coefficients $c_{jk}^{(\nu)}$ of system (\ref{02}) according to the formulae
$$
 \sqs{
 V_\nu=W_\nu, ~~~~~~S_\nu=C_\nu^+,~~~~~S_{\nu+m}^{-1}=C_\nu^-
,~~~~~\forall \nu\in \mathbb{Z} }
$$
 where $W_\nu$ {\rm is defined} by formulae (\ref{BEC1}), (\ref{BEC2}), and  $C_\nu^+$ and $C_\nu^-$  {\rm are defined} by formulae (\ref{BEC3}) and  (\ref{BEC4}).
\eth
{\bf Remark:} Here  formulae (\ref{BEC1}), (\ref{BEC2}), (\ref{BEC3}) and  (\ref{BEC4}) {\it are taken  as  the definitions} of $W_\nu$, $C_\nu^+$ and $C_\nu^-$, independently of the existence of $\Psi^*(\lambda)$. 
%
%


\bcr
\label{CHICH}
 Let  $A_1$ be any $n\times n$ matrix, with no assumptions. The following equalities hold for the monodromy matrices of $\Psi(\lambda)$ of system (\ref{02})-(\ref{03}), defined in (\ref{Pl}):
$$
\sqs{ \hbox{\rm Tr}(M_k)=n-1+e^{-2\pi i \lp_k}
}
$$
$$
\sqs{ \hbox{\rm Tr}(M_jM_k)=
\left\{
\matrix{
n-2 +e^{-2\pi i \lp_j}+e^{-2\pi i \lp_k}-e^{-2\pi i \lp_j}~[S_\nu]_{jk}[S_{\nu+m}^{-1}]_{kj} & \hbox{ if } j\prec k,
\cr
\cr
n-2 +e^{-2\pi i \lp_j}+e^{-2\pi i \lp_k}-e^{-2\pi i \lp_k}~[S_{\nu+m}^{-1}]_{jk}[S_\nu]_{kj} & \hbox{ if } j\succ k.
}
\right.
}
$$
\ecr
The corollary above is a restatement of Corollary \ref{mana}. We prove Theorem \ref{RelationSC2} in a few steps.

\section{Proof of Theorem \ref{RelationSC2}}
\label{THEPROOF}

We define 
$${}_\gamma Y(z):=z^{-\gamma} Y(z),
$$
 which yields a gauge transformation of the linear systems (\ref{01}): 
\be 
\label{01bis}
\frac{d}{dz}~({}_\gamma Y)=\left( A_0+{A_1-\gamma\over z} \right)~ {}_\gamma Y
\ee
The fundamental solutions ${}_\gamma Y_\nu(z)=z^{-\gamma} Y_\nu(z)$, have the same Stokes multiplier and Stokes matrices than $Y_\nu(z)$, and their columns are obtained as Laplace transforms of solutions of 
\be
\label{02bis}
(A_0-\lambda){d\over d\lambda} (~{}_\gamma\Psi)=  (A_1-\gamma+I)~{}_\gamma\Psi.
\ee
 If $A_1$ has diagonal entries $\lp_1,...,\lp_n$, some of which may be integers, then we can always find a sufficiently small $\gamma_0>0$ such that, for any $0<\gamma<\gamma_0$, $A_1-\gamma$ has diagonal entries $\lp_1-\gamma,...,\lp_n-\gamma$ which are not integers, and moreover has no integer eigenvalues, so that $\Psi^*$ exists. 
 In the following, we assume that $\gamma$ has this property. 
 
For system (\ref{02}), the matrix $C_\nu=(c_{jk}^{(\nu)})$ is defined by  (\ref{psik}). Consequently, the matrices   $C_\nu^+$ and $C_\nu^-$ are {\it always defined}  by the  formulae of  Proposition \ref{cpmPro}, independently of the existence of $\Psi^*$ and formulae (\ref{cpiu}) and (\ref{cmeno}).  On the other hand, for system (\ref{02bis}) the  matrices $C_\nu^+$ and $C_\nu^-$ (which depend on $\gamma$, so we write  $C_\nu^+[\gamma]$ and $C_\nu^-[\gamma]$), are well defined by formulae (\ref{cpiu}) and (\ref{cmeno}). According to Proposition \ref{cpmPro}, their entries are again given in terms of $\gamma$-dependent connection coefficients  $c_{jk}^{(\nu)}=c_{jk}^{(\nu)}[\gamma]$'s. The latter are defined  by the first equality  of (\ref{psik}) applied to the  solutions $~{}_\gamma\vPsi_k$, namely:
 \be 
 \label{ine13}
 {}_\gamma\vPsi_k(\lambda)=~{}_\gamma\vPsi_j(\lambda) ~c_{jk}^{(\nu)}[\gamma]~ +\hbox{reg}(\lambda-\lambda_j)
 \ee
 The following Proposition is the key step to prove Theorem \ref{RelationSC2}

\bpr 
\label{fundlpro}
Let $\gamma_0>0$ be small enough such that the diagonal part of $A_1-\gamma I$ has no integer entries and $A_1$ has no integer eigenvalues for any $0<\gamma<\gamma_0$. 
Let $\eta_{\nu+1}<\eta<\eta_\nu$ be fixed. Let $c_{jk}^{(\nu)}$ be the corresponding  connection coefficients  of system (\ref{02}), defined by (\ref{psik}), and  $c_{jk}^{(\nu)}[\gamma]$ be the connection coefficients of (\ref{02bis}), defined by (\ref{ine13}). Finally, let 
$$
\alpha_k=\left\{\matrix{
e^{-2\pi i \lp_k}-1,& \lp_k\not\in\mathbb{Z}
\cr
\cr
2\pi i, & \lp_k\in\mathbb{Z}
}\right.;
~~~~~~~~\alpha_k[\gamma]=e^{-2\pi i (\lp_k-\gamma)}-1
$$ 
Then, the following equalities hold
$$
\alpha_k c_{jk}^{(\nu)}= 
e^{-2\pi i \gamma}\alpha_k[\gamma]~c_{jk}^{(\nu)}[\gamma],~~~ \hbox{ if }k\succ j
$$
$$
\alpha_k c_{jk}^{(\nu)}= \alpha_k[\gamma]~c_{jk}^{(\nu)}[\gamma] ,~~~~~~~~~~~ \hbox{ if } k\prec j
$$
where the partial ordering $\prec$ refers to $\eta$. 
\epr

\bcr
\label{clcl}
Let $\gamma$ be as in Proposition \ref{fundlpro}. 
Let  $C_\nu^+[\gamma]$ and $C_\nu^-[\gamma]$ be the connection matrices defined in (\ref{cpiu}) and (\ref{cmeno}) for system (\ref{02bis}). 
 Let $C_\nu^+$ and $C_\nu^-$ be the  matrices   for system (\ref{02}) defined by (\ref{BEC3}) and (\ref{BEC4}), where the $c_{jk}^{(\nu)}$ are defined by (\ref{psik}).   Then
$$
 C_\nu^+=C_\nu^+[\gamma],~~~~~C_\nu^-=C_\nu^-[\gamma],~~~~~~~~\forall \nu\in \mathbb{Z}. 
 $$
Also, let $W_\nu$ be defined by (\ref{BEC1}) and (\ref{BEC2}) for system (\ref{02}), and $W_\nu[\gamma]$ be the matrix defined by (\ref{BEC0}) for system (\ref{02bis}). Then 
$$ 
W_\nu=W_\nu[\gamma],~~~~~~~~\forall \nu\in \mathbb{Z}.
$$
\ecr

\vskip 0.2 cm 
\noindent
{\it Proof of Corollary \ref{clcl}:} It is enough to compare the formulae of  Proposition  \ref{fundlpro} with those of 
 Propositions \ref{cpmPro} and \ref{propoW}. $\Box$

\vskip 0.3 cm 
Before proving Proposition \ref{fundlpro}, we give the proof of Theorem \ref{RelationSC2}.

\vskip 0.2 cm 
\noindent
{\it Proof of Theorem \ref{RelationSC2}:} The $S_\nu$'s are unchanged by the gauge ${}_\gamma Y(z)=z^{-\gamma} Y(z)$. Moreover, Proposition \ref{RelationSC1} applies to the system (\ref{02bis}), therefore 
$$ 
 S_\nu=C_\nu^+[\gamma],~~~~~S_{\nu+m}^{-1}=C_\nu^-[\gamma],~~~~~V_\nu=W_\nu[\gamma].
 $$
Thus, Corollary \ref{clcl} implies Theorem \ref{RelationSC2}. $\Box$.

\subsection{Proof of Proposition \ref{fundlpro}, by steps}

The idea is the same of the proof of Lemma $2^\prime$ in \cite{BJL}, but with considerably more technical efforts, do to the fact that, unlike \cite{BJL}, we do not make any assumption on the diagonal entries of $A_1$. We need a few steps, which are Propositions \ref{undici} and \ref{lab11} , and Lemmas \ref{lab9}, \ref{lab10} below. First, we introduce the $q$-primitives of vector solutions of system (\ref{03}).

\vskip 0.3 cm 
-- For  $\lp_k\in\mathbb{C}\backslash\mathbb{N}$ we have solutions
 $$
 \vPsi_k(\lambda)=
\sum_{l=0}^\infty \Gamma(\lp_k+1-l)\vec{f}_l^{~(k)} ~(\lambda-\lambda_k)^{l-\lp_k-1},~~~~~\lp_k\not\in \mathbb{Z}
$$ 
and 
$$
\vPsi_k(\lambda)=
\sum_{l=0}^\infty {(-1)^{l-N_k}\over (l-N_k-1)!}~\vec{f}_l^{~(k)}~(\lambda-\lambda_k)^{l-N_k-1},~~~~~\lp_k=N_k\in \mathbb{Z}_{-}.
$$ 
 We define the {\it $q$ primitive of $\vPsi_k$}. This is the function $(\vPsi_k)^{[-q]}(\lambda)$, $q\in \mathbb{N}$,  given by analytic continuation of the series obtained by  $q$-fold term wise integration of the corresponding term in $\vPsi_k(\lambda)$. Namely: 
 \be 
 (\vPsi_k)^{[-q]}(\lambda):=(-1)^q~\sum_{l=q }^\infty \Gamma(\lp_k+1-l)\vec{f}_{l-q}^{~(k)} ~(\lambda-\lambda_k)^{l-\lp_k-1},~~~~~\lp_k\not\in \mathbb{Z}
\label{q1}
\ee

\be 
\label{q2}
(\vPsi_k)^{[-q]}(\lambda):=(-1)^q \sum_{l=q}^\infty {(-1)^{l-N_k}\over (l-N_k-1)!}~\vec{f}_{l-q}^{~(k)}~(\lambda-\lambda_k)^{l-N_k-1},~~~~~\lp_k=N_k\leq -1.
\ee
They above converge  in a neighbourhood of $\lambda_k$, contained in ${\cal P}_\eta$,      where $\vPsi_k$ has convergent series. 
Indeed, if $\lambda_0\neq \lambda_k$ is  in the neighbourhood, then   
\be
\label{cistit1}
\int_{\lambda_0}^\lambda ds_1\int_{\lambda_0}^{s_1}ds_2 \cdots\int_{\lambda_0}^{s_{q-1}}ds_q \vPsi_k(s_q)=(\vPsi_k)^{[-q]}(\lambda)-Q_{q-1}(\lambda-\lambda_0)
\ee
where 
$$Q_{q-1}(\lambda-\lambda_0)
=(\vPsi_k)^{[-q]}(\lambda_0)+(\vPsi_k)^{[1-q]}(\lambda_0)(\lambda-\lambda_0)+~~~~~~~~~~~~~~~~~~~
$$
$$
~~~~~~~~~~~~~~~~~~~~~~~~+{(\vPsi_k)^{[2-q]}(\lambda_0)
\over 2!}(\lambda-\lambda_0)^2+    \cdots   +{(\vPsi_k)^{[-1]}(\lambda_0)
\over (q-1)!}(\lambda-\lambda_0)^{q-1}
$$
 is a polynomial in $(\lambda-\lambda_0)$ of degree $q-1$.  The path of integration is any in ${\cal P}_\eta$, such that $|\lambda-\lambda_0|$ is small enough for the series of $\vPsi_k$ to converge. 
In particular,  for $\Re \lp_k<0$, 
 $$
 (\vPsi_k)^{[-q]}(\lambda)=\int_{\lambda_k}^\lambda ds_1\int_{\lambda_k}^{s_1}ds_2 \cdots\int_{\lambda_k}^{s_{q-1}}ds_q \vPsi_k(s_q)
$$
 Once $(\vPsi_k)^{[-q]}(\lambda)$ is defined by the  convergent series, then it is analytically continued to  ${\cal P}_\eta$.

 \vskip 0.3 cm 
 -- For $\lp_k=N_k\in\mathbb{N}$ integer, consider the solution
 $$
  \vPsi_k^{(k)}(\lambda)=\vPsi_k(\lambda)\ln(\lambda-\lambda_k)+{P_{N_k}^{(k)}(\lambda)\over (\lambda-\lambda_k)^{N_k+1}}+\hbox{reg}(\lambda-\lambda_k)
$$
$$
 P_{N_k}^{(k)}(\lambda) =\sum_{l=0}^{N_k}(N_k-l)!~\vec{f}_l^{~(k)}~(\lambda-\lambda_k)^l
$$
$$
\vPsi_k(\lambda)=
\sum_{l=0}^\infty {(-1)^{l+1}\over l!} ~\vec{f}_{N_k+l+1}^{~(k)}~(\lambda-\lambda_k)^l
$$
The series is convergent in a neighbourhood of $\lambda_k$ contained in ${\cal P}_\eta$. Let  $\lambda_0$ belong to  the neighbourhood.   Let $q\geq0$ integer, and compute $q$ times  the integral of $\vPsi_k^{(k)}(\lambda) $. 
Due to convergence of the series, we can take integration term by term. We obtain:

i) For $q\leq N_k$:
$$
\int_{\lambda_0}^\lambda ds_1\int_{\lambda_0}^{s_1}ds_2 \cdots\int_{\lambda_0}^{s_{q-1}}ds_q  \vPsi_k^{(k)}(s_q)=
$$
$$
=
\sum_{l=0}^\infty
{(-1)^{l+1-q}\over l!}~\vec{f}_{N_k+1+l-q}^{~(k)}~(\lambda-\lambda_k)^l~\ln(\lambda-\lambda_k)+{P_{N_k-q}^{(k)}(\lambda)\over (\lambda-\lambda_k)^{N_k+1-q}}
+\hbox{reg}(\lambda-\lambda_k)
$$
$$
P_{N_k-q}^{(k)}(\lambda)=(-1)^q\sum_{l=0}^{N_k-q}(N_k-l-q)!~\vec{f}_l^{~(k)}~(\lambda-\lambda_k)^l
$$

ii) For $q= N_k+1$:
$$
\int_{\lambda_0}^\lambda ds_1\int_{\lambda_0}^{s_1}ds_2 \cdots\int_{\lambda_0}^{s_{N_k}}ds_{N_k+1} \vPsi_k^{(k)}(s_{N_k+1})
=\hat{\Psi}_k(\lambda)
\ln(\lambda-\lambda_k)+\hbox{reg}(\lambda-\lambda_k)
$$
where we defined
\be 
\label{hatphi}
\sqs{
\hat{\Psi}_k(\lambda):=\sum_{l=0}^\infty{(-1)^{l+N_k}\over l!}~\vec{f}_l^{~(k)}~(\lambda-\lambda_k)^l}
\ee
The function  $\hat{\Psi}_k(\lambda)$ is defined by the series,  that converges in the neighbourhood of $\lambda_k$ in ${\cal P}_\eta$, where the series of $\vPsi_k^{(k)}$ converges. Then it is  analytically continued in ${\cal P}_\eta$. Note that if all $r_j^{(k)}=0$, $\forall j$, namely when there is no logarithmic term in  $\vPsi_k^{(k)}$, then the sum in  $\hat{\Psi}_k(\lambda)$ is truncated to $\sum_{l=0}^{N_k}$, giving a polynomial of degree $N_k$.

\vskip 0.2 cm 
iii) For $q=N_k+1+\tilde{q}$, with $\tilde{q}\geq 0$: 
$$
\int_{\lambda_0}^\lambda ds_1\int_{\lambda_0}^{s_1}ds_2 \cdots\int_{\lambda_0}^{s_{q-1}}ds_q  \vPsi_k^{(k)}(s_q)=(\hat{\Psi}_k)^{[-\tilde{q}]}(\lambda)\ln(\lambda-\lambda_k)+\hbox{reg}(\lambda-\lambda_k)
$$
$$
\sqs{
(\hat{\Psi}_k)^{[-\tilde{q}]}(\lambda)=
(-1)^{\tilde{q}}\sum_{l=\tilde{q}}^\infty {(-1)^{l+N_k}\over l!}~\vec{f}_{l-\tilde{q}}^{~(k)}~(\lambda-\lambda_k)^l
}
$$
 The function $(\hat{\Psi}_k)^{[-q]}(\lambda)$ is the {\it $q$ primitive} of $\hat{\Psi}_k(\lambda)$, and the same computation of (\ref{cistit1}) yields
\be
\label{cistit2}
\int_{\lambda_0}^\lambda ds_1\int_{\lambda_0}^{s_1}ds_2 \cdots\int_{\lambda_0}^{s_{q-1}}ds_q \hat{\Psi}_k(s_q)=(\hat{\Psi}_k)^{[-q]}(\lambda)-Q_{q-1}(\lambda-\lambda_0)
\ee
where  $Q_{q-1}$ is as in (\ref{cistit1}) with substitution $\vPsi_k\mapsto \hat{\Psi}_k$.  

\bre
The computation at point i) above can be also read as follows. Let 
$$
(\hat{\Psi}_k)^{[r]}(\lambda):={d^r\over d\lambda^r}\left(\hat{\Psi}(\lambda)\right)
$$
$$=
\sum_{l=0}^\infty {(-1)^{l+N_k-r}\over l!}~\vec{f}_{l+r}^{~(k)}~(\lambda-\lambda_k)^l,~~~~~0\leq r\leq N_k+1
$$
In particular 
$$ (\hat{\Psi}_k)^{[N_k+1]}(\lambda)=\vPsi_k(\lambda).
$$
Then 
$$
 {d^r\over d\lambda^r}\Bigl( \hat{\Psi}_k(\lambda)\ln(\lambda-\lambda_k)+\hbox{reg}(\lambda-\lambda_k) \Bigr)
=(\hat{\Psi}_k)^{[r]}(\lambda)\ln(\lambda-\lambda_k)+{P_{r-1}^{(k)}(\lambda)\over (\lambda-\lambda_k)^r} +\hbox{reg}(\lambda-\lambda_k) 
$$
with 
$$
P_{r-1}^{(k)}(\lambda)=(-1)^{N_k+1-r}\sum_{l=0}^{r-1}(r-1-l)!~\vec{f}_{l}^{~(k)}~(\lambda-\lambda_k)^l .
$$
\ere

\vskip 0.2 cm 
We summarize  (\ref{cistit1})  and  (\ref{cistit2})  and the computations involving logarithmic solutions in the following

\bpr
\label{undici}
 Let $\lambda_0\neq \lambda_j$ for any $j=1,2,...,n$. 
 
 $\diamond$ For a given $k\in\{1,...,n\}$ define 
$$
\phi_k(\lambda):=
\left\{
\matrix{
\vPsi_k(\lambda), & \hbox{ if } \lp_k\in \mathbb{C}\backslash \mathbb{N},
\cr 
\cr 
\hat{\Psi}_k(\lambda), & \hbox{ if } \lp_k\in \mathbb{N},
}
\right.
$$
\be
\label{cistit3}
\Phi_k^{[-q]}(\lambda,\lambda_0):=\int_{\lambda_0}^\lambda ds_1\int_{\lambda_0}^{s_1}ds_2 \cdots\int_{\lambda_0}^{s_{q-1}}ds_q \phi_k(s_q), 
\ee
where $\vPsi_k$ are defined in (\ref{psik0}) , and $\hat{\Psi}$ is defined in (\ref{hatphi}). 
Then 
\be
\label{cistit4}
\Phi_k^{[-q]}(\lambda,\lambda_0)=(\phi_k)^{[-q]}(\lambda)-Q_{q-1}^{(k)}(\lambda-\lambda_0)
\ee
where $Q_{q-1}^{(k)}$ is a polynomial of degree $q-1$ in $(\lambda-\lambda_0)$. It follows from the definition that  
\be
\label{recu}
\int_{\lambda_0}^\lambda ds~\Phi_k^{[-q]}(s,\lambda_0)~=~\Phi_k^{[-q-1]}(\lambda,\lambda_0)
\ee

 $\diamond$ For $\lp_k\in \mathbb{Z}$, consider the singular solutions $\vPsi_k^{(sing)}$ of system (\ref{02}):
$$
 \varepsilon_k \vPsi_k(\lambda)\ln(\lambda-\lambda_k)+\hbox{\rm reg}(\lambda-\lambda_k),~~~~~\lp_k\in \mathbb{Z}_{-},
$$
$$
\vPsi_k(\lambda)\ln(\lambda-\lambda_k)+{P^{(k)}(\lambda)\over (\lambda-\lambda_k)^{N_k+1}}+\hbox{\rm reg}(\lambda-\lambda_k),~~~~~\lp_k=N_k\in \mathbb{N}.
$$
In the above, $\varepsilon_k=0$ if ~$
\vPsi_k^{(sing)}\equiv 0$, otherwise $\varepsilon_k=1$. 
Then, for $\lp_k\in \mathbb{Z}_{-}$:
$$
\int_{\lambda_0}^\lambda ds_1\int_{\lambda_0}^{s_1}ds_2 \cdots\int_{\lambda_0}^{s_{q-1}}ds_q
 \Bigl(\varepsilon_k \vPsi_k(s_q)\ln(s_q-\lambda_k)+\hbox{\rm reg}(s_q-\lambda_k)\Bigr)=
$$
\be
\label{REE1}
=
\varepsilon_k(\vPsi_k)^{[-q]}(\lambda)\ln(\lambda-\lambda_k)+\hbox{\rm reg}(\lambda-\lambda_k),~~~~~q\geq 0,
\ee
and for $\lp_k\in \mathbb{N}$:
$$
\int_{\lambda_0}^\lambda ds_1\int_{\lambda_0}^{s_1}ds_2 \cdots\int_{\lambda_0}^{s_{q-1}}ds_q
\left(
\vPsi_k(s_q)\ln(s_q-\lambda_k)+{P_{N_k}^{(k)}(s_q)\over (s_q-\lambda_k)^{N_k+1}}+\hbox{reg}(s_q-\lambda_k)
\right)=
$$
\be
\label{REE2}
=
\left\{
\matrix{
(\hat{\Psi}_k)^{[N_k+1-q]}(\lambda)\ln(\lambda-\lambda_k)+{P_{N_k-q}^{(k)}(\lambda)\over (\lambda-\lambda_k)^{N_k+1-q}}+\hbox{\rm reg}(\lambda-\lambda_k) ~,& 0\leq q\leq N_k
\cr
\cr
\cr
(\hat{\Psi}_k)^{[-q+N_k+1]}(\lambda)\ln(\lambda-\lambda_k)+\hbox{\rm reg}(\lambda-\lambda_k) ~,& q\geq  N_k+1
}
\right.
\ee

\vskip 0.2 cm 
\noindent
The above expressions hold by analytic continuation for $\lambda\in{\cal P}_\eta$. 
\epr

\bcr
\label{dodici}

Let $\lp_k=N_k\in\mathbb{N}$. Let $c_{jk}$ denote $c_{jk}(\eta)= c_{jk}^{(\nu)}$.  The vector function $\hat{\Psi}_k(\lambda)$ in (\ref{hatphi}), $\lambda\in{\cal P}_\eta$, has the following behaviours at $\lambda_j\neq \lambda_k$. 

\vskip 0.3 cm 
\noindent
For  $\lp_j\not\in\mathbb{ Z}$:
\be 
  \label{acid1} 
\hat{\Psi}_k(\lambda)
=
\vPsi_j^{[-N_k-1]}(\lambda)c_{jk}+ \hbox{\rm reg}(\lambda-\lambda_j)
\ee

 \vskip 0.3 cm 
 \noindent
 For $\lp_j=N_j\in\mathbb{N}$, 
 \be 
  \label{acid2}
\hat{\Psi}_k(\lambda)
  =
\left\{
\matrix{  
\left(
  \hat{\Psi}_j^{[N_j-N_k]}(\lambda)\log(\lambda-\lambda_j)+{{P}_{N_j-N_k-1}(\lambda)\over 
  (\lambda-\lambda_j)^{N_j-N_k}}
  \right)c_{jk} +\hbox{\rm reg}(\lambda-\lambda_j)&  N_j\geq N_k+1
\cr
\cr
 \hat{\Psi}_j^{[-N_k+N_j]}(\lambda)\log(\lambda-\lambda_j)~c_{jk}~+\hbox{\rm reg}(\lambda-\lambda_j)& 
N_k\geq N_j
}
\right.
  \ee

\vskip 0.3 cm 
\noindent
 For $\lp_j=N_j\in\mathbb{Z}_{-}$:
  \be 
  \label{acid3}
  \hat{\Psi}_k(\lambda)= \vPsi_j^{[-N_k-1]}\ln(\lambda-\lambda_j)~c_{jk}+\hbox{\rm reg}(\lambda-\lambda_j),
\ee
\ecr
Note that (\ref{acid3}) always makes sense, because $c_{jk}=0$ for any $k=1,..,n$ when $\vPsi_j^{(sing)}\equiv 0$. 

\vskip 0.2 cm
\noindent
{\it Proof of Corollary \ref{dodici}:} We use the formulae of Proposition \ref{undici}. 
We observe that 
 $$
  \int_{\lambda_0}^\lambda d\xi_{N_k+1} \int_{\lambda_0}^{\xi_{N_k}}\cdots
  \int d\xi_1 \hat{\Psi}_k^{[N_k+1]}(\xi_1)=
$$
$$= 
  \hat{\Psi}_k(\lambda)-\sum_{l=0}^{N_k}{(-1)^{l+N_k}\over l!}~\vec{f}_l^{~(k)}~(\lambda-\lambda_k)^l-Q_{N_k}(\lambda-\lambda_0)
  $$
\be 
 \label{infar1}
= \hat{\Psi}_k(\lambda)+\hbox{reg}(\lambda)
\ee
  where $Q_{N_k}$ is a polynomial in $(\lambda-\lambda_0)$ of degree $N_k$ and $\hbox{reg}(\lambda)$ is analytic of $\lambda\in\mathbb{C}$ . 
Now recall that $\hat{\Psi}_k^{[N_k+1]}=\vPsi_k$. Using (\ref{psik}), we have
  $$
  \hat{\Psi}_k^{[N_k+1]}(\xi_1)=
  \left\{
  \matrix{
  \vPsi_j(\xi_1)c_{jk}+\hbox{reg}(\xi_1-\lambda_j),&~~~\lp_j\not\in \mathbb{Z}
  \cr
  \cr
  \left(
  \vPsi_j(\xi_1)\log(\xi_1-\lambda_j)+ {P^{(j)}(\xi_1)\over (\xi_1-\lambda_j)^{N_j+1}} 
  \right)c_{jk}+\hbox{reg}(\xi_1-\lambda_j),&~~~\lp_j=N_j\in\mathbb{N}
\cr
\cr
   \vPsi_j(\xi_1)\log(\xi_1-\lambda_j)c_{jk}+\hbox{reg}(\xi_1-\lambda_j),&~~~\lp_j\in\mathbb{Z}_-
}
  \right.
  $$
When $\lp_j\not\in\mathbb{ Z}$, from  (\ref{psik}) and (\ref{infar1}), we have
$$
\hat{\Psi}_k(\lambda)= \hbox{reg}(\lambda)+\int_{\lambda_0}^\lambda d\xi_{N_k+1} \int_{\lambda_0}^{\xi_{N_k}}\cdots
  \int d\xi_1 \left[
  \vPsi_j(\xi_1)c_{jk}+\hbox{reg}(\xi_1-\lambda_j)
  \right]
  =
  $$
  $$
  = \hbox{reg}(\lambda)+\Phi_j^{[-N_k-1]}(\lambda)c_{jk}+ \hbox{reg}(\lambda-\lambda_j)
  $$
 $$
  =\hbox{reg}(\lambda)+\Bigl(\vPsi_j^{[-N_k-1]}(\lambda)-Q_{N_k}(\lambda-\lambda_0)
  \Bigr)c_{jk}+ \hbox{reg}(\lambda-\lambda_j)
  $$
$$
=\vPsi_j^{[-N_k-1]}(\lambda)c_{jk}+ \hbox{reg}(\lambda-\lambda_j).
$$
 When $\lp_j=N_j\in\mathbb{N}$, 
 $$
\hat{\Psi}_k(\lambda)= \hbox{reg}(\lambda)+ 
$$
$$
+\int_{\lambda_0}^\lambda d\xi_{N_k+1} \int_{\lambda_0}^{\xi_{N_k}}\cdots
  \int d\xi_1 \left[\left({\vPsi}_j(\xi_1)\ln(\xi_1-\lambda_j)+{P_{N_j}^{(j)}(\xi_1)\over (\xi_1-\lambda_j)^{N_j+1}}\right)
  c_{jk}+\hbox{reg}(\xi_1-\lambda_j)\right]
  $$
  $$
  =
\left\{
\matrix{  
\left(
  \hat{\Psi}_j^{[N_j-N_k]}(\lambda)\log(\lambda-\lambda_j)+{P_{N_j-N_k-1}(\lambda)\over 
  (\lambda-\lambda_j)^{N_j-N_k}}
  \right)c_{jk} +\hbox{reg}(\lambda-\lambda_j)&  N_j\geq N_k+1
\cr
\cr
 \hat{\Psi}_j^{[-N_k+N_j]}(\lambda)\log(\lambda-\lambda_j)~c_{jk}~+\hbox{reg}(\lambda-\lambda_j)& 
N_k\geq N_j
}
\right.
 $$
  where the last step follows form Proposition \ref{undici}.   
  
\noindent
 When $\lp_j=N_j\in\mathbb{Z}_{-}$, again from Proposition \ref{undici} we have: 
$$
\hat{\Psi}_k(\lambda)= \hbox{reg}(\lambda)+ 
$$
$$
+\int_{\lambda_0}^\lambda d\xi_{N_k+1} \int_{\lambda_0}^{\xi_{N_k}}\cdots
  \int d\xi_1 \left(\vPsi_j(\xi_1)\ln(\xi_1-\lambda_j)
  c_{jk}+\hbox{reg}(\xi_1-\lambda_j)\right)
  $$
  $$
  =\vPsi_j^{[-N_k-1]}\ln(\lambda-\lambda_j)~c_{jk}+\hbox{reg}(\lambda-\lambda_j)
$$

$\Box$

\vskip 0.3 cm 
Next, we introduce the  "$\gamma$ deformed"  series corresponding to ${}_\gamma\vPsi_k$. For $\gamma_0>0$ sufficiently small and $0<\gamma<\gamma_0$, $A_1-\gamma$ has non integer diagonal entries and $A_1$ no integer eigenvalues, therefore:
$$
({}_\gamma\vPsi_k)^{[-q]}(\lambda)=(-1)^q~\sum_{l=q }^\infty \Gamma(\lp_k-\gamma+1-l)\vec{f}_{l-q}^{~(k)} ~(\lambda-\lambda_k)^{l-\lp_k+\gamma-1},~~~~~\forall q\geq 0
$$
and in particular  $({}_\gamma\vPsi_k)^{[0]}(\lambda)={}_\gamma\vPsi_k(\lambda)$. Recall that the coefficients $\vec{f}_{l-q}^{~(k)}$ are the same for any $\gamma\in \mathbb{C}$.


\ble
\label{lab9}
Let $0<\gamma<\gamma_0$ be such that $(A_1-\gamma)$ has no integer diagonal entries and  no integer eigenvalues. Let $q_1, q_2\in\mathbb{N}$. Then 
$$
\sqs{
\int_{\lambda_k}^\lambda ds~(\lambda-s)^{q_1+\gamma-1}(\vPsi_k)^{[-q_2]}(s)
={\Gamma(q_1+\gamma)~\sin \pi(\lp_k-\gamma) \over \sin \pi \lp_k }
~({}_\gamma\vPsi_k)^{[-q_1-q_2]}(\lambda),~~~~~\lp_k\not\in \mathbb{Z}
}
$$
$$
\sqs{
\int_{\lambda_k}^\lambda ds~(\lambda-s)^{q_1+\gamma-1}(\vPsi_k)^{[-q_2]}(s)
={\Gamma(q_1+\gamma)~\sin \pi\gamma \over  \pi  }
~({}_\gamma\vPsi_k)^{[-q_1-q_2]}(\lambda),~~~~~\lp_k\in\mathbb{Z}_{-}
}
$$
The branch of  $(\lambda-s)^\gamma$ in the integrals, for $\lambda\in{\cal P}_\eta$, is given by $\eta-2\pi<\arg(\lambda-s)|_{s=\lambda_k}<\eta$, and the continuous change along the path of integration. The integrals  are well defined for $0<\gamma<\gamma_0$, $q_1\geq 0$ and $q_2$ sufficiently big. 
\ele

\vskip 0.2 cm
\noindent
{\it Proof:} If $\lp_k\not\in \mathbb{Z}$ the statement is proved in \cite{BJL}, Lemma $2^\prime$. The same computations bring the result also for $\lp_k=N_k\in \mathbb{Z}_{-}$. It is enough to integrate expressions ({\ref{q1}) and ({\ref{q2}) term by term (where $|\lambda-\lambda_k|$ is small enough to make the series converge). In each term, the following  integral appears
$$
\int_{\lambda_k}^\lambda 
(\lambda-s)^{q_1+\gamma-1}(s-\lambda_k)^{l-\lp_k-1}~ds = (*)
$$
Since one can integrate along a line from $\lambda_k$ to $\lambda$, we parametrize the line with parameter $x\in[0,1]$ as follows: $s=\lambda_k+x(\lambda-\lambda_k)$. This yealds the integral representation of the Beta function $B(a,b)=\Gamma(a)\Gamma(b)/\Gamma(a+b)$.  Indeed
$$
(*)= (\lambda-\lambda_k)^{q_1+\gamma+l-\lp_k-1}~\int_0^1(1-x)^{q_1+\gamma-1}x^{l-\lp_k-1}dx
$$
$$= (\lambda-\lambda_k)^{q_1+\gamma+l-\lp_k-1}~{\Gamma(q_1+\gamma)\Gamma(l-\lp_k)\over \Gamma(q_1+\gamma+l-\lp_k)}
$$ 
The formula holds for any value of $\lp_k$. Note that $l\geq q_2$, thus if $q_1$ and $q_2$ are big enough, the integrals converge. Note also that for $\lp_k=N_k\leq -1$, the integrals converge for $q_2\geq 0$. Moreover, since we have assumed $\gamma>0$, the integrals converge for any $q_1\geq 0$. 
Finally, some manipulations using $\Gamma(x)\Gamma(1-x)= \pi/\sin(\pi x)$, yield the result. For example, in case $\lp_k=N_k\leq -1$, we have
$$
\int_{\lambda_k}^\lambda ds~(\lambda-s)^{q_1+\gamma-1}(\vPsi_k)^{[-q_2]}(s)
=
$$
$$
=
(-1)^{q_2}\sum_{l\geq q_2}{(-1)^{l-N_k}\over (l-N_k-1)!}~ \vec{f}_{l-q_2}^{~(k)}~(\lambda-\lambda_k)^{q_1+\gamma+l-N_k-1}~{\Gamma(q_1+\gamma)\Gamma(l-N_k)\over \Gamma(q_1+\gamma+l-N_k)}=(**)
$$ 
We use 
{\small
$$\Gamma(l-N_k)=(l-N_k-1)!~,~~~~~{1\over \Gamma(q_1+\gamma+l-N_k)}={\Gamma(N_k+1-\gamma-l-q_1)  \sin(q_1+l-N_k+\gamma)\over \pi}
$$  
}and change $l \mapsto l-q_1$ . We get
$$
(**)= (-1)^{q_1+q_2}{\Gamma(q_1+\gamma)\sin\pi\gamma\over \pi}~\sum_{l\geq q_1+q_2}\Gamma(N_k-\gamma_1-l)\vec{f}_{l-q_1-q_2}^{~(k)}~(\lambda-\lambda_k)^{l-(N_k-\gamma)-1}
$$
$$={\Gamma(q_1+\gamma)~\sin \pi\gamma \over  \pi  }
~({}_\gamma\vPsi_k)^{[-q_1-q_2]}(\lambda).
$$
$\Box$

 \ble
\label{lab10}
Let $0<\gamma<\gamma_0$ be such that $(A_1-\gamma)$ has no integer diagonal entries and no integer eigenvalues. Then
$$
\sqs{
\int_{\lambda_k}^\lambda ds~(\lambda-s)^{\gamma-1}(\hat{\Psi}_k)^{[-q]}(s)
={\Gamma(\gamma)~\sin \pi\gamma \over  \pi  }
~({}_\gamma\vPsi_k)^{[-N_k-1-q]}(\lambda),~~~~~\lp_k=N_k\in\mathbb{N}
}
$$
The integral is well  defined for $0<\gamma<\gamma_0$ and $q\geq 0$ integer. The branch of $(\lambda-s)^\gamma$ is defined in the same way as in Lemma \ref{lab9}. 
\ele

\vskip 0.2 cm
\noindent
{\it Proof:} Integration term by term yields 
$$
(\hat{\Psi}_k)^{[-q]}(s)=(-1)^{N_k+1+q}\sum_{l\geq q} {(-1)^{l+1}\over l!} \vec{f}_{l-q}^{~(k)} (\lambda-\lambda_k)^l,
$$
$$
({}_\gamma\vPsi_k)^{[-N_k-1-q]}(\lambda)=(-1)^{N_k+1+q}\sum_{l\geq N_k+1+q}  \Gamma(N_k-\gamma+1-l)\vec{f}_{l-N_k-1-q}^{~(k)} ~(\lambda-\lambda_k)^{l-N_k+\gamma-1}.
$$
As in the previous lemma, we compute
$$
 \int_{\lambda_k}^\lambda ds~(\lambda-s)^{\gamma-1}(s-\lambda_k)^l=(\lambda-\lambda_k)^{\gamma+l}~ {\Gamma(\gamma)\Gamma(l+1)\over \Gamma(\gamma+l+1)}
 $$
 and use 
 $1/\Gamma(\gamma+l+1)=(-1)^{l+1}\sin(\pi\gamma)\Gamma(-\gamma-l)/\pi$.
 This implies that 
 $$
 \int_{\lambda_k}^\lambda ds~(\lambda-s)^{\gamma-1}(\hat{\Psi}_k)^{[-q]}(s)
=(-1)^{N_k+1+q}{\Gamma(\gamma) \sin\pi\gamma\over\pi} \sum_{l\geq q} \Gamma(-l-\gamma)\vec{f}_{l-q}^{~(k)}(\lambda-\lambda_k)^{l+\gamma}.
$$
After redefining $l^\prime=l+N_k+1$ we obtain the final result. $\Box$

\vskip 0.3 cm 
 We establish the monodromy of $\vPsi_k^{[-q]}$  and $\hat{\Psi}_k^{[-q]}$ in the following

 \bpr
 \label{lab11}
 Let $\lambda\in {\cal P}_\eta$. Let $q\geq 0$ be an integer. 
 Let   $\alpha_j=2\pi i$ when $\lp_j\in\mathbb{Z}$, and $\alpha_j=e^{-2\pi i \lp_j}-1$ when $\lp_j\not\in \mathbb{Z}$. The following transformations hold for a loop $\gamma_j$ around  a pole $\lambda_j$. 
 \vskip 0.2 cm 
 
 a) If $\lp_k\not \in \mathbb{Z}$ or $\lp_k\in \mathbb{Z}_{-}$:
 $$
  \vPsi_k^{[-q]}(\lambda)\longmapsto
   \vPsi_k^{[-q]}(\lambda)+
   \left\{
   \matrix{
   \alpha_j c_{jk} ~\vPsi_j^{[-q]}(\lambda),&~~~\lp_j\not\in \mathbb{Z} \hbox{ or } \lp_j\in\mathbb{Z}_{-}
   \cr
   \cr
   \alpha_j c_{jk}~\hat{\Psi}_j^{[-q+N_j+1]}(\lambda),& ~~~\lp_j\in\mathbb{N}
}
\right.
  $$

 \vskip 0.2 cm 
 
 b) If $\lp_k\in\mathbb{N}$:
 $$
  \hat{\Psi}_k^{[-q]}(\lambda)\longmapsto
   \hat{\Psi}_k^{[-q]}(\lambda)+
   \left\{
   \matrix{
   \alpha_j c_{jk} ~\vPsi_j^{[-q-1-N_k]}(\lambda),&~~~ \lp_j\not\in\mathbb{Z} \hbox{ or } \lp_j\in\mathbb{Z}_-
   \cr
   \cr
   \alpha_j c_{jk}~\hat{\Psi}_j^{[-q+N_j-N_k]}(\lambda),& ~~~ \lp_j\in \mathbb{N},
}
\right.
  $$

 \epr

 \vskip 0.2 cm 
 \noindent
 {\it Proof:} We consider the function $ \Phi_k^{[-q]}(\lambda,\lambda_0)$ defined in (\ref{cistit3})  for $\lambda,~\lambda_0\in{\cal P}_\eta$. For simplicity of notation, we omit $\lambda_0$, namely we write 
$$
 \Phi_k^{[-q]}(\lambda)=\int_{\lambda_0}^\lambda d\xi_q\int_{\lambda_0}^{\xi_q}d\xi_{q-1} \cdots \int_{\lambda_0}^{\xi_2}d\xi_1 
 ~\phi_k(\xi_1).
 $$
$$
\Phi_k^{[-q]}(\lambda)=\phi_k(\lambda),~~~~\hbox{ if } q=0
$$
 In the cut-plane ${\cal P}_\eta$, we  consider $\lambda$ close to $\lambda_j\neq \lambda_k$ in such a way that the series representations of $(\phi_j)^{[-q]}(\lambda)$ converge. We also consider a loop $\gamma_j$ around $\lambda_j$ in counter-clockwise direction, represented by $(\lambda-\lambda_j)\mapsto (\lambda-\lambda_j)e^{2\pi i}$. See figure \ref{fig7}.  We have the following transformation after the loop
\begin{figure}
\centerline{\includegraphics[width=0.3\textwidth]{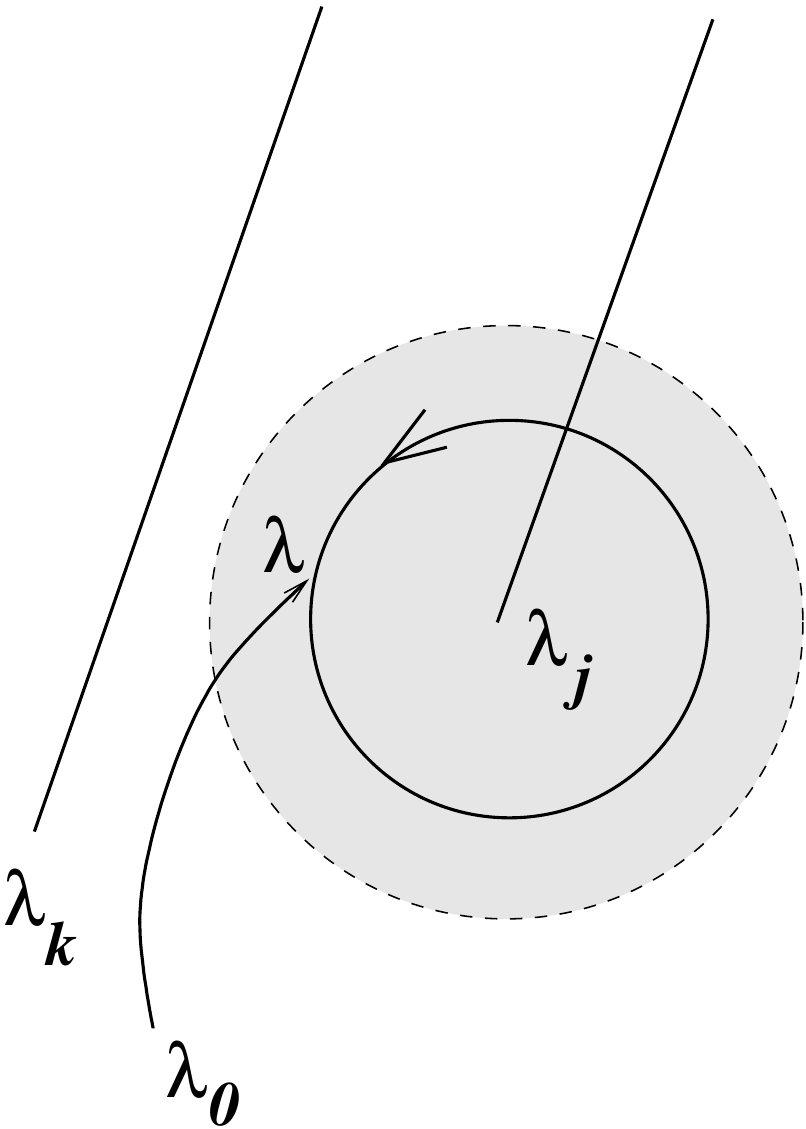}}
\caption{The paths of analytic continuation.}
\label{fig7}
\end{figure}
 $$
 \Phi_k^{[-q]}(\lambda)\longmapsto \Phi_k^{[-q]}(\lambda)+\oint_{\gamma_j}d\xi_q\int_{\lambda_0}^{\xi_q}d\xi_{q-1} \cdots \int_{\lambda_0}^{\xi_2}d\xi_1 
 ~\phi_k(\xi_1).
 $$
By formula (\ref{cistit4}) it follows that the analytic continuation of $(\phi_k)^{[-q]}(\lambda)$ along the loop $\gamma_j$ is 
 $$
 (\phi_k)^{[-q]}(\lambda)\longmapsto (\phi_k)^{[-q]}(\lambda)+\oint_{\gamma_j}d\xi_q\int_{\lambda_0}^{\xi_q}d\xi_{q-1} \cdots \int_{\lambda_0}^{\xi_2}d\xi_1 
 ~\phi_k(\xi_1),~~~~~~\forall q\geq 0.
 $$
Next, we express $\phi_k(\xi_1)$ in terms of the solutions $\phi_j$ at $\lambda_j$. We distinguish the two cases in the lemma.

\vskip 0.2 cm 
\noindent
{\bf Case a):}  $\lp_k\not \in \mathbb{Z}$ or $\lp_k \in\mathbb{Z}_{-}$. We have $
\phi_k(\xi_1)=\vPsi_k(\xi_1)$, therefore we use  (\ref{psik}), namely 
$$
\phi_k(\xi_1)=
\left\{
\matrix{
\vPsi_j(\xi_1)c_{jk}+\hbox{reg}(\xi_1-\lambda_j), &~~~~~\lp_j\not\in \mathbb{ Z}
\cr
\cr
\vPsi_j(\xi_1)\ln(\xi_1-\lambda_j)~c_{jk} + \hbox{reg}(\xi_1-\lambda_j), &~~~~~\lp_j\in\mathbb{Z}_{-}
\cr
\cr
\left(\vPsi_j(\xi_1)\ln(\xi_1-\lambda_j)+{P^{(j)}(\xi_1)\over (\xi_1-\lambda_j)^{N_j+1}}\right)c_{jk} + \hbox{reg}(\xi_1-\lambda_j), &~~~~~\lp_j=N_j\in\mathbb{ N}
}
\right.
$$

{\bf a.1)} When $\lp_j\not\in \mathbb{Z}$, we have 
 $$
 \oint_{\gamma_j}d\xi_q\int_{\lambda_0}^{\xi_q}d\xi_{q-1} \cdots \int_{\lambda_0}^{\xi_2}d\xi_1 
 ~\left(\vPsi_j(\xi_1)c_{jk}+\hbox{reg}(\xi_1-\lambda_j)\right)
 =
 $$
 $$
 =\oint_{\gamma_j}d\xi_q~ \Phi_j^{[-(q-1)]}(\xi_q)~c_{jk}+\oint_{\gamma_j}d\xi_q~\hbox{reg}(\xi_q-\lambda_j)
\equiv 
\oint_{\gamma_j}d\xi_q ~\Phi_j^{[-(q-1)]}(\xi_q)~c_{jk}~+0,
$$
because the loop integral of regular terms at $\lambda_j$ vanishes. Now, by (\ref{recu}) and (\ref{cistit4}), we have 
$$
\oint_{\gamma_j}d\xi_q~ \Phi_j^{[-(q-1)]}(\xi_q)~c_{jk}=
c_{jk}~\Phi_j^{[-q]}(\xi_q)\Big|_{\lambda-\lambda_j}^{(\lambda-\lambda_j)e^{2\pi i}}
$$
$$
\equiv
c_{jk}~\vPsi_j^{[-q]}(\xi_q)\Big|_{\lambda-\lambda_j}^{(\lambda-\lambda_j)e^{2\pi i}}~= \vPsi_j^{[-q]}(\lambda)(e^{-2\pi i\lp_j}-1)c_{jk},~~~~~q\geq 0
 $$ 
 The last step follows from the series representation (\ref{q1}).  This proves the Lemma in case a.1).

{\bf a.2)}  When $\lp_j\in\mathbb{Z}_{-}$, we use (\ref{REE1}) and compute
$$
\oint_{\gamma_j}d\xi_q\int_{\lambda_0}^{\xi_q}d\xi_{q-1} \cdots \int_{\lambda_0}^{\xi_2}d\xi_1 
\Bigl(\vPsi_j(\xi_1)\ln(\xi_1-\lambda_j)~c_{jk} + \hbox{reg}(\xi_1-\lambda_j)\Bigr)=
$$
$$
= c_{jk}~(\vPsi)_j^{[-q]}(\xi_q)\ln(\xi_q-\lambda_j)\Big|_{\lambda-\lambda_j}^{(\lambda-\lambda_j)e^{2\pi i}}=2\pi i  c_{jk} ~
(\vPsi)_j^{[-q]}(\lambda),~~~~~~q\geq 0
$$
This implies the Lemma in case a.2).

 {\bf a.3)}  When $\lp_j=N_j\in\mathbb{N}$, we use (\ref{REE2}) and compute
 $$
 \oint_{\gamma_j}d\xi_q \int_{\lambda_0}^{\xi_q}d\xi_{q-1} \cdots \int_{\lambda_0}^{\xi_2}d\xi_1  \left[
 \left(\vPsi_j(\xi_1)\ln(\xi_1-\lambda_j)+{P^{(j)}(\xi_1)\over (\xi_1-\lambda_j)^{N_j+1}}\right)c_{jk}  +\hbox{reg}(\xi_1-\lambda_j)\right]
 =
$$
$$=
c_{jk}~ \hat{\Psi}_j^{[-q+N_j+1]}(\xi_q)\ln(\xi_q-\lambda_j)\Big|_{\lambda-\lambda_j}^{(\lambda-\lambda_j)e^{2\pi i}}
=
 2\pi i c_{jk} \hat{\Psi}_j^{[-q+N_j+1]}(\lambda),~~~~~q\geq N_j+1
 $$
This proves the Lemma in case a.3).  
 
 \vskip 0.3 cm 
 \noindent
 {\bf Case b):} when $\lp_k=N_k\in \mathbb{N}$, we need to compute 
 \be 
 \label{acid}
 \oint_{\gamma_j}d\xi_q \int_{\lambda_0}^{\xi_q}d\xi_{q-1} \cdots \int_{\lambda_0}^{\xi_2}d\xi_1 \hat{\Psi}_k(\xi_1).
 \ee 
 
{\bf b.1)} In the case of $\lp_j\not\in\mathbb{Z}$, we use (\ref{acid1}), and find
$$
(\ref{acid})=\oint_{\gamma_j}d\xi_q \int_{\lambda_0}^{\xi_q}d\xi_{q-1} \cdots \int_{\lambda_0}^{\xi_2}d\xi_1 \Bigl(\vPsi_j^{[-N_k-1]}(\xi_1)c_{jk}+\hbox{reg}(\xi_1-\lambda_j)\Bigr)
$$
$$
=\Bigl(c_{jk}\vPsi_j^{[-N_k-1-q]}(\xi_q)+\hbox{reg}(\xi_q-\lambda_j)\Bigr)\Big|_{(\lambda-\lambda_j)}^{e^{2\pi i }(\lambda-\lambda_j)}
=c_{jk}(e^{-2\pi i \lp_j}-1)\vPsi_j^{[-N_k-1-q]}(\lambda).
$$

{\bf b.2)} In the case of $\lp_j\in\mathbb{N}$, we use (\ref{acid2}) and Proposition \ref{undici}, and find 
$$
(\ref{acid}) =
$$
$$
\oint_{\gamma_j}d\xi_q \int_{\lambda_0}^{\xi_q}d\xi_{q-1} \cdots \int_{\lambda_0}^{\xi_2}d\xi_1 
\left[\left(
\hat{\Psi}_j^{[N_j-N_k]}(\xi_1)\log(\xi_1-\lambda_j)+{P_{N_j-N_k-1}(\xi_1)\over 
  (\xi_1-\lambda_j)^{N_j-N_k}}
  \right)c_{jk} +\hbox{reg}(\xi_1-\lambda_j)
\right]
$$
where ${P}_{N_j-N_k-1}=0$ for $N_k\geq N_j$.  

\noindent
For $
0\leq q \leq N_j-N_k-1$, the integral is 
$$
=\left[
c_{jk} \left(
\vPsi_j^{[N_j-N_k-q]}(\xi_q)\ln(\lambda-\lambda_j)+{P_{N_j-N_k-1-q}^{(j)}(\xi_q)\over
(\xi_q-\lambda_j)^{N_j-N_k-q} }
\right)
+\hbox{reg}(\xi_q-\lambda_j)
\right]
\Big|_{(\lambda-\lambda_j)}^{e^{2\pi i }(\lambda-\lambda_j)}.
$$ 
For
 $q\geq  N_j-N_k\geq 0
$, the integral is 
$$
=\left[
c_{jk}\vPsi_j^{[N_j-N_k-q]}(\xi_q)\ln(\xi_q-\lambda_j)
+\hbox{reg}(\xi_q-\lambda_j)
\right]
\Big|_{(\lambda-\lambda_j)}^{e^{2\pi i }(\lambda-\lambda_j)}.
$$
In both cases, the above expressions yield
$$
(\ref{acid}) =2\pi i c_{jk} \vPsi_j^{[N_j-N_k-q]}(\lambda)
$$

{\bf b.3)} In case $\lp_j\in\mathbb{Z}_-$, we use (\ref{acid3}) and Proposition \ref{undici}, and find
$$
(\ref{acid}) =\oint_{\gamma_j}d\xi_q \int_{\lambda_0}^{\xi_q}d\xi_{q-1} \cdots \int_{\lambda_0}^{\xi_2}d\xi_1 
\Bigl(
\vPsi_j^{[-N_k-1]}(\xi_1)\ln(\xi_1-\lambda_j)c_{jk}+\hbox{reg}(\xi_1-\lambda_j)
\Bigr)
$$

$$
=\left[
\vPsi_j^{[-N_k-1-q]}(\xi_q)\ln(\xi_q-\lambda_j)~c_{jk}~
+\hbox{reg}(\xi_q-\lambda_j)
\right]
\Big|_{(\lambda-\lambda_j)}^{e^{2\pi i }(\lambda-\lambda_j)}=
2\pi i  c_{jk}~ \vPsi_j^{[-N_k-1-q]}(\lambda).
$$

\vskip 0.2 cm 
\noindent
   The above computations  imply  the Lemma in case b). $\Box$

\vskip 0.3 cm 
\noindent  
  {\it Proof of Proposition \ref{fundlpro}:}  Given a function $f(\lambda)$, $\lambda\in{\cal P}_\eta$, we denote with   $f_{+}(\lambda)$ the value on the left side of  $L_j$, where  $\arg(\lambda-\lambda_j)=\eta-2\pi$.  We denote with  $f_{-}(\lambda)$ the value on the right side, where $\arg(\lambda-\lambda_j)=\eta$.  
  By Lemmas \ref{lab9} and \ref{lab10} we have:
  $$
  \left({}_\gamma \vPsi_k \right)^{[-q(\lp_k)]}_{\pm}(\lambda)
  =
  F(\lp_k)\left\{
  \int_{\lambda_k}^{\lambda_j} (\lambda-s)^{\gamma-1}\phi_k^{[-q]}(s)ds~+
   \int_{\lambda_j}^{\lambda} (\lambda-s)^{\gamma-1}\left(\phi_k^{[-q]}\right)_{\pm}(s)ds
  \right\}
$$
where 
$$
q(\lp_k)=\left\{
\matrix{
q & \hbox{ if } \lp_k\not\in \mathbb{Z}  \hbox{ or } \lp_k\in\mathbb{Z}_{-},
\cr
\cr
q+N_k+1 & \hbox{ if } \lp_k\in \mathbb{ N} ,
}
\right.
$$

$$
\phi_k=\left\{
\matrix{
\vPsi_k & \hbox{ if } \lp_k\not\in \mathbb{Z}  \hbox{ or } \lp_k\in\mathbb{Z}_{-},
\cr
\cr
\hat{\Psi}_k & \hbox{ if } \lp_k\in \mathbb{ N} ,
}
\right.
$$

$$
 F(\lp_k)=\left\{
\matrix{
{\sin \pi \lp_k\over \Gamma(\gamma)\sin \pi(\lp_k-\gamma)} & \hbox{ if } \lp_k\not\in \mathbb{Z},
\cr\cr
{\pi \over \Gamma(\gamma) \sin \pi \gamma}  & \hbox{ if } \lp_k\in \mathbb{Z}.
}
\right.
$$
In the integral, $\arg(\lambda-s)$, $s\in {\cal P}_\eta$,   has the value obtained by the continuous change along the path of integration from $\lambda_k$   up  to $s$ belonging to $L_j$. 
 Change from $f_{-}$ to $f_{+}$ is obtained  along a small loop encircling only $\lambda_j$.  Therefore,  Lemma  \ref{lab11} yields
  \be 
  \label{tm1}
  \left({}_\gamma \vPsi_k \right)^{[-q(\lp_k)]}_{-}(\lambda)-\left({}_\gamma \vPsi_k \right)^{[-q(\lp_k)]}_{+}(\lambda)
  = \alpha_j[\gamma]c_{jk}[\gamma]\left({}_\gamma \vPsi_j \right)^{[-q(\lp_k)]}_{+}(\lambda)
  .
  \ee
  By Lemmas \ref{lab9} and \ref{lab10} we write
  \be 
  \label{tm2}
  \left({}_\gamma \vPsi_k \right)^{[-q(\lp_k)]}_{-}(\lambda)-\left({}_\gamma \vPsi_k \right)^{[-q(\lp_k)]}_{+}(\lambda)=
  F(\lp_k)\int_{\lambda_j}^\lambda (\lambda-s)^{\gamma-1}\left[  \left(\phi_k^{[-q]}\right)_{-}(s)-\left(\phi_k^{[-q]}\right)_+(s)\right]
  \ee

  We need to distinguish  two cases. 
  
  \vskip 0.2 cm 
 1) $\lp_k\not\in \mathbb{Z}$, or $\lp_k\in \mathbb{Z}_{-}$. In this case (\ref{tm1}), (\ref{tm2}) and Lemma \ref{lab11} applied to the integrand yield the following equalities. 
 
 1.a) for $\lp_j\not \in\mathbb{Z}$ or $\lp_j\in \mathbb{Z}_{-}$: 
 \be 
 \label{tm3}
  \alpha_j[\gamma]c_{jk}[\gamma]\left({}_\gamma \vPsi_j \right)^{[-q]}_{+}(\lambda)
  =
  F(\lp_k)\int_{\lambda_j}^\lambda ds~(\lambda-s)^{\gamma-1} \alpha_j c_{jk} \left(\vPsi_j\right)^{[-q]}_+(s),
  \ee

  1.b) for $\lp_j\in\mathbb{ N}$: 
 \be 
 \label{tm4}
    \alpha_j[\gamma]c_{jk}[\gamma]\left({}_\gamma \vPsi_j \right)^{[-q]}_{+}(\lambda)
  = F(\lp_k)\int_{\lambda_j}^\lambda ds~(\lambda-s)^{\gamma-1} \alpha_j c_{jk} \left(\hat{\Psi}_j\right)^{[-q+N_j+1]}_+(s),
  \ee
  We apply again Lemmas (\ref{lab9}) and (\ref{lab10}) to express the r.h.s. of the above equalities. To this end, we need $(\lambda-s)_+$ in the integrand. Observe that  
  $$
  (\lambda-s)^{\gamma-1} \hbox{ in the integrand }
  =\left\{
  \matrix{  
  ~~(\lambda-s)_-^{\gamma-1}=\left[e^{2\pi i}(\lambda-s)_+\right]^{\gamma -1} & \hbox{  when } k \succ j
  \cr 
 \cr 
  (\lambda-s)_{+}^{\gamma-1} & \hbox{  when } k \prec j.
}
\right.   
   $$
 Indeed,  when  $\lambda$ belongs to the left side of  $L_j$  
and  $s \in {\cal P}_\eta$, then   $\eta-2\pi<\arg(\lambda-s)<\eta$. When $s$ reaches $ L_j$ from the left, then $\arg(\lambda- s)\to \eta$ if  $L_k$  is to the left of $L_j$, namely $k\succ j$: in this case we obtain $(\lambda-s)_-$. On the other hand, $\arg(\lambda- s)\to \eta-2\pi $ if $L_k$ is to the right of $L_j$, namely $k \prec j$: in this case we obtain $(\lambda-s)_+$. See figure \ref{fig8}. 
\begin{figure}
\centerline{\includegraphics[width=0.6\textwidth]{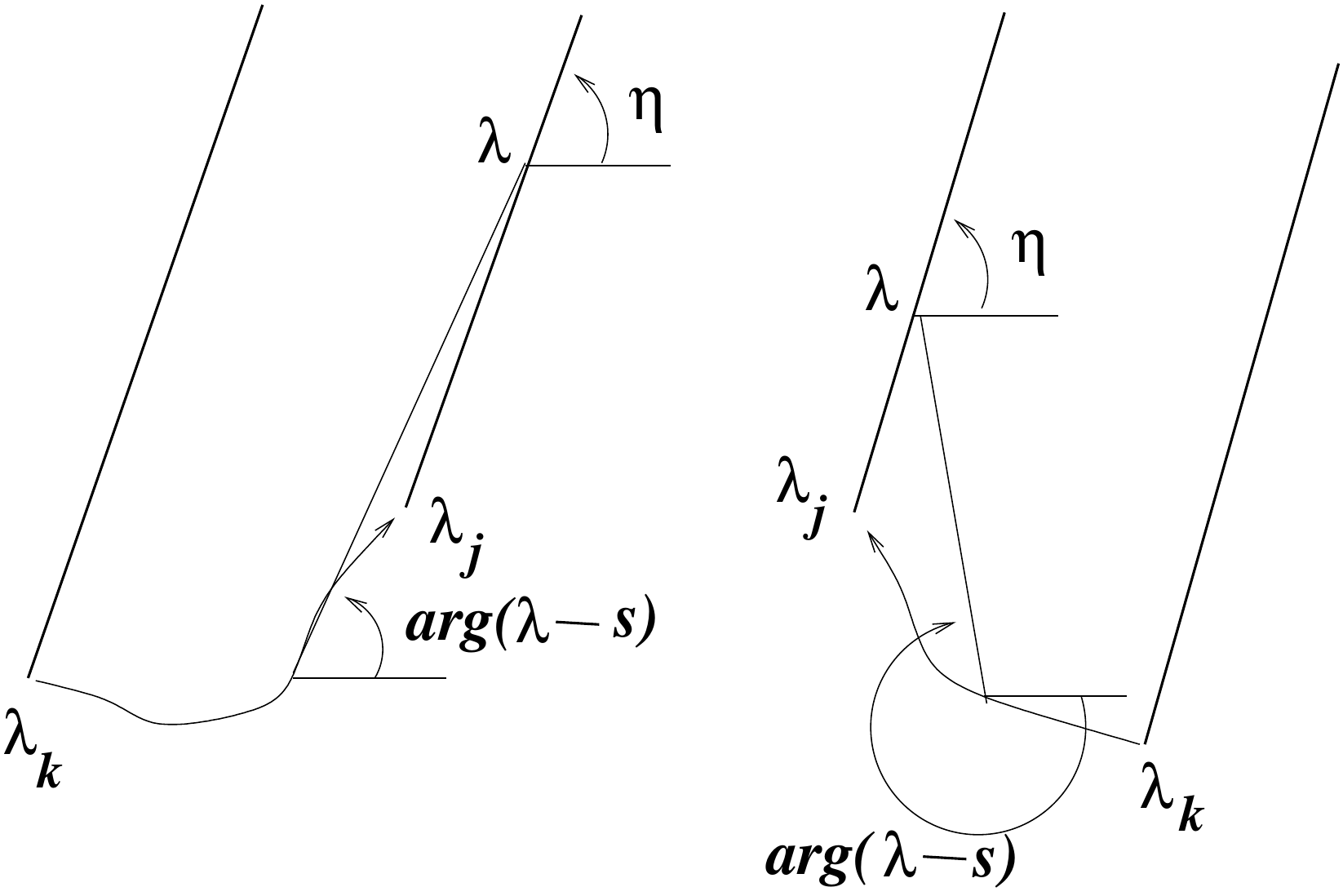}}
\caption{The figure shows $\arg(\lambda-s)$ as $s\to \lambda_j$.}
\label{fig8}
\end{figure}
Applying  Lemmas (\ref{lab9}) and (\ref{lab10}) we find
$$
 \hbox{r.h.s. of (\ref{tm3}) and (\ref{tm4}) }=  
 \left\{
  \matrix{  {F(\lp_k)\over F(\lp_j)}  e^{2\pi i \gamma} \alpha_j c_{jk} \left({}_\gamma\vPsi_j \right)_+^{[-q]}(\lambda)  & k \succ j,
  \cr
  \cr
     {F(\lp_k)\over F(\lp_j)} \alpha_j c_{jk} \left({}_\gamma\vPsi_j \right)_+^{[-q]}(\lambda)  & k \prec j.
   }
   \right.
$$
Namely:
\be 
\label{tm8}
 \alpha_j[\gamma]c_{jk}[\gamma]=  
 \left\{
  \matrix{{F(\lp_k)\over F(\lp_j)}  e^{2\pi i \gamma}  \alpha_j c_{jk}  & k \succ j,
  \cr
  \cr
     {F(\lp_k)\over F(\lp_j)} \alpha_j c_{jk}  & k \prec j.
   }
   \right.
\ee
Finally, we compute the ratio $F(\lp_k)/ F(\lp_j)$. For $\lp_j\not\in \mathbb{Z}$:
$$
 {F(\lp_k)\over F(\lp_j)}= {\sin\pi\lp_k~\sin\pi(\lp_j-\gamma)\over \sin\pi\lp_j~\sin\pi(\lp_k-\gamma)}={(1-e^{-2\pi i \lp_k})(1-e^{-2\pi i(\lp_j-\gamma)})
 \over (1-e^{-2\pi i \lp_j})(1-e^{-2\pi i(\lp_k-\gamma)})}={\alpha_k\alpha_j[\gamma]\over 
 \alpha_j\alpha_k[\gamma]}.
 $$
 For $\lp_j\in\mathbb{Z}$:
$$
 {F(\lp_k)\over F(\lp_j)}= {\sin\pi\gamma~\sin\pi\lp_k\over \pi ~\sin\pi(\lp_k-\gamma)}
 ={e^{2\pi i \gamma}-1\over 2\pi i }{e^{-2\pi i \lp_k}-1\over e^{-2\pi i (\lp_k-\gamma)}-1}={\alpha_j[\gamma]\over \alpha_j}{\alpha_k\over \alpha_k[\gamma]}
 $$
 where we have used the fact that $\alpha_j=2\pi i$ and $\alpha_j[\gamma]=e^{2\pi i\gamma}-1$.
 
 \noindent
 The above computations imply the statement of Proposition \ref{fundlpro} when $\lp_k\not\in\mathbb{Z}$ and $\lp_k\in \mathbb{Z}_{-}$.

 \vskip 0.2 cm 
 2) $\lp_k\in \mathbb{N}$. In this case (\ref{tm1}), (\ref{tm2}) and Lemma \ref{lab11} applied to the integrand yield the following equalities. 
 
 2.a) for $\lp_j\not\in \mathbb{Z}$ or $\lp_j\in\mathbb{Z}_{-}$: 
 \be 
 \label{tm5}
  \alpha_j[\gamma]c_{jk}[\gamma]\left({}_\gamma \vPsi_j \right)^{[-q-N_k-1]}_{+}(\lambda)
  =
  F(\lp_k)\int_{\lambda_j}^\lambda ds~(\lambda-s)^{\gamma-1} \alpha_j c_{jk} \left(\vPsi_j\right)^{[-q-N_k-1]}_+(s)
  \ee

   2.b) for $\lp_j\in \mathbb{N}$: 
 \be 
 \label{tm6}
    \alpha_j[\gamma]c_{jk}[\gamma]\left({}_\gamma \vPsi_j \right)^{[-q-N_k-1]}_{+}(\lambda)
  = F(\lp_k)\int_{\lambda_j}^\lambda ds~(\lambda-s)^{\gamma-1} \alpha_j c_{jk} \left(\hat{\Psi}_j\right)^{[-q-N_k-1+N_j+1]}_+(s) .
  \ee
  We apply again Lemmas (\ref{lab9}) and (\ref{lab10}) to express the r.h.s. of the above equalities, keeping into account the branch of  $(\lambda-s)^{\gamma-1}$  as before. We find
$$
 \hbox{r.h.s. of (\ref{tm5}) and (\ref{tm6}) }=  
 \left\{
  \matrix{  {F(\lp_k)\over F(\lp_j)}  e^{2\pi i \gamma} \alpha_j c_{jk} \left({}_\gamma\vPsi_j \right)_+^{[-q-N_k-1]}(\lambda)  & k \succ j,
  \cr
  \cr
     {F(\lp_k)\over F(\lp_j)} \alpha_j c_{jk} \left({}_\gamma\vPsi_j \right)_+^{[-q-N_k-1]}(\lambda)  & k \prec j.
   }
   \right.
$$
Namely, we obtain again 
\be 
\label{tm10}
 \alpha_j[\gamma]c_{jk}[\gamma]=  
 \left\{
  \matrix{{F(\lp_k)\over F(\lp_j)}  e^{2\pi i \gamma}  \alpha_j c_{jk}  & k \succ j,
  \cr
  \cr
     {F(\lp_k)\over F(\lp_j)} \alpha_j c_{jk}  & k \prec j.
   }
   \right.
\ee
  Finally, we compute the ratio $F(\lp_k)/ F(\lp_j)$. For $\lp_j\not\in \mathbb{Z}$:
$$
 {F(\lp_k)\over F(\lp_j)}= {\pi\over \sin\pi\gamma}{\sin\pi(\lp_j-1)\over \sin\pi\lp_j}={2\pi i \over e^{2\pi i \gamma}-1}{e^{-2\pi i (\lp_j-\gamma)}-1\over e^{-2\pi i \lp_j}-1}={\alpha_k\over \alpha_k[\gamma]}{\alpha_j[\gamma]\over \alpha_j}
 $$
 where we have used the fact that $\alpha_k=2\pi i$ and $\alpha_k[\gamma]=e^{2\pi i\gamma}-1$.

\noindent 
 For $\lp_j\in\mathbb{Z}$:
$$
 {F(\lp_k)\over F(\lp_j)}=1
 $$
 In this last case, observe that $\alpha_k=\alpha_j=2\pi i $ and $\alpha_k[\gamma]=\alpha_j[\gamma]=e^{2\pi i \gamma}-1$.

 \noindent
 The above computations imply the statement of Proposition \ref{fundlpro} when $\lp_k\in\mathbb{N}$. $\Box$

 
 \section{Appendix}
 
 Some statements of the main body of the paper are proved in \cite{BJL} with  the {\it assumption (i)} of \cite{BJL}, namely $\lp_1,...,\lp_n$ not integers.  Here we  prove the statements    without assumptions on $A_1$.

\vskip 0.3 cm 
\noindent
{\bf Proof of Proposition \ref{prop1}:} It goes as the proof of Proposition 1 in  \cite{BJL}, which  was done assuming  no integer $\lp_k$'s, $k=1,...,n$.  We repeat the proof with no assumptions on the $\lp_k$'s. For the proof, recall Remark \ref{remark1} and Lemma \ref{lemma1}. 

Vector solutions of the equation (\ref{02}) form a $n$ dimensional linear space. In Section \ref{locall}  we have proved that for any $j=1,...,n$, the vector solutions regular at one $\lambda_j$ form a linear space $V_j$ of dimension dim$V_j\geq n-1$. Now, choose a $\lambda_k$.  The number   of independent solutions which are regular at each of the $n-1$ poles   $\lambda_j\neq \lambda_k$ is equal to  dim$(\cap_{j\neq k} V_j)\geq 1$.     Therefore, there exist at least one vector solutions $\vPsi^*_k(\lambda)$, analytic at all $\lambda_j\neq \lambda_k$, $1\leq j\leq n$, $j\neq k$.  

If dim$(\cap_{j\neq k} V_j)>1$, then  there must exist a polynomial solution of (\ref{02}). Indeed, 
suppose that dim$(\cap_{j\neq k} V_j)>1$. So there are at least two independent solutions, regular at each
 $\lambda_j\neq \lambda_k$. Let them be $\vpsi_1(\lambda)$ and $\vpsi_2(\lambda)$. At $\lambda_k$ they
 must have representation $\vpsi_q=c_q~\vPsi_k^{(sing)}+\hbox{reg}_q(\lambda-\lambda_k)$, $q=1,2$, $c_q\in \mathbb{C}\backslash\{0\}$.  But now, $c_1^{-1}\vpsi_1-c_2^{-1}\vpsi_2$ is 
not vanishing (by independence of the two solutions) and  is analytic also at $\lambda=\lambda_k$. It follows that there
 is an analytic solution at each $\lambda_i$, $1\leq i \leq n$, which therefore must be polynomial (Remark \ref{remark1}). This can occur if and only if $A_1$ has at least one negative integer eigenvalue (Lemma \ref{lemma1}). 
Thus, if we assume that $A_1$ has no negative integer eigenvalues, so that there are no polynomial solutions, then dim$(\cap_{j\neq k}V_j)=1$. Therefore, for any $k$ there exists a unique (up to normalization) solution which is  regular at all  $\lambda_j\neq \lambda_k$ and {\it must be singular} at $\lambda_k$. It is the unique $\vPsi_k^*$ with normalization (\ref{star2}).  Obviously, $\vPsi_1^*,...,\vPsi_n^*$ are independent. 
In particular, observe that if $\lp_k\in-\mathbb{N}-2$, necessarily  the log-singular solution must exist (i.e. $\vPsi_k^{(sing)}\neq 0$).

Conversely, if $\vPsi_1^*(\lambda),...,\vPsi_n^*(\lambda)$ exist satisfying (\ref{star1}) and (\ref{star2}), then they are independent, because $\sum_{l=1}^n c_l\vPsi_l^*(\lambda)=0$ holds when  $\lambda\to \lambda_i$, $\forall i=1,...,n$, only if  $c_i=0$.  Moreover, $\sum_{l=1}^n c_l\vPsi_l^*(\lambda)$ is analytic at all $\lambda_j\neq \lambda_k$ for a chosen $k$, only if $c_j=0$ for all $j\neq k$. Therefore, there is only one solution analytic at all $\lambda_j\neq \lambda_k$, which is a multiple of $\vPsi_k^*(\lambda)$. In other words, the space $\cap_{j\neq k} V_j$ of solutions regular at all $\lambda_j\neq \lambda_k$ is one dimensional, thus there are no polynomial solutions and $A_1$ cannot have a negative integer eigenvalue.

To prove the last assertion, write $\Psi(\lambda)=\Psi^*(\lambda)\tilde{C}$. Thus $\vPsi_k(\lambda)=\vPsi_k^*(\lambda)\tilde{C}_{kk}+\sum_{j\neq k}\vPsi_j^*(\lambda)\tilde{C}_{jk}$. Now, consider the behaviours at $\lambda = \lambda_k$ of l.h.s. and r.h.s. By (\ref{psik0}), (\ref{star1})  and  (\ref{star2}) we conclude that $\tilde{C}_{kk}=1$ in case $\lp_k\not \in\mathbb{Z}$, and $\tilde{C}_{kk}=0$ in cases $\lp_k\in \mathbb{Z}$. Thus $\tilde{C}_{kk}=c_{kk}$. Now consider the behaviours at $\lambda = \lambda_j$, $j\neq k$,  of l.h.s. and r.h.s. By
 (\ref{psik}), (\ref{star1}) and (\ref{star2}) we obtain $\tilde{C}_{jk}=c_{jk}$.  $\Box$

\vskip 0.5 cm 
\noindent
{\bf Proof of Proposition \ref{etto}:} i) is immediate, because if $\lambda_k$ is accessible and $\lambda\in {\cal P}_\eta\cap{\cal P}_{\tilde{\eta}}$, then  $\vPsi_k^{(k)}(\lambda,\eta)$ and $\vPsi_k^{(k)}(\lambda,\tilde{\eta})$ are the same branch. The same holds for the solutions  $\vPsi_k$. 

ii) and iii) are proved noticing that $\lambda_1,...,\lambda_n$ are all accessible from ${\cal P}_\eta\cap{\cal P}_{\tilde{\eta}}$, therefore i) holds for any $k$. Thus $C(\eta)=C(\tilde{\eta})$.  From Proposition \ref{prop1} we see that $\Psi^*$ is uniquely defined starting from the  $\vPsi_k$'s and $C$. 
$\Box$

\vskip 0.5 cm 

\noindent
{\bf Proof of Proposition \ref{propoW}:} The proof follows that of Proposition 4 of \cite{BJL}, here  generalized  to the general case when $\lp_1,...,\lp_k$ are any complex numbers. 

 Choose a  $\lambda_k$. In the plane ${\cal P}_\eta$ add the cut $\arg(\lambda-\lambda_k)=\tilde{\eta}$. Pet $\Pi_{\eta,\tilde{\eta}}^k$ be the connected region given by the reference points relative to these cuts, namely all points satisfying both conditions $\eta-2\pi<\arg(\lambda-\lambda_j) <\eta$, $\forall j=1,..,n$,  and $\tilde{\eta}-2\pi<\arg(\lambda-\lambda_k)<\tilde{\eta}$. See figure \ref{fig2}. We have then a relation between fundamental matrices 
$$
 \vPsi_k^*(\lambda,\tilde{\eta})=\sum_{j=1}^n \vPsi_j^*(\lambda,{\eta})~W_{jk},~~~~~\lambda\in \Pi_{\eta,\tilde{\eta}}^k.
$$

\begin{figure}
\centerline{\includegraphics[width=0.5\textwidth]{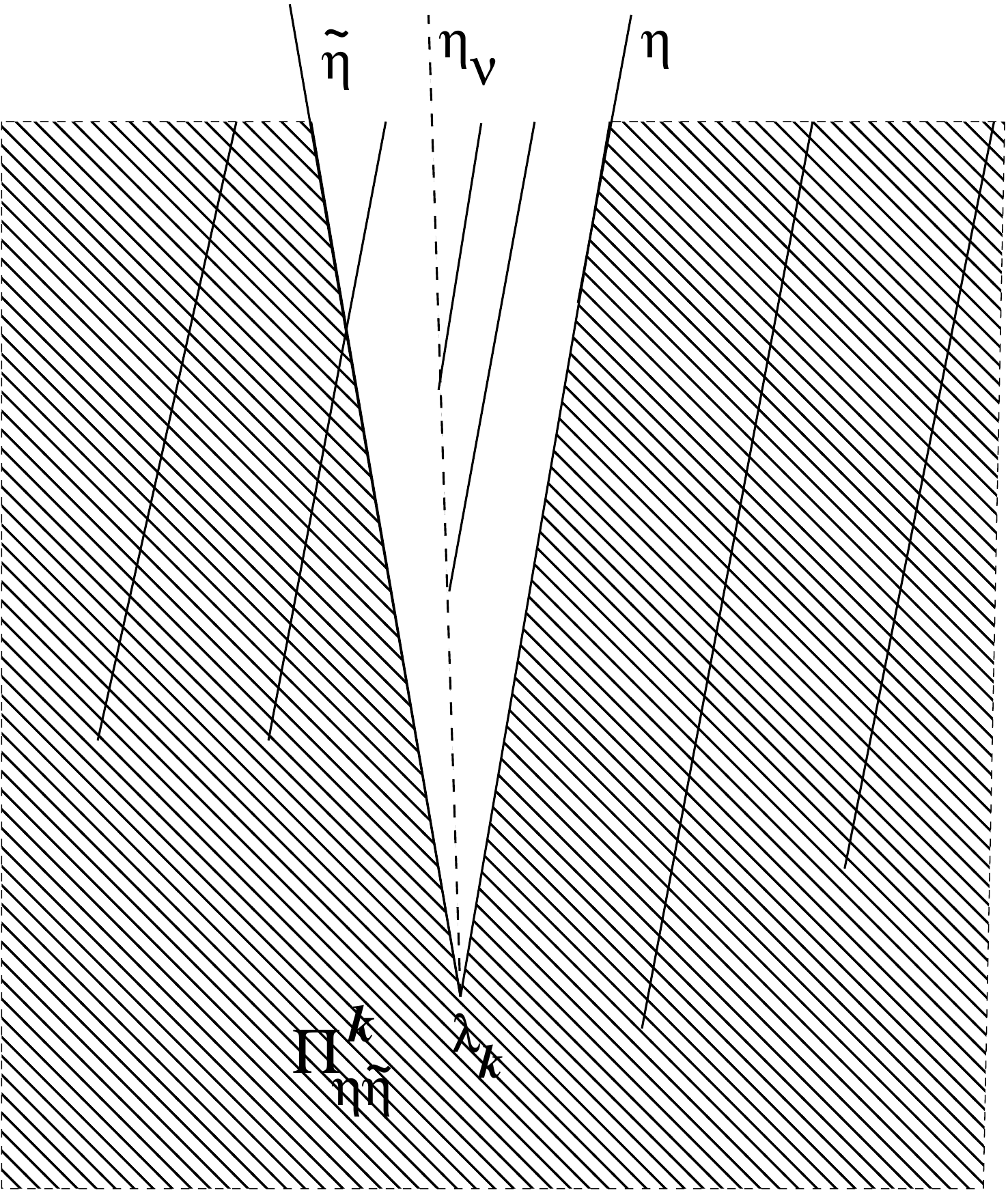}}
\caption{The connected region $\Pi_{\eta,\tilde{\eta}}^k$}
\label{fig2}
\end{figure}

\noindent
-- If $\lambda_j$ is accessible and $j\neq k$, then $W_{jk}=0$, because $\vPsi_j^*(\lambda,{\eta})$ is singular at $\lambda_j$ and $ \vPsi_k^*(\lambda,\tilde{\eta})$ is not. Both are the same branch.

\noindent
--  At $\lambda_k$, both $ \vPsi_k^*(\lambda,\tilde{\eta})$ and $ \vPsi_k^*(\lambda,\tilde{\eta})$ have the same singular behavior with the same branch of logarithm (Proposition \ref{prop1}), therefore $ \vPsi_k^*(\lambda,\tilde{\eta})= \vPsi_k^*(\lambda,\tilde{\eta})+$reg$(\lambda-\lambda_k)$. This implies $W_{kk}=1$. Thus 
\be
\label{tmp}
 \vPsi_k^*(\lambda,\tilde{\eta})=\vPsi_k^*(\lambda,{\eta})+\sum_{j ~not~ acc.} \vPsi_j^*(\lambda,{\eta})~W_{jk},~~~~~\lambda\in \Pi_{\eta,\tilde{\eta}}^k.
\ee
where the sum is on indexes of non accessible points. 

\noindent
-- If $\lambda_j$ is not accessible, in order to access it and evaluate the behaviour of l.h.s. and r.h.s. of (\ref{tmp}), we need to analytically continue the r.h.s beyond the cut $\arg(\lambda-\lambda_k)=\eta$, crossing it in anticlockwise direction. The l.h.s is already defined beyond this cut. Also, all $\vPsi_j^*(\lambda,{\eta})$ for $\lambda_j$ not accessible are  already defined beyond this cut. 
Therefore, after crossing the cut, only $\vPsi_k^*(\lambda,{\eta})$ is analytically continued by   $M_k^*$, according to Proposition \ref{promod}, and (\ref{tmp}) becomes 
$$
\vPsi_k^*(\lambda,\tilde{\eta})=m_{kk}\vPsi_k^*(\lambda,{\eta})+\alpha_k \sum_{j\neq k} c_{jk}\vPsi_j^*(\lambda,{\eta})+\sum_{j ~not~ acc.} \vPsi_j^*(\lambda,{\eta})~W_{jk}
$$
Since the l.h.s is not singular at inaccessible points, we conclude that $W_{jk}= -\alpha_k c_{jk}$, for any inaccessible $\lambda_j$. Finally, we note that the inaccessible points are those such that $\arg(\lambda_j-\lambda_k)=\eta_\nu$. We also see that $j\succ k$ w.r.t $\eta$. 

The same argument is repeated for any $k$. 

As for $W_\nu^{-1}$, the argument does not change, but this time $C$ and $M_k^*$ are referred to $\tilde{\eta}$.
$\Box$

\vskip 0.5 cm 
\noindent
{\bf Proof of Lemma \ref{pippo}:} 
 If  we write $\lambda_\eta$ if $\lambda\in{\cal P}_\eta$. Then 
$$
\lambda_\eta-\lambda_k=(\lambda_{\eta-2\pi}-\lambda_k)e^{2\pi i}
$$
 It follows from the definition of $\Psi$ that, whatever the values $\lp_1,...,\lp_n$ are,  we have
$$\Psi(\lambda,\eta-2\pi)=\Psi(\lambda,\eta)e^{2\pi i \Lambda^\prime}.
$$  
on the universal covering of $\mathbb{C}\backslash\{\lambda_1,...,\lambda_n\}$. 
 From the above and the connection relations we have
$$
 \vPsi_k(\lambda,\eta)=\vPsi_j^{(sing)}(\lambda,\eta)~c_{jk}(\eta)+\hbox{reg}(\lambda-\lambda_j),
$$
and
$$
 \vPsi_k(\lambda,\eta)=e^{-2\pi i\lp_k} \vPsi_k(\lambda,\eta-2\pi)=
$$
$$
=e^{-2\pi i\lp_k} ~\vPsi_j^{(sing)}(\lambda,\eta-2\pi)~c_{jk}(\eta-2\pi)+\hbox{reg}(\lambda-\lambda_j),
$$
Now we use the fact that $\vPsi_j^{(sing)}(\lambda,\eta-2\pi)=\vPsi_j(\lambda,\eta-2\pi)$ when $\lp_j\not\in \mathbb{Z}$, and when $ \lp_j\in \{-1\}\cup \mathbb{N}$ we have
$$
\vPsi_j^{(j)}(\lambda,\eta-2\pi)=\vPsi_j(\lambda,\eta-2\pi)\left[\ln(\lambda_{\eta-2\pi}-\lambda_j)+{P^{(j)}_{N_j}(\lambda_{\eta-2\pi})\over (\lambda-\lambda_j)^{N_j+1}}\right]+\hbox{reg}(\lambda_{\eta-2\pi}
-\lambda_j)
$$
$$=
e^{2\pi i \lp_j}\vPsi_j(\lambda,\eta)\left\{ \left[ \ln(\lambda_\eta-\lambda_j)-2\pi i\right] +{P^{(j)}_{N_j}(\lambda_{\eta})\over (\lambda-\lambda_j)^{N_j+1}}  \right\}+ \hbox{reg}(\lambda_\eta-\lambda_j)
$$
$$
=e^{2\pi i \lp_j}\vPsi_j^{(sing)}(\lambda,\eta)+ \hbox{reg}(\lambda-\lambda_j).
$$
The last step is due to the fact that $\vPsi_j(\lambda,\eta)$ is analytic at $\lambda_j$ (note that $e^{2\pi i \lp_j}=1$) and is absorbed into $ \hbox{reg}(\lambda-\lambda_j)$. 
Therefore, in all cases we have (in full generality, taking $P^{(j)}_{N_j}\equiv 0$ when $\lp_j\in\mathbb{Z}_{-}$):
$$
\vPsi_j^{(sing)}(\lambda,\eta)~c_{jk}(\eta)+\hbox{reg}(\lambda-\lambda_j)=e^{2\pi i(\lp_j-\lp_k)}~\vPsi_j^{(sing)}(\lambda,\eta)~c_{jk}(\eta-2\pi)+\hbox{reg}(\lambda-\lambda_j).
$$
Finally, since $\vPsi_j^{(sing)}(\lambda,\eta)$ is singular at $\lambda_j$, the statement of the lemma follows.
$\Box$

\vskip 0.5 cm 
\noindent
{\bf Proof of Proposition \ref{cpmPro}:} Compared to \cite{BJL}, we only need to consider more general monodromy matrices. Relation (\ref{cpiu}) is 
$$
 \vPsi_k^*(\lambda,\eta-\pi)=\sum_{j=1}^n \vPsi_j^*(\lambda,\eta) C_{jk}^+,~~~\lambda\in{\Pi}_{\eta,\eta-\pi}^k
$$
See proof of Proposition \ref{propoW} for the definition of ${\Pi}_{\eta,\eta-\pi}^k$.  See also figure \ref{fig3}.

\begin{figure}
\centerline{\includegraphics[width=0.5\textwidth]{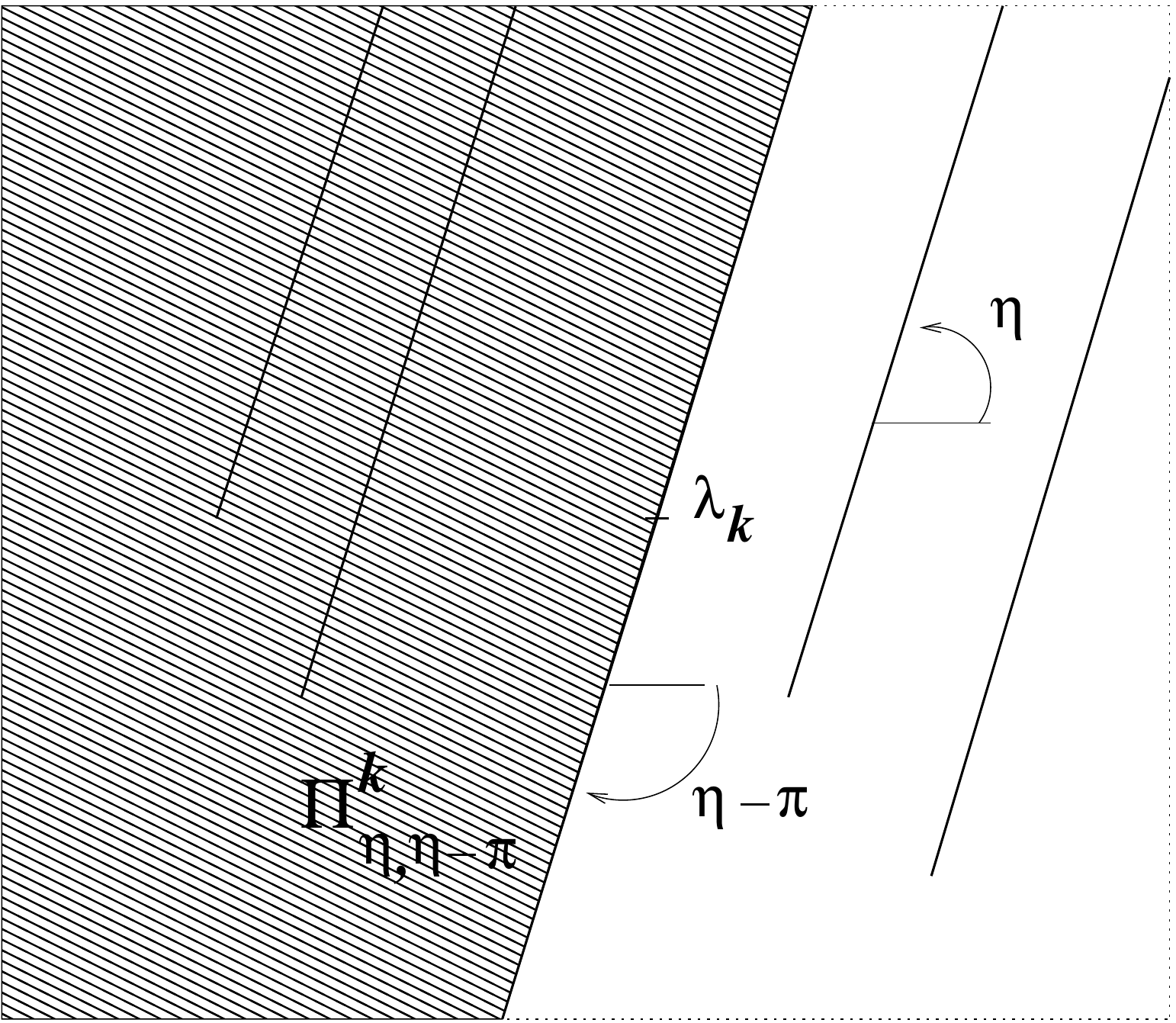}}
\caption{The connected region $\Pi_{\eta,\eta-\pi}^k$, branch cuts at angle $\eta$ and one branch cut at angle $\eta-\pi$ at $\lambda_k$}
\label{fig3}
\end{figure}

- If $\lambda_j\neq \lambda_k$ is accessible,  it means that $j\succ k$. At $\lambda_j$,  $\vPsi_j^*(\lambda,\eta)$ is singular, while $\vPsi_k^*(\lambda,\eta-\pi)$ is not. Both have the same branch of $\ln(\lambda-\lambda_j)$, so we must conclude that $ C_{jk}^+=0$. 

- At $\lambda_k$,  $\vPsi_k^*(\lambda,\eta)$ and  $\vPsi_k^*(\lambda,\eta-\pi)$ are both singular with the same branch of $\ln(\lambda-\lambda_k)$, so we must conclude that $ C_{kk}^+=1$. 

- If $\lambda_j\neq \lambda_k$ is not accessible, it means that $j\prec k$. The l.h.s. is define also at $\lambda_j$, while the r.h.s. needs analytic continuation in order to reach $\lambda_j$. The continuation is given by $ \Psi^*(\lambda,\eta)\mapsto \Psi^*(\lambda,\eta)(M_k^*(\eta))^{-1}$.  therefore 
$$
 \vPsi_k^*(\lambda,\eta-\pi)= \vPsi_k^*(\lambda,\eta)e^{2\pi i \lp_k}+\sum_{l\neq k}\beta_k c_{lk}(\eta)~\vPsi_l^*(\lambda,\eta)
+
\sum_{j~ not ~access.} \vPsi_j^*(\lambda,\eta)C_{jk}^+.
$$
Since $ \vPsi_k^*(\lambda,\eta-\pi)$ is not singular at a non accessible $\lambda_j$ , while $\vPsi_l^*(\lambda,\eta)$ is, we conclude that $C_{jk}^++\beta_kc_{jk}(\eta)=0$. 

The same reasoning is repeated for (\ref{cmeno}). This time the monodromy which gives analytic continuation is $M_k^*(\eta-2\pi)$. Thus, we obtain that $C_{jk}^-=-\alpha_k c_{jk}(\eta-2\pi)$ when $j\succ k$ (non accessible case), $C_{jk}^-=0$ for $j\prec  k$ (accessible case) and $C_{kk}^-=1$ . Then, we use Lemma \ref{pippo}.
$\Box$


\end{document}